\documentclass[11pt]{amsart}
\usepackage{amscd}
\usepackage{amssymb}
\usepackage{amsfonts}
\usepackage{amsmath}

\allowdisplaybreaks[1]

\newtheorem{thm}{Th\'eor\`eme}[section]
\newtheorem{rema}[thm]{Remarque}

\newtheorem{lem}[thm]{Lemme}
\newtheorem{prop}[thm]{Proposition}

\newtheorem{defn}[thm]{D\'efinition}
\theoremstyle{definition}

\numberwithin{equation}{section}

\newcommand{\n}{\noindent}

\newcommand{\resumename}{R\'esum\'e}
\newenvironment{resume}{\narrower\footnotesize\bf
\noindent\resumename.\quad\footnotesize\rm}{\par\bigskip}

\font\hb=cmbx12
\newcommand\pt{\hbox{\hb .}}

\keywords{Alg\`ebres homotopiques, cog\`ebres, alg\`ebres diff\'erentielles gradu\'ees, cohomologie}

\begin{document}

\title[Alg\`ebres de Gerstenhaber \`a homotopie pr\`es]{ Alg\`ebres et cog\`ebres de Gerstenhaber et cohomologie de Chevalley-Harrison}
\author[W. Aloulou, D. Arnal et R. Chatbouri]{Walid Aloulou, Didier Arnal et
Ridha Chatbouri}
\address{
Institut de Math\'ematiques de Bourgogne\\
UMR CNRS 5584\\
Universit\'e de Bourgogne\\
U.F.R. Sciences et Techniques
B.P. 47870\\
F-21078 Dijon Cedex\\France} \email{Didier.Arnal@u-bourgogne.fr}
\address{
D\'epartement de Math\'ematiques\\
Unit\'e de Recherche Physique Math\'ematique\\
Fa\-cult\'e des Sciences de Monastir\\
Avenue de l'environnement\\
5019 Monastir\\
Tuni\-sie} \email{Walid.Aloulou@ipeim.rnu.tn}
           \email{Ridha.Chatbouri@ipeim.rnu.tn}

\thanks{
Ce travail a \'et\'e effectu\'e dans le cadre de l'accord CMCU 06 S 1502.
W. Aloulou et R. Chatbouri remercient l'Universit\'e de Bourgogne pour l'accueil
dont ils ont b\'en\'efici\'e au cours de leurs s\'ejours, D. Arnal remercie la
Facult\'e des Sciences de Monastir pour l'accueil dont il a b\'en\'efici\'e au
cours de ses s\'ejours.}


\begin{abstract}
The fundamental example of Gerstenhaber algebra is the space $T_{poly}({\mathbb R}^d)$ of polyvector fields on $\mathbb{R}^d$, equipped with the wedge product and the Schouten bracket.

In this paper, we explicitely describe what is the enveloping $G_\infty$ algebra of a Gerstenhaber algebra $\mathcal{G}$. This structure gives us a definition of the Chevalley-Harrison cohomology operator for $\mathcal{G}$.

We finally show the nontriviality of a Chevalley-Harrison cohomology group for a natural Gerstenhaber subalgebra in $T_{poly}({\mathbb R}^d)$.
\end{abstract}


\

\maketitle
\begin{resume}
Un prototype des alg\`ebres de Gerstenhaber est l'espace $T_{poly}({\mathbb R}^d)$ des champs de tenseurs sur ${\mathbb R}^d$
muni du produit ext\'erieur et du crochet de Schouten.

Dans cet article, on d\'ecrit explicitement la structure de la $G_\infty$ alg\`ebre enveloppante d'une alg\`ebre de Gerstenhaber. Cette structure permet de d\'efinir une cohomologie de Chevalley-Harrison sur cette alg\`ebre.

On montre que cette cohomologie \`a valeur dans $\mathbb{R}$ n'est pas triviale dans le cas de la sous alg\`ebre de Gerstenhaber des tenseurs homog\`enes $T_{poly}^{hom}({\mathbb R}^d)$.
\end{resume}

\


\

\section{Introduction et motivation}\label{sec0}

\

L'espace $T_{poly}(\mathbb{R}^d)$ des champs de tenseurs antisym\'etriques est une alg\`ebre de Lie gradu\'ee pour le crochet de Schouten. Afin d'\'etudier la cohomologie de Chevalley de cette alg\`ebre pour la repr\'esentation adjointe, on peut se resteindre comme dans \cite{[AAC1]} \`a des cocha\^ines tr\`es simples: les cocha\^ines lin\'eaires ou vectorielles d\'efinies sur les champs de vecteurs (respectivement les tenseurs lin\'eaires). Dans les deux cas, la cohomologie est donn\'ee par les m\^emes cocha\^ines caract\'eris\'ees par leur valeur sur les champs de vecteurs lin\'eaires
$$
\alpha=\sum\alpha^i(x)\partial_i,\qquad\alpha^i(x)\text{ lin\'eaire}.
$$

\

Mais $\mathcal{G}=T_{poly}(\mathbb{R}^d)$ poss\`ede une structure plus riche: c'est une alg\`ebre de Gerstenhaber pour le produit ext\'erieur et le crochet de Schouten. Gr\^ace \`a cette structure, on peut d\'efinir une cohomologie de Chevalley-Harrison dont les cocha\^ines sont les applications lin\'eaires de $S^+\big((\underline{\bigotimes}^+\mathcal{G}[1])[1]\big)$ dans $\mathcal{G}[1]$ (la d\'efinition de cette bicog\`ebre est donn\'ee dans les sections 6 et 7). Cette derni\`ere cohomologie est triviale (voir \cite{[GH]} ou \cite{[T]}).

\

Les champs de vecteurs ou les tenseurs lin\'eaires ne sont que des sous alg\`ebres de Lie de $T_{poly}(\mathbb{R}^d)$. Une sous alg\`ebre de Gerstenhaber simple (et int\'eressante) de $T_{poly}(\mathbb{R}^d)$ est form\'ee par l'espace not\'e $T_{poly}^{hom}(\mathbb{R}^d)$  des tenseurs homog\`enes, c'est \`a dire des tenseurs
$$
\alpha=\sum_{i_1<\dots<i_k}\alpha^{i_1\dots i_k}(x)\partial_{i_1}\wedge\dots\wedge\partial_{i_k},\qquad\text{ avec }\alpha^{i_1\dots i_k}(x)\text{ polyn\^ome homog\`ene de degr\'e }k.
$$

En particulier, cette sous alg\`ebre est de dimension finie et elle contient toutes les structures de Poisson quadratiques. Les cocycles de Chevalley fondamentaux sur les champs de vecteurs et les tenseurs lin\'eaires d\'ecrits dans \cite{[AAC1]} et \cite{[AAC2]} ne sont pas nuls sur $T_{poly}^{hom}(\mathbb{R}^d)$.

\

Le but de cet article est d'\'etudier plus en d\'etail cette situation. Tout d'abord, nous reprenons explicitement et compl\`etement la construction de la structure de cog\`ebre de Gerstenhaber et de $G_\infty$-alg\`ebre induites par celle d'alg\`ebre de Gerstenhaber sur l'espace $S^+\big((\underline{\bigotimes}^+\mathcal{G}[1])[1]\big)$, ceci nous permet de pr\'eciser \`a chaque \'etape les signes apparaissant dans les prolongements des op\'erateurs d\'efinis sur $\mathcal{G}$. Plus pr\'ecisement, une alg\`ebre de Gerstenhaber est un espace vectoriel gradu\'e $\mathcal{G}$ muni de deux op\'erations $\wedge$ et $[~,~]$ de degr\'es respectifs $0$ et $-1$ et respectivement commutatif et antisym\'etrique gradu\'es. Malheureusement les axiomes usuels de cette structure ne satisfont pas la r\`egle des signes de Koszul sur $\mathcal{G}$. On proc\`ede donc \`a un premier d\'ecalage en consid\'erant l'espace $\mathcal{G}[1]$. On obtient deux op\'erations $\mu_2$ et $[~,~]$, dont la sym\'etrie est l'oppos\'ee de celle de $\wedge$ et $[~,~]$ et dont les axiomes v\'erifient bien la r\`egle des signes de Koszul.

On effectue, alors, le prolongement en $\mu$ et $[~,~]$ de ces deux structures sur le quotient $\underline{\bigotimes}^+(\mathcal{G}[1])$ de $T(\mathcal{G}[1])$ par l'espace engendr\'e par les images de toutes les applications battements.

Sur $\underline{\bigotimes}^+(\mathcal{G}[1])$, on a un cocrochet naturel $\delta$. Le produit $\mu$ est une cod\'erivation de $\delta$ telle que $\mu\circ \mu=0$.

Revenant au crochet $[~,~]$, on l'\'etend \`a $\underline{\bigotimes}^+(\mathcal{G}[1])$ de fa\c{c}on \`a en faire une alg\`ebre de Lie. La construction classique consiste \`a consid\'erer l'espace $S^+\big((\underline{\bigotimes}^+\mathcal{G}[1])[1]\big)$ et \`a le munir d'un coproduit $\Delta$ et d'une cod\'erivation $\ell$ telle que $\ell\circ\ell=0$. Il reste \`a d\'ecaler $\delta$ et $\mu$ respectivement en $\kappa$ et $m$ pour les \'etendre aussi \`a $S^+\big((\underline{\bigotimes}^+\mathcal{G}[1])[1]\big)$. Cependant, $m$ et $\ell$ sont des cod\'erivations de $\Delta$ et $\kappa$ telles que $(m+\ell)\circ(m+\ell)=0$.
 Cette construction est explicitement d\'ecrite dans les sections 5 et 6.

\

Finalement, on montre que le $3$-cocycle de Chevalley fondamental $C$ d\'efini sur les champs de vecteurs \`a valeurs dans $C^\infty(\mathbb{R}^d)$ donne sur $T_{poly}^{hom}(\mathbb{R}^d)$ une cocha\^ine \`a valeurs dans $\mathbb{R}$ not\'ee $f$. Nous montrons que $f$ est un cocycle de Chevalley-Harrison non trivial.

\


\section{$A_\infty$ alg\`ebres et cohomologie de Hochschild}\label{sec1}\

Soit $A$ une alg\`ebre associative $|~|$-gradu\'ee, son produit $A\otimes A\longrightarrow A$, $\alpha\otimes\beta\longmapsto\alpha.\beta$ est associatif de degr\'e $0$. On consid\`ere l'espace $A[1]$ et la graduation $deg(\alpha)=|\alpha|-1$ qu'on note simplement par la lettre $\alpha$. On construit un nouveau produit $m_2$ sur $A[1]$ d\'efini par $m_2(\alpha\otimes\beta)=(-1)^\alpha\alpha.\beta$. Alors, $m_2$ devient un produit antiassocitif de degr\'e $1$ sur $A[1]$: $$m_2(m_2(\alpha,\beta),\gamma)=-(-1)^\alpha m_2(\alpha,m_2(\beta,\gamma)).$$

On consid\`ere, maintenant, l'alg\`ebre tensorielle de $A[1]$ sans unit\'e: $T^+(A[1])=\displaystyle\bigoplus_{p\geq1}$$(\bigotimes^pA[1])$.

Cette alg\`ebre munie du coproduit de d\'econcat\'enation
\begin{align*}\triangle(\alpha_1\otimes\dots\otimes\alpha_p)=\displaystyle\sum_{k=1}^{p-1}\big(\alpha_1\otimes\dots\otimes\alpha_k\big)\bigotimes\big(
\alpha_{k+1}\otimes\dots\otimes\alpha_p\big)\end{align*}
est une cog\`ebre coassociative.
En effet, on a
\begin{align*}
(id\otimes\bigtriangleup)\circ\bigtriangleup(\alpha_1\otimes\dots\otimes\alpha_p)&=(id\otimes\bigtriangleup)\Big(\sum_{k=1}^{p-1}\big
(\alpha_1\otimes\dots\otimes\alpha_k\big)\bigotimes\big(
\alpha_{k+1}\otimes\dots\otimes\alpha_p\big)\Big)\cr&=\sum_{1\leq k<j\leq p-1}(\alpha_1\otimes\dots\otimes\alpha_k)\bigotimes(\alpha_{k+1}\otimes\dots\otimes\alpha_j)\bigotimes(\alpha_{j+1}\otimes\dots\otimes\alpha_p)
\end{align*}et
\begin{align*}
(\bigtriangleup\otimes id)\circ\bigtriangleup(\alpha_1\otimes\dots\otimes\alpha_p)&=(\bigtriangleup\otimes id)\Big(\sum_{j=1}^{p-1}\big
(\alpha_1\otimes\dots\otimes\alpha_j\big)\bigotimes\big(
\alpha_{j+1}\otimes\dots\otimes\alpha_p\big)\Big)\cr&=\sum_{1\leq k<j\leq p-1}(\alpha_1\otimes\dots\otimes\alpha_k)\bigotimes(\alpha_{k+1}\otimes\dots\otimes\alpha_j)\bigotimes(\alpha_{j+1}\otimes\dots\otimes\alpha_p).
\end{align*}
Donc, $(id\otimes\bigtriangleup)\circ\bigtriangleup=(\bigtriangleup\otimes id)\circ\bigtriangleup$ et la cog\`ebre $\big(T^+(A[1]),\triangle\big)$ est alors coassociative. Cette cog\`ebre est colibre, ce qui permet de prolonger toute application lin\'eaire $Q_k:\bigotimes^pA[1]\longrightarrow A[1]$ en une cod\'erivation de fa\c{c}on unique. En particulier, on prolonge le produit $m_2$ \`a $\big(T^+(A[1]),\triangle\big)$ comme une cod\'erivation $m$ de cette cog\`ebre ($(m\otimes id+id\otimes m)\circ\bigtriangleup=\bigtriangleup\circ m$) en posant:
$$m(\alpha_1\otimes\dots\otimes\alpha_p)=\sum_{j=1}^{p-1}(-1)^{\sum_{i<j} \alpha_i}\alpha_1\otimes\dots\otimes\alpha_{j-1}\otimes m_2(\alpha_j,\alpha_{j+1})\otimes\alpha_{j+2}\otimes\dots\otimes\alpha_p.$$
 La cod\'erivation $m$ est de degr\'e $1$ dans $T^+(A[1])$, elle v\'erifie $m\circ m=0$.
\vskip0.3cm

En effet, d'une part on a
\begin{align*}
(m\otimes id+id\otimes& m)\circ\bigtriangleup(\alpha_1\otimes\dots\otimes\alpha_p)=(m\otimes id+id\otimes m)\Big(\sum_{k=1}^{p-1}\alpha_1\otimes\dots\otimes\alpha_k\bigotimes
\alpha_{k+1}\otimes\dots\otimes\alpha_p\Big)\cr&=\sum_{k=1}^{p-1}\Big(\sum_{j=1}^{k-1}(-1)^{\sum_{i<j} \alpha_i}\big(\alpha_1\otimes\dots\otimes m_2(\alpha_j,\alpha_{j+1})\otimes\dots\otimes\alpha_k\big)\bigotimes\big(\alpha_{k+1}\otimes\dots\otimes\alpha_p\big)\cr&+\sum_{j=k+1}^{p-1}
(-1)^{\sum_{i<j}\alpha_i}\big(\alpha_1\otimes\dots\otimes\alpha_k\big)\bigotimes\big(\alpha_{k+1}\otimes \dots\otimes m_2(\alpha_j,\alpha_{j+1})\otimes\dots\otimes\alpha_p\big)\Big)\cr&=\sum_{1\leq j<k\leq p-1}(-1)^{\sum_{i<j} \alpha_i}\big(\alpha_1\otimes\dots\otimes m_2(\alpha_j,\alpha_{j+1})\otimes\dots\otimes\alpha_k\big)\bigotimes\big(\alpha_{k+1}\otimes\dots\otimes\alpha_p\big)\cr&+\sum_{1\leq k<j\leq p-1}
(-1)^{\sum_{i<j}\alpha_i}\big(\alpha_1\otimes\dots\otimes\alpha_k\big)\bigotimes\big(\alpha_{k+1}\otimes \dots\otimes m_2(\alpha_j,\alpha_{j+1})\otimes\dots\otimes\alpha_p\big).
\end{align*}

D'autre part, on a
\begin{align*}
\bigtriangleup\circ m&(\alpha_1\otimes\dots\otimes\alpha_p)=\bigtriangleup\Big(\sum_{j=1}^{p-1}(-1)^{\sum_{i<j} \alpha_i}\alpha_1\otimes\dots\otimes\alpha_{j-1}\otimes m_2(\alpha_j,\alpha_{j+1})\otimes\alpha_{j+2}\otimes\dots\otimes\alpha_p\Big)\cr&=\sum_{1\leq j<k\leq p-1}(-1)^{\sum_{i<j} \alpha_i}\big(\alpha_1\otimes\dots\otimes m_2(\alpha_j,\alpha_{j+1})\otimes\dots\otimes\alpha_k\big)\bigotimes\big(\alpha_{k+1}\otimes\dots\otimes\alpha_p\big)\cr&+\sum_{1\leq k<j\leq p-1}
(-1)^{\sum_{i<j}\alpha_i}\big(\alpha_1\otimes\dots\otimes\alpha_k\big)\bigotimes\big(\alpha_{k+1}\otimes \dots\otimes m_2(\alpha_j,\alpha_{j+1})\otimes\dots\otimes\alpha_p\big)\cr&=(m\otimes id+id\otimes m)\circ\bigtriangleup(\alpha_1\otimes\dots\otimes\alpha_p).
\end{align*}

$m$ \'etant une cod\'erivation de degr\'e impair, $m\circ m=0$ est aussi une cod\'erivation. Avec les notations pr\'ec\'edentes, $m\circ m$ est l'unique prolongement de $(m\circ m)_3$ et puisque
$$
(m\circ m)_3(\alpha_1\underline\otimes\alpha_2\underline\otimes\alpha_3)=m_2(m_2(\alpha_1\underline\otimes\alpha_2)\underline\otimes\alpha_3+(-1)^{a_1} \alpha_1\underline\otimes m_2(\alpha_2\underline\otimes\alpha_3))=0.
$$
Par unicit\'e de la cod\'erivation qui prolonge les $(m\circ m)_k$, on en d\'eduit que $m\circ m=0$.\\

\begin{defn} {\rm ($A_\infty$ alg\`ebre)}

\

Une $A_\infty$ alg\`ebre est une cog\`ebre codiff\'erentielle gradu\'ee coassociative de la forme $(T^+(A[1]),\triangle,m)$ o\`u $\triangle$ est le coproduit de d\'econcat\'enation et $m$ est une cod\'erivation  de $\triangle$ de degr\'e $1$ et de carr\'e nul.
\end{defn}

Soit $F:(T^+(A[1]),\triangle)\longrightarrow(T^+(B[1]),\triangle')$ un morphisme de cog\`ebres ($(F\otimes F)\circ\bigtriangleup=\bigtriangleup'\circ F$).
On d\'efinit la projection $F_n$ sur $B[1]$ parall\`element \`a $\displaystyle\bigoplus_{n\geq1}$$(\bigotimes^nB[1])$ de la restriction de $F$ \`a $T^n(A[1])$. L'application $F_n:T^n(A[1])\longrightarrow B[1]$ est $n-$lin\'eaire. Si on connait la suite des $(F_n)_n$, on montre qu'on peut reconstruire $F$ de fa\c{c}on unique, plus pr\'ecis\'ement:
$$F(\alpha_1\otimes\dots\otimes\alpha_n)=\sum_{k=1}^{n-1}\sum_{0<r_1<\dots<r_k<n}F_{r_1}(\alpha_1\otimes\dots\otimes\alpha_{r_1})\otimes\dots\otimes
F_{r_{k}}(\alpha_{r_{k-1}+1}\otimes\dots\otimes\alpha_n)
$$

Nous exposons en d\'etail cette preuve dans la section $4$ dans le cas commutatif.

\begin{defn} {\rm (morphisme de $A_\infty$ alg\`ebres)}

\

Un morphisme de $A_\infty$ alg\`ebres $A$ et $B$ est un morphisme de cog\`ebres coassociatives codiff\'erentielles $F:\big(T^+(A[1]),m^A\big)\longrightarrow \big(T^+(B[1]),m^B\big)$ tel que $m^B\circ F=F\circ m^A$.
\end{defn}

L'\'equation de $A_\infty$ morphisme $m^B\circ F=F\circ m^A$ \'ecrite sur les applications $F_n:T^n(A[1])\longrightarrow B[1]$ d\'efinissant $F$ prend la forme suivante: Posons $\alpha_{\{1,\dots,n\}}=\alpha_1\otimes\dots\alpha_n$.

\n D'une part, on a
$$
\begin{aligned}
&m^B\circ F(\alpha_1\otimes\dots\otimes\alpha_n)=\sum_{k,0<r_1<\dots<r_k<n}m^B\Big(F_{r_1}(\alpha_1\otimes\dots\otimes\alpha_{r_1})\otimes\dots\otimes F_{r_k}(\alpha_{r_{k-1}+1}\otimes\dots\otimes\alpha_n)
\Big)\cr&=\sum_{0<j<k}\sum_{0<r_1<\dots<r_k<n}(-1)^{\sum_{i<r_j}\alpha_i}F_{r_1}(\alpha_{\{1,\dots,r_1\}})\otimes\dots\otimes F_{r_j-1}(\alpha_{\{r_{j-2}+1,\dots,r_{j-1}\}})
\otimes\cr&m^B\big(F_{r_j}(\alpha_{\{r_{j-1}+1,\dots,r_{j}\}})\otimes F_{r_{j+1}}(\alpha_{\{r_{j}+1,\dots,r_{j+1}\}})\big)
\otimes F_{r_{j+2}}(\alpha_{\{r_{j+1}+1,\dots,r_{j+2}\}})\otimes\dots\otimes F_{r_{k}}(\alpha_{\{r_{k-1}+1,\dots,r_{k}\}}).
\end{aligned}
$$

\n D'autre part, on a
$$
\begin{aligned}
&F\circ m^A(\alpha_{\{1,\dots,n\}})=F\big(\sum_{j=1}^{n-1}(-1)^{\alpha_{\{1,\dots,j-1\}}}\alpha_{\{1,\dots,j-1\}}\otimes m^A(\alpha_j\otimes\alpha_{j+1})\otimes\alpha_{\{j+2,\dots,n\}}\big)\cr&=\sum_{j=1}^{n-1}(-1)^{\alpha_{\{1,\dots,j-1\}}}\sum_{k, 0<r_1<\dots<r_k<n}\cr&F_{r_1}(\alpha_{\{1,\dots,r_1\}})\otimes\dots\otimes F_{r_t}\Big(\alpha_{\{r_{t-1}+1,\dots,j-1\}}\otimes m^A
(\alpha_j\otimes\alpha_{j+1})\otimes\alpha_{\{j+2,\dots,r_t\}}\Big)\otimes\dots\otimes F_{r_k}(
\alpha_{\{r_{k-1}+1,\dots,r_k\}}).
\end{aligned}
$$

\

\n Il n'y a pas de $F_n$ dans l'\'equation $(m^B\circ F-F\circ m^A)(\alpha_{\{1,\dots,n\}})=0$, cherchons les termes o\`u $F_{n-1}$ appara\^it: ce sont
$$
\begin{aligned}
&m^B\big(F_{n-1}(\alpha_{\{1,\dots,n-1\}})\otimes F_1(\alpha_n)\big)+m^B\big(F_1(\alpha_1)\otimes F_{n-1}(\alpha_{\{2,\dots,n\}})\big)\cr&
-\sum_{j=1}^{n-1}(-1)^{\alpha_{\{1,\dots,j-1\}}}F_{n-1}\big(\alpha_{\{1,\dots,j-1\}}\otimes m^A
(\alpha_j\otimes\alpha_{j+1})\otimes\alpha_{\{j+2,\dots,n\}}\big)\cr&=(-1)^{\alpha_{\{1,\dots,n-1\}}+1}F_{n-1}(\alpha_{\{1,\dots,n-1\}}). F_1(\alpha_n)+(-1)^{\alpha_1+1}F_1(\alpha_1). F_{n-1}(\alpha_{\{2,\dots,n\}})\cr&-\sum_{j=1}^{n-1}(-1)^{\alpha_{\{1,\dots,j\}+1}}F_{n-1}\big(\alpha_{\{1,\dots,j-1\}}\otimes
(\alpha_j.\alpha_{j+1})\otimes\alpha_{\{j+2,\dots,n\}}\big)\cr&=(d_{H}F_{n-1})(\alpha_{\{1,\dots,n\}})
\end{aligned}
$$

On retrouve l'op\'erateur de cobord de Hochschild $d_{H}$.\\

\vskip0.15cm

Finalement, si $V$ un $A$ bimodule gradu\'e, on note $V[p]$ l'espace gradu\'e $V$ tel que si $v$ est de degr\'e $i$ dans $V$, il sera de degr\'e $i-p$ dans $V[p]$. L'espace $B=A\oplus\sum_{p>0}V[p]$ muni du produit
$$
m'\big((\alpha+\sum u_p),(\beta+\sum v_q)\big)=\big(m_2(\alpha,\beta)+\sum \alpha v_p+u_p\beta\big)
$$
est une alg\`ebre associative et l'application $f:A\longrightarrow B$, d\'efinie par $f(\alpha)=(\alpha,0)$ est un morphisme d'alg\`ebres.\vskip0.15cm

Un morphisme de cog\`ebres coassociatives diff\'erentielles $F=f+C$ sera appel\'e une $A_\infty$ formalit\'e de module s'il est d\'efini par des $F_n$ homog\`enes de degr\'e 0, de la forme $F_1=f+C_1$, o\`u $C_1$ est lin\'eaire de $A[1]$ dans $V$, et pour $p>1$, $F_p=C_p$, o\`u les $C_p$ sont $p$-lin\'eaires de $\bigotimes^pA[1]$ dans $V$.

\vskip0.15cm

On retrouve la cohomologie de Hochschild des alg\`ebres associatives \`a valeurs dans un module $V$.

Plus pr\'ecis\'ement, cette formalit\'e est dite triviale s'il existe un morphisme $G$ tel que $C=m^B\circ G+G\circ m^A$, $G$ \'etant de degr\'e $-1$ et $G=\sum B_p$ avec $B_p:\underline\bigotimes^p(A[1])\longrightarrow V[p]$.\vskip0.15cm

On retrouve ainsi la cohomologie de Hochschild des alg\`ebres associatives, puisque

\vskip 0.3cm
\begin{prop} {\rm ($A_\infty$ formalit\'es et cohomologie de Hochschild)}

\

Avec les notations pr\'ec\'edentes, $F$ est une $A_\infty$ formalit\'e si et seulement si
$$
d_{H}C_k=0\quad\text{pour tout }k>0.
$$

$F$ est triviale si et seulement si
$$
C_1=0\quad\text{et}\quad C_k=d_{H}B_k\quad\text{pour tout }k>1.
$$
\end{prop}

\


\section{$L_\infty$ alg\`ebres et cohomologie de Chevalley}

\

Soit $\mathfrak g$ une alg\`ebre de Lie sur un corps $\mathbb K$. Soit $V$ un $\mathfrak g$ module. La cohomologie de Chevalley de $\mathfrak g$ \`a valeurs dans $V$ est d\'efinie de la fa\c con suivante:

Une $n$-cocha\^\i ne $C$ est une application $n$-lin\'eaire altern\'ee de $\mathfrak g^n$ dans $V$:
$$
C\in C^n(\mathfrak g,V)=Hom(\wedge^n\mathfrak g,V),
$$
son cobord de Chevalley $d_CC\in C^{n+1}(\mathfrak g,V)$ est la cocha\^\i ne d\'efinie comme:
$$\begin{aligned}
d_CC&(X_0,\dots,X_n)=\\&=\sum_{j=0}^n(-1)^jX_jC(X_0,\dots,\hat{j},\dots,X_n)+\sum_{i<j}(-1)^{i+j}C\left([X_i,X_j],X_0,\dots,\hat{\imath},\dots,\hat{\jmath},\dots,X_n\right).
\end{aligned}
$$
On a $d_C\circ d_C=0$ et le $n^{eme}$ groupe de cohomologie $H^n(\mathfrak g,V)$ est le quotient de l'espace $Z^n(\mathfrak g,V)$ des $n$ cocycles (les cocha\^\i nes $C$ telles que $d_CC=0$) par l'espace $B^n(\mathfrak g,V)$  des $n$ cobords (les cocha\^\i nes $C$ telles qu'il existe $b\in C^{n-1}(\mathfrak g,V)$ tel que $C=d_Cb$).\\

Afin de pr\'esenter cette cohomologie de fa\c con plus alg\'ebrique et intrins\`eque, on regarde $\mathfrak g$ comme une $L_\infty$ alg\`ebre. Cela permettra entre autres de g\'en\'eraliser im\-m\'e\-dia\-tement la cohomologie au cas des alg\`ebres de Lie diff\'erentielles et gradu\'ees.\\

Rappelons d'abord la r\`egle des signes de Koszul: si dans les axiomes d'une structure alg\`ebrique on a une somme de quantit\'es qui sont des compositions ou des produits d'objets $X_1,\dots,X_n$ de degr\'es respectifs $x_1,\dots,x_n$, dans divers ordres, lorsqu'on veut d\'efinir la structure gradu\'ee correspondante, on ajoute devant la quantit\'e compos\'ee des objets dans l'ordre $X_{i_1},\dots,X_{i_n}$ la signature de la permutation gradu\'ee $\left(\begin{smallmatrix}x_1&\dots&x_n\\ x_{i_1}&\dots&x_{i_n}\end{smallmatrix}\right)$, c'est \`a dire le signe
$$
\varepsilon\left(\begin{smallmatrix}x_1&\dots&x_n\\ x_{i_1}&\dots&x_{i_n}\end{smallmatrix}\right)
$$
$\varepsilon$ est un morphisme et sur la transposition $(x_i,x_{i+1})$, on a $\varepsilon(x_i,x_{i+1})=(-1)^{x_ix_{i+1}}$.

Par exemple les axiomes d'une alg\`ebre de Lie sont:
$$
[X,Y]=-[Y,X],\qquad \left[[X,Y],Z\right]+\left[[Y,Z],X\right]+\left[[Z,X],Y\right]=0.
$$
Les axiomes d'une alg\`ebre de Lie gradu\'ee seront donc:
$$\aligned
{[X,Y]}&=-\varepsilon\left(\begin{smallmatrix}xy\\ yx\end{smallmatrix}\right)[Y,X]=(-1)^{xy}[Y,X],\\
0&=\left[[X,Y],Z\right]+\varepsilon\left(\begin{smallmatrix}xyz\\ yzx\end{smallmatrix}\right)\left[[Y,Z],X\right]+\varepsilon\left(\begin{smallmatrix}xyz\\ zxy\end{smallmatrix}\right)\left[[Z,X],Y\right]\\
&=\left[[X,Y],Z\right]+(-1)^{x(y+z)}\left[[Y,Z],X\right]+(-1)^{z(x+y)}\left[[Z,X],Y\right].
\endaligned
$$
Pour une alg\`ebre de Lie gradu\'ee diff\'erentielle, on ajoute la diff\'erentielle qui est une application $d:\mathfrak g\longrightarrow\mathfrak g$ de degr\'e 1 et telle que:
$$
d[X,Y]=[dX,Y]+(-1)^x[X,dY].
$$

La premi\`ere \'etape de notre construction consiste en un d\'ecalage des degr\'es. Les signes apparaissant dans la formule de $d_C$ pour une alg\`ebre de Lie usuelle seront alors directement donn\'es par la r\`egle de Koszul et la g\'en\'eralisation $\ell$ de $d_C$ sera de degr\'e 1. On munit donc les vecteurs $X$ de $\mathfrak g$ du degr\'e $degr\acute{e}(X)=x=-1$. On note $\mathfrak g[1]$ cet espace gradu\'e. Le crochet devient une application gradu\'ee symm\'etrique $\ell:S^2\left(\mathfrak g[1]\right)\longrightarrow\mathfrak g[1]$ homog\`ene de degr\'e 1. De m\^eme l'alg\`ebre $\bigwedge\mathfrak g$ est isomorphe en tant qu'espace vectoriel \`a l'alg\`ebre $S\left(\mathfrak g[1]\right)$. Il n'y a pas d'isomorphisme d'alg\`ebre entre ces deux espaces, il n'y a pas non plus d'isomorphisme lin\'eaire canonique. Nous choisissons l'isomorphisme donn\'e dans \cite{[AAC1]}:
$$
X_{i_1}\wedge\dots\wedge X_{i_n}\longrightarrow(-1)^{\sum_j(n-j)x_{i_j}}X_{i_1}\pt\dots\pt X_{i_n}.
$$
Alors $\ell_2(X,Y)=(-1)^x[X,Y]$ et si $\mathfrak g$ est diff\'erentielle, on posera $\ell_1(X)=dX$.\\

On consid\`ere $S^+\left(\mathfrak g[1]\right)=\sum_{n>0}S^n\left(\mathfrak g[1]\right)$ comme une cog\`ebre pour la comultiplication $\Delta$ d\'eduite de la d\'econcat\'enation de $T^+(\mathfrak g[1])$, d\'efini de la mani\`ere suivante:

soit $I=\{i_1<i_2<\dots<i_k\}$ une partie de $\{1,\dots,n\}$, on note $X_I$ le produit $X_{i_1}\pt\dots\pt X_{i_k}$. $\Delta$ est alors la comultiplication de degr\'e 0 d\'efinie par:
$$
\Delta(X_1\pt\dots\pt X_n)=\sum_{\begin{smallmatrix}I\sqcup J=\{1,\dots n\}\\ |I|>0,|J|>0\end{smallmatrix}}\varepsilon\left(\begin{smallmatrix}x_{\{1,\dots,n\}}\\ x_I,x_J\end{smallmatrix}\right)X_I\otimes X_J.
$$
Remarquons que lorsque chaque $x_i=-1$, le signe est simplement la signature de la permutation $\left(\begin{smallmatrix}1,\dots,n\\ I,J\end{smallmatrix}\right)$.

La cog\`ebre ainsi obtenue est une cog\`ebre cocommutative et coassociative : on note $\tau$ la volte gradu\'ee
$$
\tau(X\otimes Y)=\varepsilon\left(\begin{smallmatrix}xy\\yx\end{smallmatrix}\right)Y\otimes X.
$$
alors
$$
\begin{aligned}
\tau\circ\Delta&=\Delta\\
(id\otimes\Delta)\circ\Delta&=(\Delta\otimes id)\circ\Delta.
\end{aligned}
$$
En effet, on a:
$$
\begin{aligned}
\tau\circ\Delta(X_{\{1,\dots,n\}})&=\tau\left(\sum_{I,J}\varepsilon\left(\begin{smallmatrix}x_{\{1,\dots,n\}}\\x_I,x_J\end{smallmatrix}\right)X_I\otimes X_J\right)\\
&=\sum_{I,J}\varepsilon\left(\begin{smallmatrix}x_I,x_J\\x_J,x_I\end{smallmatrix}\right)\varepsilon\left(\begin{smallmatrix}x_{\{1,\dots,n\}}\\x_I,x_J\end{smallmatrix}\right)X_J\otimes X_I\\
&=\sum_{I,J}\varepsilon\left(\begin{smallmatrix}x_{\{1,\dots,n\}}\\x_J,x_I\end{smallmatrix}\right)X_J\otimes X_I=\Delta(X_{\{1,\dots,n\}})\\
\end{aligned}
$$
et
$$
\begin{aligned}
(id\otimes\Delta)\circ\Delta(X_{\{1,\dots,n\}})&=\sum_{I,J}(id\otimes\Delta)\varepsilon\left(\begin{smallmatrix}x_{\{1,\dots,n\}}\\x_I,x_J\end{smallmatrix}\right)X_I\otimes X_J\\
&=\sum_{I\sqcup J=\{1\dots n\}}\sum_{K\sqcup L=J}\varepsilon\left(\begin{smallmatrix}x_{\{1,\dots,n\}}\\x_I,x_J\end{smallmatrix}\right) \varepsilon\left(\begin{smallmatrix}x_J\\x_K,x_L\end{smallmatrix}\right)X_I\otimes X_K\otimes X_L\\
&=\sum_{I\sqcup K\sqcup L=\{1\dots n\}}\varepsilon\left(\begin{smallmatrix}x_{\{1,\dots,n\}}\\x_I,x_K,x_L\end{smallmatrix}\right)X_I\otimes X_K\otimes X_L\\
&=(\Delta\otimes id)\circ\Delta(X_{\{1,\dots,n\}}).
\end{aligned}
$$
Toute application lin\'eaire $f:(\mathcal{C},c)\longrightarrow (\mathfrak{g}[1],\Delta)$ o\`u $(\mathcal{C},c)$ est une cog\`ebre cocommutative, coassociative et nilpotente (c'est \`a dire que pour tout $c$,
$$
(\Delta\otimes id^{\otimes n})\circ(\Delta\otimes id^{\otimes n-1})\circ\dots\circ\Delta c=0
$$
pour $n$ assez grand) se prolonge d'une fa\c{c}on unique en un morphisme de cog\`ebre $F:(\mathcal{C},c)\longrightarrow (S^+(\mathfrak{g}[1]),\Delta)$.

On dira que c'est la cog\`ebre cocommutative et coassociative libre (sans co-unit\'e) engendr\'ee par $\mathfrak g[1]$. On peut donc prolonger de fa\c con unique l'application $\ell_1+\ell_2$ en une cod\'erivation de degr\'e 1 de notre cog\`ebre. Ce prolongement est donn\'e par (voir \cite{[AMM]}):
$$
\begin{aligned}
\ell(X_1\pt\dots\pt X_n)&=\sum_j(-1)^{\sum_{i<j}x_i}X_1\pt\dots\pt\ell_1(X_j)\pt\dots\pt X_n\\
&\hskip 2cm+\sum_{i<j}\varepsilon\left(\begin{smallmatrix}x_1\ \dots\ x_n\\ x_i x_j x_1\dots\hat{\imath}\hat{\jmath}\dots x_n\end{smallmatrix}\right)\ell_2(X_i,X_j)\pt X_1\pt\dots\hat\imath\dots\hat\jmath\dots\pt X_n.
\end{aligned}
$$
$\ell$ est une cod\'erivation veut dire que, en tenant compte de la r\`egle des signes de Koszul dans la d\'efinition du produit tensoriel des applications,
$$
(id\otimes\ell+\ell\otimes id)\circ\Delta=\Delta\circ\ell.
$$
Elle est de carr\'e nul $\ell\circ\ell=0$ (voir \cite{[AMM]}, \cite{[K]}).

\begin{defn} {\rm ($L_\infty$ alg\`ebre)}

\

Une $L_\infty$ alg\`ebre est une cog\`ebre diff\'erentielle de la forme $\left(S^+(\mathfrak g[1]),\Delta,\ell\right)$ o\`u $\Delta$ est d\'efini ci-dessus et $\ell$ est une codiff\'erentielle de $\Delta$ est de carr\'e nul.

Si $\left(\mathfrak g,[~,~],d\right)$ est une alg\`ebre de Lie gradu\'ee diff\'erentielle, la $L_\infty$ alg\`ebre $L(\mathfrak g)=\left(S^+(\mathfrak g[1]),\Delta,\ell\right)$ telle que
$$
\ell_1(X)=dX,\quad\ell_2(X\pt Y)=(-1)^x[X,Y],\quad\ell_k=0,\quad k=3,4,\dots
$$
s'appelle la $L_\infty$ alg\`ebre enveloppante de $\mathfrak g$.
\end{defn}

Un morphisme de $L_\infty$ alg\`ebre entre $S^+\left(\mathfrak g[1]\right)$ et $S^+\left(\mathfrak h[1]\right)$ est une application $F$ entre ces espaces qui est un morphisme de cog\`ebres diff\'erentielles. Puisque $S^+\left(\mathfrak h[1]\right)$ est libre, un tel morphisme est caract\'eris\'e par la donn\'ee d'une suite d'aplications:
$$
F_n:S^n\left(\mathfrak g[1]\right)\longrightarrow\mathfrak h[1],
$$
homog\`ene de  degr\'e 0. $F$ est un morphisme de cog\`ebre si et seulement si, pour tout $n$,
$$
F(X_1\pt\dots\pt X_n)=\sum_{j>0}\frac{1}{j!}\sum_{\begin{smallmatrix}I_1\sqcup\cdots\sqcup I_j=\{1,\ldots,n\}\\ I_1\dots I_j\neq\emptyset\end{smallmatrix}}\varepsilon\left(\begin{smallmatrix}x_1\dots x_n\\ x_{I_1}\dots x_{I_j}\end{smallmatrix}\right)F_{|I_1|}(X_{I_1})\pt\dots\pt F_{|I_j|}(X_{I_j}).
$$
Enfin, $F$ est un morphisme de cog\`ebres diff\'erentielles si et seulement si $\ell^\mathfrak{h}\circ F=F\circ\ell^\mathfrak{g}$. Ceci donne une \'equation sur les $F_n$, appel\'ee \'equation de formalit\'e. Si $\ell^{\mathfrak g}$ (resp $\ell^{\mathfrak h}$) est la cod\'erivation caract\'eris\'ee par les applications $\ell^{\mathfrak g}_p: S^p({\mathfrak g}[1])\longrightarrow{\mathfrak g}[1]$ (resp. $\ell^{\mathfrak h}_q:S^q({\mathfrak h}[1])\longrightarrow{\mathfrak h}[1]$) et si $F$ est caract\'eris\'ee par les applications $F_n: S^n({\mathfrak g}[1])\longrightarrow{\mathfrak h}[1]$, cette \'equation s'\'ecrit:
$$
0=\hskip-0.7cm\sum_{\begin{smallmatrix}1\leq p\leq n\\ I_1\sqcup\dots\sqcup I_p=\{1,\dots,n\}\\ 0<|I_1|,\dots,|I_p|<n\end{smallmatrix}}\hskip-0.4cm\varepsilon\left(\begin{smallmatrix}x_{\{1,\dots,n\}}\\ x_{I_1}\dots x_{I_p}\end{smallmatrix}\right) \ell^{\mathfrak h}_p\left(F_{|I_1|}(X_{I_1})\pt\dots\pt F_{|I_p|}(X_{I_p})\right)-\hskip-0.4cm\sum_{\begin{smallmatrix}1\leq p\leq n\\ I\sqcup J=\{1,\dots,n\}\\ |J|=p-1\end{smallmatrix}}\hskip-0.4cm \varepsilon\left(\begin{smallmatrix}x_{\{1,\dots,n\}}\\ x_{I} x_{J}\end{smallmatrix}\right)F_p\left(\ell^{\mathfrak g}_{|I|}(X_I)\pt X_J\right).
$$

\

Supposons maintenant que ${\mathfrak g}$ et ${\mathfrak h}$ soient deux alg\`ebres de Lie gradu\'ees et $\varphi:{\mathfrak g}\longrightarrow{\mathfrak h}$ un morphisme de degr\'e 0 d'alg\`ebres de Lie. Cherchons tous les morphismes de $L_\infty$ alg\`ebres $F:S^+({\mathfrak g}[1])\longrightarrow S^+({\mathfrak h}[1])$ tels que $F_1=\varphi$. Cela revient \`a chercher toutes les suites d'applications $(F_n)$ ($F_n:S^n({\mathfrak g}[1])\longrightarrow {\mathfrak h}[1]$) telles que:
$$
\begin{aligned}
0=&\hskip-0.7cm\sum_{\begin{smallmatrix}I\sqcup J=\{1,\dots,n\}\\ 0<|I|,|J|<n\end{smallmatrix}}\hskip-0.4cm\frac{1}{2}\varepsilon\left(\begin{smallmatrix}x_{\{1,\dots,n\}}\\ x_Ix_J\end{smallmatrix}\right) \ell^{\mathfrak h}_2\left(F_{|I|}(X_I)\pt F_{|J|}(X_J)\right)\\&\hskip 1.5cm-\sum_{0<i<j<n+1}\hskip-0.4cm \varepsilon\left(\begin{smallmatrix}x_{\{1,\dots,n\}}\\ x_ix_jx_1\dots\hat{\imath}\hat{\jmath}\dots x_n\end{smallmatrix}\right)F_{n-1}\left(\ell^{\mathfrak g}_2(X_i\pt X_j)\pt X_1\dots\hat{\imath}\hat{\jmath}\dots\pt X_n\right).
\end{aligned}
\eqno{(n)}
$$
Cette suite d'\'equations peut se r\'esoudre par r\'ecurrence sur $n$. L'\'equation est v\'eriffi\'ee pour $n=2$ puisque $\varphi$ est un morphisme. Si on a r\'esolu les \'equations $(2),\dots,(n)$ qui portent sur $F_1=\varphi,F_2,\dots,F_{n-1}$, l'\'equation $(n)$ devient une \'equation en $F_n$ de la forme:
$$
\begin{aligned}
\hskip-0.7cm-\frac{1}{2}\sum_{\begin{smallmatrix}I\sqcup J=\{1,\dots,n+1\}\\ 1<|I|,|J|<n\end{smallmatrix}}&\hskip-0.4cm\varepsilon\left(\begin{smallmatrix}x_{\{1,\dots,n+1\}}\\ x_Ix_J\end{smallmatrix}\right) \ell^{\mathfrak h}_2\left(F_{|I|}(X_I)\pt F_{|J|}(X_J)\right)\\&=\sum_i\varepsilon\left(\begin{smallmatrix}x_1\dots x_{n+1}\\ x_ix_1\dots\hat\imath\dots x_{n+1}\end{smallmatrix}\right)\ell^{\mathfrak h}_2(\varphi(X_i)\pt F_n(X_1\pt\dots\hat{\imath}\dots\pt X_{n+1}))-\\
&-\sum_{i<j}\varepsilon\left(\begin{smallmatrix}x_1\dots x_{n+1}\\ x_ix_jx_1\dots\hat\imath\hat\jmath\dots x_{n+1}\end{smallmatrix}\right) F_n(\ell^{\mathfrak g}_2(X_i\pt X_j)\pt X_1\dots\hat\imath\dots\hat\jmath\dots\pt X_{n+1}).
\end{aligned}
$$

On peut donc d\'efinir la $L_\infty$ cohomologie.

\vskip 0.3cm
\begin{defn} {\rm ($L_\infty$ cohomologie)}

\

Soit $\mathfrak g$ et $\mathfrak h$ deux alg\`ebres de Lie gradu\'ees et $\varphi$ un homomorphisme d'alg\`ebres de Lie de degr\'e 0 de $\mathfrak g$ dans $\mathfrak h$. On appelle $n$ cocha\^\i ne sur $\mathfrak g$ \`a valeurs dans $\mathfrak h$ une application $F_n$ de degr\'e $f_n$ de $S^n({\mathfrak g}[1])$ dans $\mathfrak h[1]$. L'op\'erateur de cobord $d_L$ associe \`a cette cocha\^\i ne $F_n$ la cocha\^\i ne
$$
\begin{aligned}
d_LF_n(X_1\pt\dots\pt X_{n+1})&=\ell^{\mathfrak h}\circ F_n-(-1)^{f_n}F_n\circ\ell^{\mathfrak g}\\
&=\sum_i\varepsilon\left(\begin{smallmatrix}x_1~~~\dots~~~ x_{n+1}\\ x_ix_1\dots\hat\imath\dots x_{n+1}\end{smallmatrix}\right)\ell^{\mathfrak h}_2(\varphi(X_i)\pt F_n(X_1\pt\dots\hat{\imath}\dots\pt X_{n+1}))-\\
&-(-1)^{f_n}\sum_{i<j}\varepsilon\left(\begin{smallmatrix}x_1~~~\dots~~~ x_{n+1}\\ x_ix_jx_1\dots\hat\imath\hat\jmath\dots x_{n+1}\end{smallmatrix}\right) F_n(\ell^{\mathfrak g}_2(X_i\pt X_j)\pt X_1\dots\hat\imath\dots\hat\jmath\dots\pt X_{n+1}).
\end{aligned}
$$
\end{defn}

\

Si $F_n$ est de degr\'e $f_n$, $d_LF_n$ est de degr\'e $f_n+1$, donc
$$
\begin{aligned}
d_L\circ d_L(F_n)&=\ell^{\mathfrak h}\circ d_LF_n-(-1)^{f_n+1}d_LF_n\circ\ell^{\mathfrak g}\\
&=\ell^{\mathfrak h}\circ\ell^{\mathfrak h}\circ F_n-(-1)^{f_n}\ell^{\mathfrak h}\circ F_n\circ\ell^{\mathfrak g}-(-1)^{f_n+1}\ell^{\mathfrak h}\circ F_n\circ\ell^{\mathfrak g}-F_n\circ\ell^{\mathfrak g}\circ\ell^{\mathfrak g}=0.
\end{aligned}
$$

On retrouve en particulier la cohomologie de Chevalley usuelle. En effet soit $\mathfrak g$ une alg\`ebre de Lie et $V$ un $\mathfrak g$ module. $V$ permet de construire imm\'ediatement une alg\`ebre de Lie gradu\'ee en posant
$$
\mathfrak h=\mathfrak g\oplus\sum_{p=-1}^\infty V[p], \quad [X+\sum_p u_p,Y+\sum_q v_q]=[X,Y]+\sum_p Xv_p-Yu_p
$$
et un morphisme $\varphi:\mathfrak g\longrightarrow\mathfrak h$ d\'efini par $\varphi(X)=X$.

Un morphisme de cog\`ebres diff\'erentielles $F=f+C$ sera appel\'e une formalit\'e de module s'il est d\'efini par des $F_n$ homog\`enes de degr\'e 0, de la forme $F_1=\varphi+C_1$, $C_1$ lin\'eaire de $\mathfrak g$ dans $V$, et pour $p>1$, $F_p=C_p$, $C_p$ $p$ lin\'eaire de $\mathfrak g$ dans $V$.

De m\^eme une formalit\'e de module $F$ est dite triviale s'il y a un morphisme de cog\`ebres $B$ d\'efini par $B_p$ $p$ lin\'eaire de $\mathfrak g$ dans $V$, de degr\'e -1 tel que $F=\varphi+\ell^\mathfrak h\circ B+B\circ\ell^\mathfrak g$.

\vskip 0.3cm
\begin{prop} {\rm (Expression de la cohomologie de Chevalley)}

L'\'equation de formalit\'e de module pour un morphisme $F=\varphi+C$ est
$$
d_CC_n=0\quad\forall n>0.
$$

De plus, $F=\varphi+\ell_\mathfrak h\circ B+B\circ\ell_\mathfrak g$ si et seulement si
$$
C_1=0\quad\text{et}\quad C_n=d_CB_{n-1}\quad\forall n>1.
$$
\end{prop}

\vskip 0.3cm
\begin{rema}

\

L'espace $S^+(\mathfrak g[1])$ a, en plus de sa structure de cog\`ebre libre, une structure d'alg\`ebre commutative gradu\'ee libre. Nous n'utiliserons pas cette structure qui peut s'interpr\'eter comme issue de l'op\'erade $Com$ qui est duale de l'op\'erade $Lie$.
\end{rema}

\
\section{$C_\infty$ alg\`ebres et cohomologie de Harrison}

\

\subsection{Cohomologie de Harrison des alg\`ebres commutatives}

\

\

La m\`ethode de la section pr\'ec\'edente s'applique aussi aux alg\`ebres commutatives. Soit $A$ une alg\`ebre associative et commutative et $V$ un $A$ module vu comme un bimodule tel que $av=va$ pour tout $v$ de $V$ et tout $a$ de $A$. La cohomologie de Harrison de $A$ \`a valeurs dans $V$ est d\'efinie de la fa\c con suivante.

On d\'efinit d'abord les $p,q$ battements de $n=p+q$ lettres comme les permutatiosn $\sigma$ de $\{1,\dots,n\}$ telles que
$\sigma(1)<\dots<\sigma(p)$ et $\sigma(p+1)<\dots<\sigma(p+q)$. On appelle $Bat(p,q)$ l'ensemble de tous ces battements et on d\'efinit le produit battement de deux tenseurs $\alpha=\alpha_1\otimes\dots\otimes\alpha_p$ et $\beta=\alpha_{p+1}\otimes\dots\otimes\alpha_{p+q}$ par
$$
bat_{p,q}(\alpha,\beta)=\sum_{\sigma\in Bat(p,q)}\varepsilon(\sigma^{-1}) \alpha_{\sigma^{-1}(1)}\otimes\dots\otimes \alpha_{\sigma^{-1}(n)}.
$$
Ceci repr\'esente la somme sign\'ee de tous les tenseurs $\alpha_{i_1}\otimes\dots\otimes \alpha_{i_n}$ dans lesquels les vecteurs $\alpha_1,\dots,\alpha_p$ et les vecteurs $\alpha_{p+1},\dots,\alpha_{p+q}$ apparaissent rang\'es dans leur ordre naturel.\\

Par d\'efinition, l'espace $\underline{A^{\otimes n}}$ est le quotient de $A^{\otimes n}$
par la somme de toutes les images des applications lin\'eaires $bat_{p,n-p}$ ($0<p<n$) (voir \cite{[G]}, \cite{[L]}).
Une $n$ cocha\^\i ne $C$ est une application lin\'eaire de $\underline{A^{\otimes n}}$ dans $V$.
L'espace de ces cocha\^\i nes est not\'e $C^n(A,V)$. L'op\'erateur de cobord de Harrison est l'application
$d_{Ha}:C^{n-1}(A,V)\longrightarrow C^n(A,V)$ d\'efinie par
$$\begin{aligned}
d_{Ha}C(\underline{\alpha_1\otimes\dots\otimes\alpha_n})=&\alpha_1C(\underline{\alpha_2\otimes\dots\otimes \alpha_n})-C(\underline{\alpha_1\alpha_2\otimes\dots\otimes \alpha_n})+\\
&+ C(\underline{\alpha_1\otimes \alpha_2\alpha_3\otimes\dots\otimes \alpha_n})+\dots+(-1)^{n-1}C(\underline{\alpha_1\otimes\dots\otimes \alpha_{n-1}\alpha_n})+\\ &+(-1)^nC(\underline{\alpha_1\otimes\dots\otimes \alpha_{n-1}})\alpha_n.
\end{aligned}
$$
(On a bien s\^ur not\'e $\underline{\alpha_1\otimes\dots\otimes \alpha_n}$ la classe de $\alpha_1\otimes\dots\otimes \alpha_n$ dans $\underline{A^{\otimes n}}$.)

On a $d_{Ha}\circ d_{Ha}=0$, le noyau de $d_{Ha}:C^n(A,V)\longrightarrow C^{n+1}(A,V)$ est not\'e $Z^n(A,V)$, c'est l'espace des $n$ cocycles, l'image de $d_{Ha}:C^{n-1}(A,V)\longrightarrow C^n(A,V)$ est not\'e $B^n(A,V)$, c'est l'espace des $n$ cobords. Le $n^{eme}$ espace de cohomologie de Harrison de $A$ \`a valeurs dans $V$ est le quotient $H^n(A,V)$ de $Z^n(A,V)$ par $B^n(A,V)$.\\


\subsection{La cog\`ebre $(\underline\bigotimes^+(A[1]),\delta)$}

\

\

La construction de la section pr\'ec\'edente a un \'equivalent ici. On commence comme plus haut par d\'ecaler les degr\'es et consid\'erer l'espace $A[1]$. La construction suivante est valable lorsque $A$ est gradu\'ee. Le deg\'e de $\alpha$, $\beta$ dans $A[1]$ est not\'e $a$, $b$.

Le produit battement dans $T(A[1])$ est d\'efini par:
$$
bat_{p,q}(\alpha,\beta)=\sum_{\sigma\in Bat(p,q)}\varepsilon\left(\begin{smallmatrix}a_1&\dots&a_n\\ a_{\sigma^{-1}(1)}&\dots&a_{\sigma^{-1}(n)}\end{smallmatrix}\right)\alpha_{\sigma^{-1}(1)}\otimes\dots\otimes\alpha_{\sigma^{-1}(n)}.
$$
On notera $\underline{\bigotimes}^n(A[1])$ le quotient de $A[1]^{\otimes n}$ par la somme des images des applications $bat_{p,n-p}$ ($0<p<n$) et de fa\c con abusive $\alpha_{[1,n]}=\alpha_1\underline{\otimes}\dots\underline{\otimes}\alpha_n$ la classe de $\alpha_1\otimes\dots\otimes \alpha_n$, lorsque les $\alpha_i$ appartiennent \`a $A[1]$. Ceci ne veut pas dire que $\underline\otimes$ soit une multiplication associative dans $\underline{\bigotimes}^+(A[1])= \sum_{n>0} \underline{ \bigotimes}^n(A[1])$.

\vskip 0.3cm
\begin{rema}\label{r1}
\begin{itemize}

\vskip0.25cm

\item[1.] En fait ce dernier espace peut \^etre muni d'une structure d'alg\`ebre de Lie libre mais nous
n'utiliserons pas cette structure.\\

\item[2.] Une base de l'espace $\underline{\bigotimes}^n(A[1])$. Prenons une base $(e_i)_i$ de $A[1]$ compos\'ee d'\'el\'ements homog\`enes.

Pour chaque suite croissante ${\bf i}=(i_1\leq\dots\leq i_p)$, on pose $e_{\bf i}=e_{i_1}\otimes\dots\otimes e_{i_p}$ et pour chaque $\sigma$ de $S_p$, $e_{\sigma(\bf{i})}=e_{\sigma(i_1)}\otimes\dots\otimes e_{\sigma(i_p)}$.

On note $V(e_{\bf i})$ l'espace engendr\'e par ces vecteurs dans $A[1]^{\otimes p}$ et $W(e_{\bf i})$ le sous-espace:
$$
W(e_{\bf i})=Vect\Big(\sum_\sigma bat_{r,s}\big((e_{i_{\sigma(1)}},\dots,e_{i_{\sigma(r)}}),(e_{i_{\sigma(r+1)}},\dots,e_{i_{\sigma(r+s)}})\big)\Big).
$$
On choisit enfin pour chaque $e_{\bf i}$, une base d'un suppl\'ementaire de $W(e_{\bf i})$ dans $V(e_{\bf i})$ de la forme:
$$
\mathcal B(e_{\bf i})=\left\{e_{\sigma(\bf i)},~~\sigma\in\Sigma(e_{\bf i})\right\}.
$$
Une base de $\underline\bigotimes^n(A[1])$ est donn\'ee par
$$
\mathcal B=\bigcup_{|{\bf i}|=n}\bigcup_{\sigma\in\Sigma(e_{\bf i})}\left\{e_{\sigma(\bf{i})}\right\}.
$$
\end{itemize}
\end{rema}

\vskip 0.3cm
Rappelons maintenant les propri\'et\'es du produit battement.

\begin{lem} {\rm (Associativit\'e de $bat$)}

\

Le produit battement est associatif et commutatif gradu\'e de degr\'e 0: pour tout $\alpha\in A[1]^{\otimes p}$,
$\beta\in A[1]^{\otimes q}$, $\gamma\in A[1]^{\otimes r}$,
$$
\begin{aligned}
(i)&\quad bat_{p,q}(\alpha,\beta)=(-1)^{ab}bat_{q,p}(\beta,\alpha),\\
(ii)&\quad bat_{p+q,r}(bat_{p,q}(\alpha,\beta),\gamma)=bat_{p,q+r}(\alpha,bat_{q,r}(\beta,\gamma)).
\end{aligned}
$$
\end{lem}

\noindent
{\bf Preuve}

\noindent
$(i)$ Soient $\alpha=\alpha_1\otimes\dots\otimes \alpha_p$ et $\beta=\alpha_{p+1}\otimes\dots\otimes \alpha_{p+q}$. On a
$$
bat_{p,q}(\alpha,\beta)=\sum_{\sigma\in Bat(p,q)}\varepsilon\left(\begin{smallmatrix}a_1&\dots&a_{p+q}\\ a_{\sigma^{-1}(1)}&\dots&a_{\sigma^{-1}(p+q)}\end{smallmatrix}\right)\alpha_{\sigma^{-1}(1)}\otimes\dots\otimes \alpha_{\sigma^{-1}(p+q)}.
$$

Pour chaque $\sigma\in Bat(p,q)$ , on construit deux permutations $\tau$ et $\rho$ de $S_{p+q}$ en posant:
$$\tau(k)=\left\{
\begin{array}{ll}
k+p,&\hbox{si $1\leq k\leq q$}\\
k-q,&\hbox{si $q< k\leq q+p$.}
\end{array}
\right.\qquad
\text{et}\qquad\rho=\sigma\circ\tau.
$$

On v\'erifie que $\rho$ appartient \`a $Bat(q,p)$ et que l'application $\sigma\mapsto\rho$ est une bijection de $Bat(p,q)$ sur $Bat(q,p)$.

Posons $\beta_k=\alpha_{\tau(k)}$, pour tout $k$ ($1\leq k\leq p+q$), on a $\beta_{\rho^{-1}(k)}= \alpha_{\tau\circ\rho^{-1}(k)}=\alpha_{\sigma^{-1}(k)}$ et
$$\begin{aligned}
\varepsilon\left(\begin{smallmatrix}b_1&\dots&b_{p+q}\\ b_{\rho^{-1}(1)}&\dots&b_{\rho^{-1}(p+q)}\end{smallmatrix}\right)
&=\varepsilon\left(\begin{smallmatrix}b_1&\dots&b_q&b_{q+1}&\dots&b_{p+q}\\ b_{q+1}&\dots&b_{p+q}&b_1&\dots&b_q\end{smallmatrix}\right)\varepsilon\left(\begin{smallmatrix}a_1&\dots&a_{p+q}\\ a_{\sigma^{-1}(1)}&\dots&a_{\sigma^{-1}(p+q)}\end{smallmatrix}\right)\\
&=\varepsilon\left(\begin{smallmatrix}a_{p+1}&\dots&a_{p+q}&a_1&\dots&a_p\\ a_1&\dots&a_q&a_{q+1}&\dots&a_{p+q}\end{smallmatrix}\right)\varepsilon\left(\begin{smallmatrix}a_1&\dots&a_{p+q}\\ a_{\sigma^{-1}(1)}&\dots&a_{\sigma^{-1}(p+q)}\end{smallmatrix}\right)\\
&=(-1)^{ab}\varepsilon\left(\begin{smallmatrix}a_1&\dots&a_{p+q}\\ a_{\sigma^{-1}(1)}&\dots&a_{\sigma^{-1}(p+q)}\end{smallmatrix}\right).
\end{aligned}
$$

Donc
$$
bat_{p,q}(\alpha,\beta)=(-1)^{ab}bat_{q,p}(\beta,\alpha).
$$

\noindent
$(ii)$ Disons qu'une permutation $\sigma$ de $S_{p+q+r}$ est un $(p,q,r)$-battement si elle v\'erifie
$$
\sigma(1)<\dots<\sigma(p),~\sigma(p+1)<\dots<\sigma(p+q)~~\text{et}~~\sigma(p+q+1)<\dots<\sigma(p+q+r).
$$
Notons $Bat(p,q,r)$ l'ensemble des $(p,q,r)$-battements.

Soient $\alpha=\alpha_1\otimes\dots\otimes \alpha_p$, $\beta=\alpha_{p+1}\otimes\dots\otimes\alpha_{p+q}$ et $\gamma=\alpha_{p+q+1}\otimes\dots\otimes \alpha_{p+q+r}$. On d\'efinit le produit $bat_{p,q,r}(\alpha,\beta,\gamma)$ par:
$$
bat _{p,q,r}(\alpha,\beta,\gamma)=\sum_{\rho\in Bat(p,q,r)}\varepsilon\left(\begin{smallmatrix}a_1&\dots&a_{p+q+r}\\ a_{\rho^{-1}(1)}&\dots&a_{\rho^{-1}(p+q+r)}\end{smallmatrix}\right)\alpha_{\rho^{-1}(1)}\otimes\dots\otimes \alpha_{\rho^{-1}(p+q+r)}.
$$

On a en fait
$$
bat_{p,q,r}(\alpha,\beta,\gamma)=bat_{p+q,r}(bat_{p,q}(\alpha,\beta),\gamma)=bat_{p,q+r}(\alpha,bat_{q,r}(\beta,\gamma)).
$$
Il suffit de montrer que $bat_{p,q,r}(\alpha,\beta,\gamma)=bat_{p+q,r}(bat_{p,q}(\alpha,\beta),\gamma)$, l'autre \'egalit\'e se prouvant de la m\^eme fa\c con.\\

Fixons $\sigma_1\in Bat(p,q)$, on construit une permutation $(\sigma_1\times id)$ sur $\{1,\dots,p+q+r\}$ en posant:
$$
(\sigma_1\times id)(k)=\left\{
\begin{array}{ll}
\sigma_1(k),&\text{si }1\leq k\leq p+q,\\
k,&\text{si }p+q+1\leq k\leq p+q+r.
\end{array}
\right.
$$
Par construction, $(\sigma_1\times id)$ appartient \`a $Bat(p,q,r)$.

Soit maintenant $\sigma_2$ une permutation de $Bat(p+q,r)$, on d\'efinit la permutation $\rho$ de $S_{p+q+r}$ par: $\rho=\sigma_2\circ(\sigma_1\times id)$. On v\'erifie que $\rho$ appartient \`a $Bat(p,q,r)$, que l'application $\varphi:(\sigma_2,\sigma_1)\mapsto\rho=\sigma_2\circ(\sigma_1\times id)$ est une bijection de $Bat(p+q,r)\times Bat(p,q)$ sur $Bat(p,q,r)$ et que
$$
\varepsilon\left(\begin{smallmatrix}a_1&\dots&a_{p+q+r}\\ a_{\rho^{-1}(1)}&\dots&a_{\rho^{-1}(p+q+r)}\end{smallmatrix}\right)=\varepsilon\left(\begin{smallmatrix}a_1&\dots&a_{p+q+r}\\ a_{\sigma_2^{-1}(1)}&\dots&a_{\sigma_2^{-1}(p+q+r)}\end{smallmatrix}\right)\varepsilon\left(\begin{smallmatrix}a_1&\dots&a_{p+q}\\ a_{\sigma_1^{-1}(1)}&\dots&a_{\sigma_1^{-1}(p+q)}\end{smallmatrix}\right).
$$

On a donc bien $bat_{p+q,r}(bat_{p,q}(\alpha,\beta),\gamma)=bat_{p,q,r}(\alpha,\beta,\gamma)$.\\

\hfill$\square$

\

Maintenant on introduit un cocrochet de Lie $\delta$ sur $\underline{\bigotimes}^+(A[1])$ en posant d'abord:
$$\begin{aligned}
\delta(\alpha_1\otimes\dots\otimes\alpha_n)&=\sum_{j=1}^{n-1}\alpha_1\otimes\dots\otimes \alpha_j\bigotimes \alpha_{j+1}\otimes\dots\otimes \alpha_n\\
&-\varepsilon\left(\begin{smallmatrix}a_1\dots a_{n-j}&a_{n-j+1}\dots a_n\\ a_{j+1}\dots a_n&a_1\dots a_j\end{smallmatrix}\right)\alpha_{j+1}\otimes\dots\otimes \alpha_n\bigotimes \alpha_1\otimes\dots\otimes \alpha_j.
\end{aligned}
$$

Cette formule permet de d\'efinir $\delta$ sur l'espace quotient $\underline{\bigotimes}^n(A[1])$. En effet, si $p+q=n$, on a, en posant $\alpha=\alpha_1\otimes\dots\otimes \alpha_p$ et $\beta=\alpha_{p+1}\otimes\dots \alpha_{p+q}$,
$$
\begin{aligned}
\delta(bat_{p,q}(\alpha,\beta))&=\hskip-0.5cm\sum_{\begin{smallmatrix}\sigma\in Bat(p,q)\\0<j<n\end{smallmatrix}}\hskip-0.5cm \varepsilon\left(\begin{smallmatrix}a_1&\dots&a_n\\ a_{\sigma^{-1}(1)}&\dots &a_{\sigma^{-1}(n)}\end{smallmatrix}\right)\\
&\hskip-1cm\Big(\alpha_{\sigma^{-1}(1)}\otimes\dots\otimes\alpha_{\sigma^{-1}(j)}\bigotimes\alpha_{\sigma^{-1}(j+1)}\otimes\dots\otimes \alpha_{\sigma^{-1}(n)}\\
&\hskip-1cm-\varepsilon\left(\begin{smallmatrix}a_{I_\sigma}&a_{J_\sigma}\\ a_{J_\sigma}&a_{I_\sigma}\end{smallmatrix}\right) \alpha_{\sigma^{-1}(j+1)}\otimes\dots\otimes \alpha_{\sigma^{-1}(n)}\bigotimes \alpha_{\sigma^{-1}(1)}\otimes\dots\otimes \alpha_{\sigma^{-1}(j)}\Big).
\end{aligned}
$$
Dans cette formule, on a pos\'e $I_\sigma=\{\sigma^{-1}(1),\dots,\sigma^{-1}(j)\}$ et $J_\sigma=\{\sigma^{-1}(j+1),\dots,\sigma^{-1}(n)\}$.

Posons maintenant $I_\sigma^k=I_\sigma\cap\{1,\dots,p\}$, $I_\sigma^{j-k}=I_\sigma\cap\{p+1,\dots,n\}$ et de m\^eme  $J_\sigma^r=J_\sigma\cap\{1,\dots,p\}$, $J_\sigma^{n-j-r}=J_\sigma\cap\{p+1,\dots,n\}$, ($k=|I_\sigma^k|$ et $r=|J_\sigma^r|$). On peut alors \'ecrire:
$$
\begin{aligned}
\delta(bat_{p,q}(\alpha,\beta))=\sum_{\begin{smallmatrix}0<j<n\\ I\sqcup J=\{1,\dots,n\}\\ I,J\neq\emptyset\end{smallmatrix}}\sum_{\begin{smallmatrix}\sigma\in Bat(p,q)\\ I_\sigma=I\end{smallmatrix}}\varepsilon\left(\begin{smallmatrix}a_{\{1,\dots,n\}}\\ a_{I_\sigma}~a_{J_\sigma}\end{smallmatrix}\right)\Big(\alpha_{I_\sigma}\bigotimes \alpha_{J_\sigma} -\varepsilon\left(\begin{smallmatrix}a_{I_\sigma}&a_{J_\sigma}\\ a_{J_\sigma}&a_{I_\sigma}\end{smallmatrix}\right)\alpha_{J_\sigma}\bigotimes \alpha_{I_\sigma}\Big).
\end{aligned}
$$

\begin{itemize}
\item{Cas 1}\quad $I\neq\{1,\dots,p\}$ ou $I\neq\{p+1,\dots,n\}$.

On v\'erifie qu'alors l'application $(\sigma|_{I_\sigma},\sigma|_{J_\sigma})\mapsto\sigma=\sigma|_{I_\sigma}\otimes\sigma|_{J_\sigma}$ est une bijection entre $Bat(k,j-k)\times Bat(r,n-j-r)$ et $\{\sigma\in Bat(p,q)/~I_\sigma=I\}$.

Dans ce cas, la seconde somme est un produit $\bigotimes$ de produits battements. Elle est donc nulle lorsque l'on passe au quotient.\\

\item{Cas 2}\quad $I=\{1,\dots,p\}$ ou $I=\{p+1,\dots,n\}$.

Les termes restant s'\'ecrivent
$$
\begin{aligned}
\delta(bat_{p,q}(\alpha,\beta))&=\alpha\bigotimes\beta-(-1)^{ab}\beta\bigotimes\alpha\\
&+(-1)^{ab}(\beta\bigotimes\alpha-(-1)^{ab}\alpha\bigotimes\beta)=0.
\end{aligned}
$$
\end{itemize}

On notera $\alpha_{[i,j]}=\alpha_i\underline\otimes\dots\underline\otimes\alpha_j$, $\delta$ est donc le cocrochet de $\underline\bigotimes^+(A[1])$ donn\'ee par
$$
\delta(\alpha_{[1,n]})=\sum_{0<j<n}\alpha_{[1,j]}\bigotimes\alpha_{[j+1,n]}-(-1)^{a_{[1,j]}a_{[j+1,n]}}\alpha_{[j+1,n]}\bigotimes\alpha_{[1,j]}.
$$

\begin{prop} {\rm (La structure de cog\`ebre)}

\

L'espace $\underline\bigotimes^+(A[1])$ \'equipp\'e de $\delta$ est une cog\`ebre de Lie, c'est \`a dire que $\delta$ est coantisym\'etrique de degr\'e $0$ et v\'erifie l'identit\'e de coJacobi: si $\tau$ est la volte,
$$
\tau\circ\delta=-\delta,\qquad\Big(id^{\otimes3}+(\tau\otimes id)\circ(id\otimes\tau)+(id\otimes\tau)\circ(\tau\otimes id)\Big)\circ(\delta\otimes id)\circ\delta=0.
$$
\end{prop}

\vskip 0.3cm
\noindent
{\bf Preuve}\\

D'une part, en notant toujours $\tau$ la volte, on a
$$
\begin{aligned}
\delta(\alpha)&=\delta(\alpha_1\underline{\otimes}\dots\underline{\otimes}\alpha_n)=\sum_{j=1}^{n-1}\Big(
\alpha_1\underline{\otimes}\dots\underline{\otimes}\alpha_j\bigotimes\alpha_{j+1}\underline{\otimes}\dots\underline{\otimes}\alpha_n\\
&\hskip 1cm-\varepsilon\left(\begin{smallmatrix}a_{\{1,\dots,n\}}\\ a_{\{j+1,\dots,n\}}~a_{\{1,\dots,j\}}\end{smallmatrix}\right) \alpha_{j+1}\underline{\otimes}\dots\underline{\otimes}\alpha_n\bigotimes \alpha_1\underline{\otimes}\dots\underline{\otimes}\alpha_j\Big)\\
&=\sum_{0<j<n}\Big(\alpha_{[1,j]}\bigotimes \alpha_{[j+1,p]}-\tau\big(\alpha_{[1,j]}\bigotimes \alpha_{[j+1,p]}\big)\Big).
\end{aligned}
$$
Donc
$$
\tau\circ\delta(\alpha)=\sum_{j=1}^{n-1}\tau\big(\alpha_{[1,j]}\bigotimes \alpha_{[j+1,p]}\big)-\alpha_{[1,j]}\bigotimes \alpha_{[j+1,p]}=-\delta(\alpha).
$$
D'autre part, on a
$$
\begin{aligned}
(\delta\otimes id)\circ\delta(\alpha)=&\sum_{0<i<j<n}\alpha_{[1,i]}\bigotimes \alpha_{[i+1,j]}\bigotimes \alpha_{[j+1,n]}\\
&\quad-\varepsilon\left(\begin{smallmatrix}a_{[1,i]}&a_{[i+1,j]}&a_{[j+1,n]}\\ a_{[i+1,j]}&a_{[1,i]}&a_{[j+1,n]}\end{smallmatrix}\right)\alpha_{[i+1,j]}\bigotimes \alpha_{[1,i]}\bigotimes \alpha_{[j+1,n]}\\
&\quad-\varepsilon\left(\begin{smallmatrix}a_{[1,i]}&a_{[i+1,j]}&a_{[j+1,n]}\\
a_{[i+1,j]}&a_{[j+1,n]}&a_{[1,i]}\end{smallmatrix}\right) \alpha_{[i+1,j]}\bigotimes\alpha_{[j+1,n]}\bigotimes\alpha_{[1,i]}\\
&\quad+\varepsilon\left(\begin{smallmatrix}a_{[1,i]}&a_{[i+1,j]}&a_{[j+1,n]}\\ a_{[j+1,n]}&a_{[i+1,j]}&a_{[1,i]}\end{smallmatrix}\right) \alpha_{[j+1,n]}\bigotimes \alpha_{[i+1,j]}\bigotimes \alpha_{[1,i]}.
\end{aligned}
$$

Donc, en notant $(i)=[1,i]$, $(j)=[i+1,j]$ et $(n)=[j+1,n]$,
$$
\begin{aligned}
\Big(id^{\otimes3}&+(\tau\otimes id)\circ(id\otimes\tau)+(id\otimes\tau)\circ(\tau\otimes id)\Big)\circ(\delta\otimes id)\circ\delta(\alpha)\\
&=\sum_{0<i<j<n}\alpha_{(i)}\bigotimes \alpha_{(j)}\bigotimes \alpha_{(n)}-\varepsilon\left(\begin{smallmatrix}a_{(i)}&a_{(j)}&a_{(n)}\\ a_{(j)}&a_{(i)}&a_{(n)}\end{smallmatrix}\right) \alpha_{(j)}\bigotimes \alpha_{(i)}\bigotimes \alpha_{(n)}\\
&-\varepsilon\left(\begin{smallmatrix}a_{(i)}&a_{(j)}&a_{(n)}\\ a_{(j)}&a_{(n)}&a_{(i)}\end{smallmatrix}\right) \alpha_{(j)}\bigotimes \alpha_{(n)}\bigotimes \alpha_{(i)}+\varepsilon\left(\begin{smallmatrix}a_{(i)}&a_{(j)}&a_{(n)}\\ a_{(n)}&a_{(j)}&a_{(i)}\end{smallmatrix}\right)\alpha_{(n)}\bigotimes
\alpha_{(j)}\bigotimes \alpha_{(i)}\\
&+\varepsilon\left(\begin{smallmatrix}a_{(i)}&a_{(j)}&a_{(n)}\\ a_{(n)}&a_{(i)}&a_{(j)}\end{smallmatrix}\right) \alpha_{(n)}\bigotimes \alpha_{(i)}\bigotimes \alpha_{(j)}-\varepsilon\left(\begin{smallmatrix}a_{(i)}&a_{(j)}&a_{(n)}\\ a_{(n)}&a_{(j)}&a_{(i)}\end{smallmatrix}\right)\alpha_{(n)}\bigotimes
\alpha_{(j)}\bigotimes \alpha_{(i)}\\
&-\alpha_{(i)}\bigotimes \alpha_{(j)}\bigotimes \alpha_{(n)}+\varepsilon\left(\begin{smallmatrix}a_{(i)}&a_{(j)}&a_{(n)}\\ a_{(i)}&a_{(n)}&a_{(j)}\end{smallmatrix}\right) \alpha_{(i)}\bigotimes \alpha_{(n)}\bigotimes \alpha_{(j)}\\
&+\varepsilon\left(\begin{smallmatrix}a_{(i)}&a_{(j)}&a_{(n)}\\ a_{(j)}&a_{(n)}&a_{(i)}\end{smallmatrix}\right) \alpha_{(j)}\bigotimes \alpha_{(n)}\bigotimes \alpha_{(i)}-\varepsilon\left(\begin{smallmatrix}a_{(i)}&a_{(j)}&a_{(n)}\\ a_{(i)}&a_{(n)}&a_{(j)}\end{smallmatrix}\right)\alpha_{(i)}\bigotimes
\alpha_{(n)}\bigotimes \alpha_{(j)}\\
&-\varepsilon\left(\begin{smallmatrix}a_{(i)}&a_{(j)}&a_{(n)}\\ a_{(n)}&a_{(i)}&a_{(j)}\end{smallmatrix}\right) \alpha_{(n)}\bigotimes \alpha_{(i)}\bigotimes \alpha_{(j)}+\varepsilon\left(\begin{smallmatrix}a_{(i)}&a_{(j)}&a_{(n)}\\ a_{(j)}&a_{(i)}&a_{(n)}\end{smallmatrix}\right)\alpha_{(j)}\bigotimes
\alpha_{(i)}\bigotimes \alpha_{(n)}\\
&=0.
\end{aligned}
$$

\subsection{Morphismes et cod\'erivations}

\

\

La structure de cog\`ebre de Lie de $(\underline\bigotimes^+(A[1]),\delta)$ est libre. C'est \`a dire que si $(\mathcal{C},c)$ est une cog\`ebre de Lie nilpotente quelconque, tout $f:(\mathcal{C},c)\longrightarrow A[1]$ lin\'eaire se prolonge en $F:(\mathcal{C},c)\longrightarrow \underline\bigotimes^+(A[1])$ qui est un morphisme de cog\`ebre. Nous montrons ici comment d\'efinir des cod\'erivations $Q$ et des morphismes $F$ de cette structure \`a partir de leurs `s\'erie de Taylor'.

Soit $F:\underline\bigotimes^+(A[1])\longrightarrow\underline\bigotimes^+(B[1])$ un morphisme de cog\`ebres de Lie. On suppose toujours $F$ homog\`ene de degr\'e 0. On appelle $F_n$ la projection sur $B[1]$ parall\`element \`a $\underline\bigoplus^{k>1}B[1]$ de la restriction de $F$ \`a $\underline\bigotimes^n(A[1])$: $F_n$ est une application lin\'eaire de $\underline\bigotimes^n(A[1])$ dans $B[1]$.

De m\^eme soit $Q:\underline\bigotimes^+(A[1])\longrightarrow\underline\bigotimes^+(A[1])$ une cod\'erivation de cog\`ebres de Lie. On suppose $Q$ homog\`ene de degr\'e $q$. On appelle $Q_n$ la projection sur $A[1]$ parall\`element \`a $\underline\bigoplus^{k>1}A[1]$ de la restriction de $Q$ \`a $\underline\bigotimes^n(A[1])$: $Q_n$ est une application lin\'eaire de $\underline\bigotimes^n(A[1])$ dans $A[1]$.\\

\begin{prop} {\rm (Reconstruction de $F$ et $Q$)}

\

La suite d'applications $(F_n)$ (resp. $(Q_n)$) permet de reconstruire $F$ (resp. $Q$) de fa\c con unique. On a explicitement
$$
F(\alpha_{[1,n]})=\sum_{\begin{smallmatrix}k>0,~0<r_1,\dots,r_k\\ r_1+\dots+r_k=n\end{smallmatrix}}F_{r_1}(\alpha_{[1,r_1]})\underline\otimes F_{r_2}(\alpha_{[r_1+1,r_1+r_2]})\underline\otimes\dots\underline\otimes F_{r_k}(\alpha_{[n-r_k+1,n]})
$$
et
$$
Q(\alpha_{[1,n]})=\sum_{\begin{smallmatrix}1\leq r\leq n\\0\leq
j\leq
n-r\end{smallmatrix}}(-1)^{qa_{[1,j]}}\alpha_{[1,j]}\underline\otimes
Q_r(\alpha_{[j+1,j+r]})\underline\otimes\alpha_{[j+r+1,n]}.
$$
Plus pr\'ecis\'ement, toute suite d'applications $(\varphi_n)$ peut se relever d'une seule fa\c con en un morphisme (resp. une cod\'erivation).
\end{prop}

\vskip 0.3cm
\noindent
{\bf Preuve}\\

La preuve est semblable \`a celle de \cite{[AMM]}. Pour un morphisme, si tous les $F_n$ sont nuls, $F(\alpha_1)$ est nul pour tout $\alpha_1\in A[1]$, et si tous les $F(\alpha_{[1,p]})$ sont nuls quelque soit $p<n$, alors
$$
\begin{aligned}
\delta\circ F(\alpha_{[1,n]})&=(F\otimes F)\circ\delta(\alpha_{[1,n]})\\
&=\sum_{0<j<n}F(\alpha_{[1,j]})\otimes F(\alpha_{[j+1,n]})-(-1)^{a_{[1,j]}} F(\alpha_{[j+1,n]})\otimes F(\alpha_{[1,j]})=0
\end{aligned}
$$
donc $F(\alpha_{[1,n]})=F_n(\alpha_{[1,n]})\in A[1]$ est aussi nul et par induction, $F$ est nul. Le m\^eme argument s'applique pour une d\'erivation $Q$. Ceci prouve l'unicit\'e.

Il reste juste \`a montrer que les formules de la proposition d\'efinissent bien un morphisme (resp. une cod\'erivation).\\

$F$ est bien d\'efini.\\

D'abord $F$ est bien d\'efini, c'est \`a dire l'application $\overline{F}$ d\'efinie par la m\^eme formule mais sur $\alpha_{\{1,\dots,n\}}=\alpha_1\otimes\dots\otimes \alpha_n$ passe bien au quotient. Il suffit pour cela de montrer que
$$
\overline{F}(bat_{p,q}(\alpha_{\{1,\dots,p\}},\alpha_{\{p+1,\dots,p+q\}}))=0.
$$

En fait, si on pose
$$
\tilde{F}^k(\alpha_{\{1,\dots,n\}})=\sum_{\begin{smallmatrix}0<r_1,\dots,r_k\\ r_1+\dots+r_k=n\end{smallmatrix}}F_{r_1}(\alpha_{[1,r_1]})\otimes F_{r_2}(\alpha_{[r_1+1,r_1+r_2]})\otimes\dots\otimes F_{r_k}(\alpha_{[n-r_k+1,n]}),
$$
alors $\tilde{F}^k (bat_{p,q}(\alpha_{\{1,\dots,p\}},\alpha_{\{p+1,\dots,p+q\}}))$ est une somme de produits battements de la forme
$$
\tilde{F}^k (bat_{p,q}(\alpha_{\{1,\dots,p\}},\alpha_{\{p+1,\dots,p+q\}}))=\sum bat_{r,s}(\tilde{F}^r(\alpha_I),\tilde{F}^s (\alpha_J)),
$$
ce qui prouve que $F$ est bien d\'efini. Prouvons cette derni\`ere relation.

On a
$$
bat_{p,q}(\alpha_{\{1,\dots,p\}},\alpha_{\{p+1,\dots,p+q\}})=\sum_{\sigma\in Bat(p,q)}\varepsilon\left(\begin{smallmatrix}a_{\{1,\dots,p+q\}}\\
a_{\{\sigma^{-1}(1),\dots,\sigma^{-1}(p+q)\}}\end{smallmatrix}\right)\alpha_{\{\sigma^{-1}(1),\dots,\sigma^{-1}(p+q)\}}.
$$
Donc
$$
\begin{aligned}
\tilde{F}^k(bat_{p,q}(\alpha_{\{1,\dots,p\}},&\alpha_{\{p+1,\dots,p+q\}}))=\\&=\hskip-0.5cm\sum_{\begin{smallmatrix}\sigma\in Bat(p,q)\\ r_1,\dots,r_k\end{smallmatrix}}\varepsilon\left(\begin{smallmatrix}a_{\{1,\dots,p+q\}}\\
a_{\{\sigma^{-1}(1),\dots,\sigma^{-1}(p+q)\}}\end{smallmatrix}\right)F_{r_1}(\alpha_{\{\sigma^{-1}(1),\dots,\sigma^{-1}(r_1)\}})\otimes\dots\\
&\hskip 3cm\dots\otimes F_{r_k}(\alpha_{\{\sigma^{-1}(p+q-r_k+1),\dots ,\sigma^{-1}(p+q)\}}).
\end{aligned}
$$
Fixons un battement $\sigma\in Bat(p,q)$. Posons $s_j=r_1+\dots+r_j$. S'il existe un $j$ tel que $I_j=\{\sigma^{-1}(s_{j-1}+1),\dots,\sigma^{-1}(s_j)\}$ n'est inclus ni dans $\{1,\dots,p\}$ ni dans $\{p+1,\dots,p+q\}$, alors cet ensemble peut s'\'ecrire $I_j=\{i_1,\dots,i_{t_j},i_{t_j+1},\dots,i_{r_j}\}$ avec $0<i_1<\dots<i_{t_j}<p+1$ et $p<i_{t_j+1}<\dots<i_{r_j}<p+q+1$ et bien s\^ur $0<t_j<r_j$. L'ensemble $Bat_{I_j}(p,q)$ des battements $\rho\in Bat(p,q)$ tels que $\rho^{-1}(u)=\sigma^{-1}(u)$ pour tous les $u$ de $\{1,\dots,p+q\}\setminus\{s_{j-1}+1,\dots,s_j\}$ est en bijection avec $Bat(t_j,r_j-t_j)$ puisque chaque battement $\mu$ de $Bat(t_j,r_j-t_j)$ peut se prolonger en un battement $\tilde\mu$ de $Bat(p,q)$ d\'efini par $\tilde\mu^{-1}(u)=\sigma^{-1}(u)$ si $u$ est dans $\{1,\dots,p+q\}\setminus\{s_{j-1}+1,\dots,s_j\}$ et $\tilde\mu(i_v)=i_{\mu(v)}$ pour $0<v<r_j+1$. Alors:
$$
\begin{aligned}
\sum_{\rho\in Bat_{I_j}(p,q)}&\varepsilon\left(\begin{smallmatrix}a_{\{1,\dots,p+q\}}\\
a_{\{\rho^{-1}(1),\dots,\rho^{-1}(p+q)\}}\end{smallmatrix}\right)F_{r_1}(\alpha_{\{\rho^{-1}(1),\dots,\rho^{-1}(r_1)\}})\otimes\dots\\
&\hskip 3cm\dots\otimes F_{r_k}(\alpha_{\{\rho^{-1}(p+q-r_k+1), \rho^{-1}(p+q)\}})\\
=\pm&\varepsilon\left(\begin{smallmatrix}a_{\{1,\dots,p+q\}\setminus\{s_{j-1}+1,\dots,s_j\}}\\
a_{\{\sigma^{-1}(1),\dots,\sigma^{-1}(p+q)\}\setminus\{\sigma^{-1}(s_{j-1}+1),\dots,\sigma^{-1}(s_j)\}}\end{smallmatrix}\right)F_{r_1}(\alpha_{I_1})\otimes\dots\otimes F_{r_{j-1}}(\alpha_{I_{r_{j-1}}})\otimes\\
&\otimes\sum_{\mu\in Bat(t_j,r_j-t_j)}\varepsilon\left(\begin{smallmatrix}a_{\{i_1,\dots,i_{r_j}\}}\\ a_{\{i_{\mu^{-1}(1)},\dots,i_{\mu^{-1}(r_j)}\}}\end{smallmatrix}\right)F_{r_j}(\alpha_{i_{\mu^{-1}(1)}},\dots,\alpha_{i_{\mu^{-1}(r_j)}})\otimes\\
&\otimes F_{r_{j+1}}(\alpha_{I_{r_{j+1}}})\otimes\dots\otimes F_{r_k}(\alpha_{I_k})\\
=0&
\end{aligned}
$$
en posant $\beta_v=\alpha_{i_v}$ et en remarquant que $F_{r_j}$ s'annule sur
$$
bat_{(t_j,r_j-t_j)}(\beta_{\{1,\dots,t_j\}},\beta_{\{t_j+1,\dots,r_j\}}).
$$
Il ne reste donc que la somme sur les battements $\sigma$ tels que pour tout $j$, on ait soit $I_j\subset\{1,\dots,p\}$ soit $I_j\subset\{p+1,\dots,p+q\}$. Dans ce cas, les nombres $\sigma^{-1}(s_{j-1}+1),\dots,\sigma^{-1}(s_j)$ sont rang\'es dans leur ordre naturel, puisque $\sigma$ est un battement. Notons maintenant $J_1,\dots,J_r$ les ensembles $I_j$ tels que $I_j\subset\{1,\dots,p\}$ et $J_{s+1},\dots,J_{s+r}$ les autres. Posons $J_j=\{g^j_1<g^j_2<\dots<g^j_{h_j}\}$. La somme sur tous les battements $\sigma$ tels que $\{I_1,\dots,I_k\}=\{J_1,\dots,J_{r+s}\}$ (\`a l'ordre pr\`es) est isomorphe \`a une somme sur tous les $(r,s)$ battements $\nu$: \'etant donn\'e un tel battement $\nu$, on construit le $(p,q)$ battement $\tilde{\nu}$ en posant, pour tout $j$, $1\leq j\leq k$ et tout $t$, $1\leq t\leq r_j$:
$$
\tilde{\nu}^{-1}(s_{j-1}+t)=g^{\nu^{-1}(j)}_t.
$$
Alors, si $s_j=|J_j|$,
$$
\begin{aligned}
\hskip-1cm\sum_{\begin{smallmatrix}\sigma\in Bat(p,q)\\ \{I_1,\dots,I_k\}=\{J_1,\dots,J_k\}\end{smallmatrix}}&\hskip-0.5cm\varepsilon\left(\begin{smallmatrix}a_{\{1,\dots,p+q\}}\\
a_{\{\sigma^{-1}(1),\dots,\sigma^{-1}(p+q)\}}\end{smallmatrix}\right)F_{r_1}(\alpha_{\{\sigma^{-1}(1),\dots,\sigma^{-1}(r_1)\}})\otimes\dots\\
&\hskip 4cm\dots\otimes F_{r_k}(\alpha_{\{\sigma^{-1}(p+q-r_k+1), \sigma^{-1}(p+q)\}})=\\
&=\sum_{\mu\in Bat(r,s)}\varepsilon\left(\begin{smallmatrix}a_{J_1}\dots a_{J_k}\\ a_{J_{\mu^{-1}(1)}}\dots a_{J_{\mu^{-1}(k)}}\end{smallmatrix}\right)
F_{s_{\mu^{-1}(1)}}(\alpha_{J_{\mu^{-1}(1)}})\otimes\dots\otimes F_{s_{\mu^{-1}(k)}}(\alpha_{J_{\mu^{-1}(k)}})\\
&=bat_{r,s}\left(F_{s_1}(\alpha_{J_1})\otimes\dots\otimes F_{s_r}(\alpha_{J_r}),F_{s_{r+1}}(\alpha_{J_{r+1}})\otimes\dots\otimes F_{s_{r+s}}(\alpha_{J_{r+s}})\right).
\end{aligned}
$$
C'est l'\'egalit\'e annonc\'ee.\\

$F$ est un morphisme.\\

En gardant nos notations et en ajoutant $(r_j)=[s_{j-1}+1,s_j]$, on a par d\'efinition:
$$
\begin{aligned}
\delta\circ&F(\alpha_{[1,n]})=\sum_{r_1,\dots,r_k}\sum_{0<j<k}F_{r_1}(\alpha_{(r_1)})\underline\otimes\dots\underline\otimes F_{r_j}(\alpha_{(r_j)})
\bigotimes F_{r_{j+1}}(\alpha_{(r_{j+1})})\underline\otimes\dots\underline\otimes F_{r_k}(\alpha_{(r_k)})\\
&-(-1)^{a_{[1,s_j]}a_{[s_j+1,n]}}F_{r_{j+1}}(\alpha_{(r_{j+1})})\underline\otimes\dots\underline\otimes F_{r_k}(\alpha_{(r_k)})\bigotimes F_{r_1}(\alpha_{(r_1)})\underline\otimes\dots\underline\otimes F_{r_j}(\alpha_{(r_j)}).
\end{aligned}
$$
D'autre part,
$$
\begin{aligned}
(F\otimes F)&\circ\delta(\alpha_{[1,n]})=(F\otimes F)\big(\sum_{0<s<n}\alpha_{[1,s]}\bigotimes\alpha_{[s+1,n]}-(-1)^{a_{[1,s]}a_{[s+1,n]}} \alpha_{[s+1,n]}\bigotimes\alpha_{[1,s]}\big)\\
&\hskip-1cm=\sum_{0<s<n}\sum_{\begin{smallmatrix}r_1,\dots,r_j\\ r_1+\dots+r_j=s\end{smallmatrix}}\sum_{\begin{smallmatrix}r_{j+1},\dots,r_k\\r_1+\dots+r_k=n\end{smallmatrix}}\\ &F_{r_1}(\alpha_{(r_1)})\underline\otimes\dots\underline\otimes F_{r_j}(\alpha_{(r_j)})
\bigotimes F_{r_{j+1}}(\alpha_{(r_{j+1})})\underline\otimes\dots
\underline\otimes F_{r_k}(\alpha_{(r_k)})-\\
-&(-1)^{a_{[1,s]}a_{[s+1,n]}}F_{r_{j+1}}(\alpha_{(r_{j+1})})\underline\otimes\dots\underline\otimes F_{r_k}(\alpha_{(r_k)})\bigotimes F_{r_1}(\alpha_{(r_1)})\underline\otimes\dots\underline\otimes F_{r_j}(\alpha_{(r_j)}).
\end{aligned}
$$
On v\'erifie ais\'ement que chaque terme de la premi\`ere expression appara\^\i t une fois et une seule dans la seconde et r\'eciproquement. On a donc
$$
\delta\circ F=(F\otimes F)\circ\delta.
$$

$Q$ est bien d\'efini.\\

Il s'agit l\`a encore de montrer que la d\'efinition de $\tilde Q$:
$$
\tilde{Q}(\alpha_{\{1,\dots,n\}})=\sum_{\begin{smallmatrix}0<r\\1\leq j\leq n-r\end{smallmatrix}} (-1)^{qa_{\{1,\dots,j\}}}\alpha_{\{1,\dots,j\}}\otimes Q_r(\alpha_{\{j+1,\dots,j+r\}})\otimes\alpha_{\{j+r+1,\dots,n\}}
$$
passe au quotient. On calcule donc $Q\big(bat_{p,q}(\alpha_{\{1,\dots,p\}},\alpha_{\{p+1,\dots,p+q\}})\big)$. Le m\^eme argument que ci-desssus, nous dit d'abord qu'il ne reste que
$$
\begin{aligned}
\sum_{\begin{smallmatrix}0<r\\1\leq j\leq n-r\end{smallmatrix}}\sum_{\begin{smallmatrix}I\subset\{1,\dots,p\}\\
\text{ou}\\
I\subset\{p+1,\dots,p+q\}\end{smallmatrix}}&\sum_{\begin{smallmatrix}\sigma\in Bat(p,q)\\
\sigma^{-1}(\{j+1,\dots,j+r\})=I\end{smallmatrix}}(-1)^{qa_{\{\sigma^{-1}(1),\dots,\sigma^{-1}(j)\}}}\varepsilon\left(\begin{smallmatrix}a_{\{1,\dots,p+q\}}\\ a_{\{\sigma^{-1}(1),\dots,\sigma^{-1}(p+q)\}}\end{smallmatrix}\right)\\
&\alpha_{\{\sigma^{-1}(1),\dots,\sigma^{-1}(j)\}}\otimes Q_r(\alpha_I)\otimes\alpha_{\{\sigma^{-1}(j+r+1),\dots,\sigma^{-1}(p+q)\}}.
\end{aligned}
$$
Ensuite, comme $\sigma$ est un battement, les \'el\'ements de $I$ sont rang\'es dans leur ordre naturel $I=\{t,t+1,\dots,t+r-1\}$. Supposons (par exemple) que $I\subset\{1,\dots,p\}$, alors les $\alpha$ d'indices dans $\{1,\dots,p\}$ et pr\'ec\'edant ceux de $I$ apparaissent avant le terme en $Q_r$. On pose donc $t=\sigma^{-1}(s_{j-1}+1)$ et:
$$
\left\{\begin{aligned}
\beta_i&=\alpha_i~\text{ si }~\sigma^{-1}(i)\notin I\\
\beta_t&=Q_r(\alpha_I)\end{aligned}\right.
$$
Il y a donc $p+q-r+1$ $\beta$, indic\'es par $\{1,\dots,t-1,t,t+r,\dots,p,p+1,\dots,p+q\}$. On a $b_i=a_i$ si $\sigma^{-1}(i)\notin I$ et $b_t=q+a_I$. Les battements $\sigma$ consid\'er\'es sont en bijection avec les battements $\rho$ qu'ils induisent sur les $\beta$. On a
$$
\begin{aligned}
(-1)^{qa_{\{\sigma^{-1}(1),\dots,\sigma^{-1}(j)\}}}&\varepsilon\left(\begin{smallmatrix}a_{\{1,\dots,p+q\}}\\ a_{\{\sigma^{-1}(1),\dots,\sigma^{-1}(p+q)\}}\end{smallmatrix}\right)=\\
=&(-1)^{qa_{\{1,\dots,t-1\}}}\varepsilon\left(\begin{smallmatrix}b_{\{1,\dots,t,t+r,\dots,p+q\}}\\ b_{\{\rho^{-1}(1),\dots,\rho^{-1}(t),\rho^{-1}(t+r),\dots,\rho^{-1}(p+q)\}}\end{smallmatrix}\right).
\end{aligned}
$$
Donc
$$
\begin{aligned}
Q\big(bat_{p,q}&(\alpha_{\{1,\dots,p\}},\alpha_{\{p+1,\dots,p+q\}})\big)=\\
=&\sum_{0<r}\sum_{t,~t+r-1\leq p}(-1)^{qa_{\{1,\dots,t-1\}}}bat_{p-r+1,q}(\beta_{\{1,\dots,t,t+r,\dots,p\}},\beta_{\{p+1,\dots,p+q\}})\\
&+\sum_{p<t}(-1)^{qa_{\{p+1,\dots,t-1\}}}bat_{p,q-r+1}(\beta_{\{1,\dots,p\}},\beta_{\{p+1,\dots,t,t+r,\dots,p+q\}}).
\end{aligned}
$$
Comme ci-dessus, $Q$ est donc bien d\'efinie.\\

$Q$ est une cod\'erivation.\\

On a
$$
\delta\circ Q(\alpha_{[1,n]})=\delta\left(\sum_{r,j}(-1)^{qa_{[1,j]}}\alpha_{[1,j]}\underline\otimes Q_r(\alpha_{[j+1,j+r]})\underline\otimes \alpha_{[j+r+1,n]}\right).
$$
Donc
$$
\begin{aligned}
\delta\circ Q&(\alpha_{[1,n]})=\sum_{r,0<k<j}(-1)^{qa_{[1,j]}}\alpha_{[1,k]}\bigotimes\alpha_{[k+1,j]}\underline\otimes Q_r(\alpha_{[j+1,j+r]})\underline\otimes \alpha_{[j+r+1,n]}-\\
&-(-1)^{qa_{[1,j]}+a_{[1,k]}(a_{[k+1,n]}+q)}\alpha_{[k+1,j]}\underline\otimes Q_r(\alpha_{[j+1,j+r]})\underline\otimes \alpha_{[j+r+1,n]}\bigotimes \alpha_{[1,k]}+\\
&+\sum_{r,0<j<k-r}(-1)^{qa_{[1,j]}}\alpha_{[1,j]}\underline\otimes Q_r(\alpha_{[j+1,j+r]})\underline\otimes\alpha_{[j+r+1,k]}\bigotimes \alpha_{[k+1,n]}-\\
&-(-1)^{qa_{[1,j]}+(a_{[1,k]}+q)a_{[k+1,n]}}\alpha_{[k+1,n]}\bigotimes\alpha_{[1,j]}\underline\otimes Q_r(\alpha_{[j+1,j+r]})\underline\otimes \alpha_{[j+r+1,k]}.
\end{aligned}
$$
Par ailleurs,
$$
\begin{aligned}
(id\otimes Q&+Q\otimes id)\circ\delta(\alpha_{[1,n]})=\\
&=(id\otimes Q+Q\otimes id)\sum_{0<k<n}\alpha_{[1,k]}\bigotimes\alpha_{[k+1,n]}-(-1)^{a_{[1,k]} a_{[k+1,n]}} \alpha_{[k+1,n]}\bigotimes\alpha_{[1,k]}
\end{aligned}
$$
Donc
$$
\begin{aligned}
(id\otimes& Q+Q\otimes id)\circ\delta(\alpha_{[1,n]})=\\
=&\hskip-0.3cm\sum_{0<k<j<n-r+1}\hskip-0.3cm(-1)^{q(a_{[1,k]}+a_{[k+1,j]})}\alpha_{[1,k]}\bigotimes\alpha_{[k+1,j]} \underline\otimes Q_r(\alpha_{[j+1,j+r]})\underline\otimes\alpha_{[j+r+1,n]}\\
&-\hskip -1cm\sum_{0<j<k-r+1<n-r+1}\hskip-0.3cm(-1)^{q(a_{[k+1,n]}+a_{[1,j]})+a_{[1,k]}a_{[k+1,n]}}\alpha_{[k+1,n]}\bigotimes\\
&\hskip 6cm\bigotimes\alpha_{[1,j]}\underline\otimes Q_r(\alpha_{[j+1,j+r]})\underline\otimes\alpha_{[j+r+1,k]}\\
&+\hskip-0.3cm\sum_{0<j<k-r+1<n-r+1}\hskip-0.5cm(-1)^{qa_{[1,j]}}\alpha_{[1,j]}\underline\otimes Q_r(\alpha_{[j+1,j+r]})\underline\otimes\alpha_{[j+r+1,k]}\bigotimes \alpha_{[k+1,n]}\\
&-\hskip-0.3cm\sum_{0<k<j<n-r+1}\hskip-0.3cm(-1)^{qa_{[k+1,j]}+a_{[1,k]}a_{[k+1,n]}}\alpha_{[k+1,j]} \underline\otimes Q_r(\alpha_{[j+1,j+r]})\underline\otimes\alpha_{[j+r+1,n]}\bigotimes \alpha_{[1,k]}.
\end{aligned}
$$
$Q$ est donc une cod\'erivation et la proposition est prouv\'ee.\hfill$\square$

\


\subsection{$C_\infty$ alg\`ebre, morphismes de $C_\infty$ alg\`ebre}

\

\

Lorsque $A$ est une alg\`ebre commutative, $A[1]$ est muni d'un produit $m_2$ d\'efini par $m_2(\alpha,\beta)=(-1)^a\alpha\beta$ qui devient de degr\'e 1, anticommutatif et antiassociatif:
$$
m_2(\beta,\alpha)=-(-1)^{ab}m_2(\alpha,\beta),\qquad m_2(m_2(\alpha,\beta),\gamma)=-(-1)^am_2(\alpha,m_2(\beta,\gamma)).
$$

Le produit $m_2$ \'etant anticommutatif est d\'efini de $\underline{\bigotimes}^2(A[1])$ dans $A[1]$. Il se prolonge donc, gr\^ace \`a la proposition pr\'ec\'edente, en une unique cod\'erivation $m$ du cocrochet $\delta$. Comme $m$ est de degr\'e 1, le prolongement \`a $\underline\bigotimes^n(A[1])$ est
$$
m(\alpha_{[1,n]})=\sum_{0<k<n}(-1)^{a_{[1,k-1]}}\alpha_{[1,k-1]}\underline\otimes m_2(\alpha_k,\alpha_{k+1})\underline\otimes\alpha_{[k+2,n]}.
$$
\vskip 0.3cm
\begin{lem} {\rm (Propri\'et\'es de $m$)}

\

$m$ est l'unique cod\'erivation de $\delta$ qui prolonge $m_2$ sur $\underline{\bigotimes}^2(A[1])$. Elle v\'erifie $m\circ m=0$.\\
\end{lem}

\vskip 0.3cm
\noindent
{\bf Preuve}

Remarquons que, $m$ \'etant une cod\'erivation de degr\'e impair, $m\circ m$ est aussi une cod\'erivation. Avec les notations pr\'ec\'edentes, $(m\circ m)_k=0$ si $k\neq3$ et, puisque $m_2$ est antiassociative,
$$
(m\circ m)_3(\alpha_1\underline\otimes\alpha_2\underline\otimes\alpha_3)=m_2(m_2(\alpha_1\underline\otimes\alpha_2)\underline\otimes\alpha_3+(-1)^{a_1} \alpha_1\underline\otimes m_2(\alpha_2\underline\otimes\alpha_3))=0.
$$
Par unicit\'e de la cod\'erivation qui prolonge les $(m\circ m)_k$, on en d\'eduit que $m\circ m=0$.\hfill$\square$\\

\

\begin{defn} {\rm ($C_\infty$ alg\`ebre)}

\

Une $C_\infty$ alg\`ebre est une cog\`ebre diff\'erentielle de la forme $(\underline\bigotimes^+(A[1]),\delta,m)$ o\`u $\delta$ est le cocrochet de Lie d\'efini ci-dessus et $m$ est une cod\'erivation  de $\delta$ de carr\'e nul.

Si $A$ est une alg\`ebre commutative gradu\'ee, la cog\`ebre de Lie diff\'erentielle $C(A)=\left(\underline\bigotimes^+(A[1]),\delta,m\right)$ o\`u
$m_k=0$ pour tout $k\neq2$ et $m_2(\alpha\underline\otimes\beta)=(-1)^a\alpha\wedge\beta$ s'appelle la $C_\infty$ alg\`ebre enveloppante de $(A,\wedge)$.

Un morphisme de $C_\infty$ alg\`ebres $A$ et $B$ est un morphisme de cog\`ebres de Lie $F:\underline\bigotimes^+(A[1])\longrightarrow\underline\bigotimes^+(B[1])$ tel que $m^B\circ F=F\circ m^A$.
\end{defn}

Puisqu'une cod\'erivation $m$ est caract\'eris\'ee par la suite des $m_k$, on voit qu'une $C_\infty$ alg\`ebre est la $C_\infty$ alg\`ebre enveloppante d'une alg\`ebre commutative si et seulement si elle est telle que $m_k=0$ pour tout $k\neq2$.

De m\^eme, un morphisme d'alg\`ebres commutatives $f:A\longrightarrow B$ se rel\`eve d'une fa\c con et une seule en un morphisme de leur $C_\infty$ alg\`ebres enveloppantes $F:\underline\bigotimes^+(A[1])\longrightarrow\underline\bigotimes^+(B[1])$ tel que $F_1=f$ et $F_k=0$ si $k>1$.

L'\'equation de $C_\infty$ morphisme $m^B\circ F=F\circ m^A$ pour un morphisme $F$ de cog\`ebres de Lie, \'ecrite sur les applications $F_n:\underline\bigotimes^n(A[1])\longrightarrow B[1-n]$, s'appelle l'\'equation de $C_\infty$ formalit\'e.\\

\n D'une part, on a
$$
\begin{aligned}
m^B\circ F(\alpha_{[1,n]})&=\sum_{k,r_1+\dots+r_k=n}m^B\big(F_{r_1}(\alpha_{[1,s_1]})\underline{\otimes}\dots\underline{\otimes}F_{r_k}(\alpha_{[s_{k-1}+1,s_k]})
\big)\cr&=\sum_{0<j<k}(-1)^{a_{[1,s_{j-1}]}}F_{r_1}(\alpha_{[1,s_1]})\underline{\otimes}\dots\underline{\otimes}F_{r_j-1}(\alpha_{[s_{j-2}+1,s_{j-1}]})
\underline{\otimes}\cr&m^B\big(F_{r_j}(\alpha_{[s_{j-1}+1,s_{j}]})\underline{\otimes}F_{r_{j+1}}(\alpha_{[s_{j}+1,s_{j+1}]})\big)
\underline{\otimes}F_{r_{j+2}}(\alpha_{[s_{j+1}+1,s_{j+2}]})\underline{\otimes}\dots\underline{\otimes}F_{r_{k}}(\alpha_{[s_{k-1}+1,s_{k}]}).
\end{aligned}
$$

\n D'autre part, on a
$$
\begin{aligned}
F\circ m^A(\alpha_{[1,n]})&=F\big(\sum_{j=1}^{n-1}(-1)^{a_{[1,j-1]}}\alpha_{[1,j-1]}\underline{\otimes}m^A(\alpha_j\underline{\otimes}\alpha_{j+1})\underline
{\otimes}\alpha_{[j+2,n]}\big)\cr&=\sum_{j=1}^{n-1}(-1)^{a_{[1,j-1]}}\sum_{k, r_1+\dots+r_k=n-1}\cr&F_{r_1}(\alpha_{[1,s_1]}\underline{\otimes}\dots\underline{\otimes}F_{r_t}\big(\alpha_{[s_{t-1}+1,j-1]\underline{\otimes}}m^A
(\alpha_j\underline{\otimes}\alpha_{j+1})\underline{\otimes}\alpha_{[j+1,s_t]}\big)\underline{\otimes}\dots\underline{\otimes}F_{r_k}(
\alpha_{[s_{k-1}+1,s_k]}).
\end{aligned}
$$

\

\n En \'ecrivant $(m^B\circ F-F\circ m^A)(\alpha_{[1,n]})$ \`a l'ordre $n-1$, on trouve les termes

$$
\begin{aligned}
&m^B\big(F_{n-1}(\alpha_{[1,n-1]})\underline{\otimes}F_1(\alpha_n)\big)+m^B\big(F_1(\alpha_1)\underline{\otimes}F_{n-1}(\alpha_{[2,n]})\big)\cr&
-\sum_{j=1}^{n-1}(-1)^{a_{[1,j-1]}}F_{n-1}\big(\alpha_{[1,j-1]\underline{\otimes}}m^A
(\alpha_j\underline{\otimes}\alpha_{j+1})\underline{\otimes}\alpha_{[j+2,n]}\big)\cr&=(-1)^{a_{[1,n-1]}+1}F_{n-1}(\alpha_{[1,n-1]})\wedge F_1(\alpha_n)+(-1)^{a_1+1}F_1(\alpha_1)\wedge F_{n-1}(\alpha_{[2,n]})\cr&-\sum_{j=1}^{n-1}(-1)^{a_{[1,j]+1}}F_{n-1}\big(\alpha_{[1,j-1]\underline{\otimes}}
(\alpha_j\wedge\alpha_{j+1})\underline{\otimes}\alpha_{[j+2,n]}\big)\cr&=(d_{Ha}F_{n-1})(\alpha_{[1,n-1]})
\end{aligned}
$$

\

On retrouve l'op\'erateur de cobord de Harrison $d_{Ha}$.\\

Finalement, comme pour les alg\`ebres de Lie, si $A$ est une alg\`ebre commutative et $V$ un $A$ module vu comme un bimodule, l'espace $B=A\oplus\sum_{p>0}V[p]$ muni du produit
$$
(\alpha+\sum u_p)(\beta+\sum v_q)=(\alpha\beta+\sum \alpha v_p+u_p\beta)
$$
est une alg\`ebre commutative et l'application $f:A\longrightarrow B$, d\'efinie par $f(\alpha)=(\alpha,0)$ est un morphisme d'alg\`ebres.\\

Comme pour les alg\`ebres de Lie, un morphisme de cog\`ebres de Lie $F$ tel que $F_1=f+C_1$, $F_k=C_k$ ($k>1$) avec $C_k:\underline\bigotimes^k(A[1])\longrightarrow V[k]$ est un morphisme de $C_\infty$ alg\`ebres sera appel\'e une $C_\infty$ formalit\'e.

Cette formalit\'e est dite triviale s'il existe un morphisme $G$ tel que $C=m^B\circ G+G\circ m^A$, $G$ \'etant de degr\'e -1 et $G=\sum B_p$ avec $B_p:\underline\bigotimes^p(A[1])\longrightarrow V[p]$.

On retrouve ainsi la cohomologie de Harrison des alg\`ebres commutatives, puisque

\vskip 0.3cm
\begin{prop} {\rm ($C_\infty$ formalit\'es et cohomologie de Harrison)}

\

Avec les notations pr\'ec\'edentes, $F$ est une $C_\infty$ formalit\'e si et seulement si
$$
d_{Ha}C_k=0\quad\text{pour tout }k>0.
$$

$F$ est triviale si et seulement si
$$
C_1=0\quad\text{et}\quad C_k=d_{Ha}B_k\quad\text{pour tout }k>1.
$$
\end{prop}

\

\

\section{Alg\`ebre de Lie diff\'erentielle associ\'ee \`a une alg\`ebre de Gerstenhaber}\label{1}

\

\

\subsection{Alg\`ebres de Gerstenhaber}

\

\

Le prototype des alg\`ebres de Gerstenhaber est l'espace $T_{poly}(\mathbb R^d)$ des champs de tenseurs sur $\mathbb R^d$. Cet espace gradu\'e est muni du produit ext\'erieur $\wedge$ et du crochet $[~,~]_S$ de Schouten. Les axiomes usuels d'une alg\`ebre de Gerstenhaber sont donc les suivants.\\

Une alg\`ebre de Gerstenhaber est un espace vectoriel gradu\'e $\mathcal G$ muni d'une multiplication commutative gradu\'ee et associative $\wedge: \mathcal G\otimes\mathcal G\longrightarrow\mathcal G$ de degr\'e $0$ et d'un crochet $[~~,~~]: \mathcal G\otimes\mathcal G\longrightarrow\mathcal G$
de degr\'e $-1$ tel que $\left(\mathcal G[1], [~~,~~]\right)$ soit une alg\`ebre de Lie gradu\'ee et que, pour tout $\alpha$ homog\`ene, l'application $[\alpha, ~.~]$ soit une d\'erivation gradu\'ee pour la multiplication $\wedge$. En notant $|\alpha|$ le degr\'e d'un \'el\'ement homog\`ene $\alpha$ de $\mathcal G$, on a donc:
$$
\begin{aligned}
\alpha\wedge\beta&=(-1)^{|\alpha||\beta|}\beta\wedge\alpha,\\
\alpha\wedge(\beta\wedge\gamma)&=(\alpha\wedge\beta)\wedge\gamma,\\
[\alpha,\beta]&=-(-1)^{(|\alpha|-1)(|\beta|-1)}[\beta,\alpha],\\
0=(-1)^{(|\alpha|-1)(|\gamma|-1)}\big[[\alpha,\beta],\gamma\big]&+(-1)^{(|\beta|-1)(|\alpha|-1)}\big[[\beta,\gamma],\alpha\big]+(-1)^{(|\gamma|-1)(|\beta|-1)}\big[[\gamma,\alpha],\beta\big],\\
[\alpha,\beta\wedge\gamma]&=[\alpha,\beta]\wedge\gamma+(-1)^{|\beta|(|\alpha|-1)}\beta\wedge[\alpha,\gamma]
\end{aligned}
$$
et donc aussi:
$$
[\alpha\wedge\beta,\gamma]=\alpha\wedge[\beta,\gamma]+(-1)^{|\beta|(|\gamma|-1)}[\alpha,\gamma]\wedge\beta.
$$

Remarquons qu'il n'y a pas d'\'equivalent non gradu\'e \`a la structure d'alg\`ebre de Gerstenhaber. En particulier, une alg\`ebre de Poisson gradu\'ee n'est pas une alg\`ebre de Gerstenhaber.\\

La derni\`ere identit\'e ne v\'erifie malheureusement pas la r\`egle des signes de Koszul qui s'\'ecrirait ici:
$$
\begin{aligned}
{[\alpha\wedge\beta,\gamma]}&=[~,~]\big(\wedge(\alpha,\beta),\gamma\big)\\ &=\varepsilon\left(\begin{smallmatrix}1&0&|\alpha|&|\beta|&|\gamma|\\ 0&|\alpha|&1&|\beta|&|\gamma|\end{smallmatrix}\right)\alpha\wedge[\beta,\gamma]+\varepsilon\left(\begin{smallmatrix}1&0&|\alpha|&|\beta|&|\gamma|\\ 0&1&|\alpha|&|\gamma|&|\beta|\end{smallmatrix}\right)[\alpha,\gamma]\wedge\beta\\
&=(-1)^{|\alpha|}\alpha\wedge[\beta,\gamma]+(-1)^{|\beta||\gamma|}[\alpha,\gamma]\wedge\beta.
\end{aligned}
$$

Pour \'eviter ce probl\`eme, on \'ecrit les axiomes des alg\`ebres de Gerstenhaber dans $\mathcal G[1]$, apr\`es un d\'ecalage de degr\'e.\\

On note comme ci-dessus $a,b,\dots$ les degr\'es de $\alpha,\beta,\dots$. Le produit $\wedge$ donne un produit $\mu_2$ :
$$
\mu_2(\alpha,\beta)=(-1)^a\alpha\wedge\beta.
$$
On a vu que $\mu_2$ est anticommutatif et antiassociatif et de degr\'e 1. Le crochet $[~,~]$ est un crochet d'alg\`ebre de Lie gradu\'ee sur $\mathcal G[1]$ (de degr\'e 0.) De plus, l'application $[\alpha,~.~]$ est une d\'erivation gradu\'ee pour la multiplication $\mu_2$ qui v\'erifie bien la r\`egle des signes de Koszul:
$$
[\alpha,\mu_2(\beta,\gamma)]=(-1)^a\mu_2([\alpha,\beta],\gamma)+(-1)^{a(b+1)}\mu_2(\beta,[\alpha,\gamma])
$$
et
$$
[\mu_2(\alpha,\beta),\gamma]=\mu_2(\alpha,[\beta,\gamma])+(-1)^{bc}\mu_2([\alpha,\gamma],\beta).
$$

\


\subsection{La $C_\infty$ alg\`ebre vue comme une alg\`ebre de Lie}

\

\

Comme dans la section pr\'ec\'edente, on peut construire la $C_\infty$ alg\`ebre naturellement associ\'ee \`a l'alg\`ebre commutative $(\mathcal G,\wedge)$. On notera cette cog\`ebre diff\'erentielle:
$$
(\mathcal H,\delta,\mu)=\left(\underline\bigotimes^+(\mathcal G[1]),\delta,\mu\right).
$$
Avec, comme plus haut,
$$
\begin{aligned}
\delta(\alpha_{[1,n]})&=\sum_{0<j<n}\Big(\alpha_{[1,j]}\bigotimes \alpha_{[j+1,n]}-\tau\big(\alpha_{[1,j]}\bigotimes \alpha_{[j+1,n]}\big)\Big),\\
\mu(\alpha_{[1,n]})&=\sum_{0<j<n}(-1)^{a_{[1,j-1]}}\alpha_{[1,j-1]}\underline\otimes\mu_2(\alpha_j,\alpha_{j+1})\underline\otimes\alpha_{[j+2,n]}.
\end{aligned}
$$

Maintenant, il faut prolonger la d\'efinition du crochet de $\mathcal G[1]$ \`a $\mathcal H$. On pose donc:

On prolonge d'abord l'op\'erateur $[~,~]$ sur l'espace $\mathcal G^{\otimes ^p}\wedge\mathcal G^{\otimes ^q}$ en d\'efinissant le crochet de
$$
\alpha=\alpha_{\{1,\dots,p\}}=\alpha_1\otimes\dots\otimes\alpha_p\quad\text{ et}\quad \beta=\alpha_{\{p+1,\dots,p+q\}}=\alpha_{p+1}\otimes\dots\otimes\alpha_{p+q}
$$
par:
$$
\begin{aligned}
{[\alpha,\beta]}&=\hskip-0.5cm\sum_{\begin{smallmatrix}\sigma\in Bat(p,q)\\ \sigma^{-1}(k)\leq p<\sigma^{-1}(k+1)\end{smallmatrix}}\hskip-0.5cm
\varepsilon\left(\begin{smallmatrix}a_1~~\dots~~a_{p+q}\\ a_{\sigma^{-1}(1)}\dots a_{\sigma^{-1}(p+q)}\end{smallmatrix}\right) \alpha_{\sigma^{-1}(1)}\otimes\dots\otimes[\alpha_{\sigma^{-1}(k)},\alpha_{\sigma^{-1}(k+1)}]\otimes\dots\otimes\alpha_{\sigma^{-1}(p+q)}\\ &=\hskip-0.5cm\sum_{\begin{smallmatrix}\sigma\in Bat(p,q)\\ \sigma^{-1}(k)\leq p<\sigma^{-1}(k+1)\end{smallmatrix}}\hskip-0.5cm[\sigma.(\alpha\otimes \beta)]_k.
\end{aligned}
$$
Maintenant, on peut passer au quotient par les battements puisque:

\begin{lem}\

Soient $\alpha\in\mathcal G[1]^{\otimes p}$, $\beta\in\mathcal G[1]^{\otimes q}$ et $\gamma\in\mathcal G[1]^{\otimes r}$. On a alors:
$$
[bat_{p,q}(\alpha,\beta),\gamma]=bat_{p,q+r-1}(\alpha,[\beta,\gamma])+(-1)^{ab}bat_{q,p+r-1}(\beta,[\alpha,\gamma]).
$$

\end{lem}

\vskip 0.3cm
\noindent
{\bf Preuve}

Notons $\alpha=\alpha_1\otimes\dots\otimes\alpha_p$, $\beta=\alpha_{p+1}\otimes\dots\otimes\alpha_{p+q}$ et $\gamma=\alpha_{p+q+1}\otimes\dots\otimes\alpha_{p+q+r}$. On a

\begin{align*}
[bat_{p,q}(\alpha,\beta),\gamma]&=\sum_{\begin{smallmatrix}\sigma_2\in Bat(p+q,r)\\ \sigma_2^{-1}(k)\leq p+q<\sigma_2^{-1}(k+1) \end{smallmatrix}} [\sigma_2.(bat_{p,q}(\alpha,\beta)\otimes\gamma)]_k\\
&=\sum_{\begin{smallmatrix}\sigma_2\in Bat(p+q,r),~\sigma_1\in Bat(p,q)\\ \sigma_2^{-1}(k)\leq p+q<\sigma_2^{-1}(k+1) \end{smallmatrix}} \hskip-0.5cm [\sigma_2\circ(\sigma_1\otimes id_{\{p+q+1,\dots,p+q+r\}}).(\alpha\otimes\beta\otimes\gamma)]_k\\
&=\hskip-0.5cm\sum_{\begin{smallmatrix}\rho\in Bat(p,q,r)\\ \rho^{-1}(k)\leq p+q<\rho^{-1}(k+1)\end{smallmatrix}}\hskip-0.5cm [\rho.(\alpha\otimes \beta\otimes\gamma)]_k\\
&=\hskip-0.5cm\sum_{\begin{smallmatrix}\rho\in Bat(p,q,r)\\ p<\rho^{-1}(k)\leq p+q<\rho^{-1}(k+1)\end{smallmatrix}}\hskip-0.5cm [\rho.(\alpha\otimes \beta\otimes\gamma)]_k+\hskip-0.5cm\sum_{\begin{smallmatrix}\rho\in Bat(p,q,r)\\ \rho^{-1}(k)\leq p< p+q<\rho^{-1}(k+1)\end{smallmatrix}}\hskip-0.5cm [\rho.(\alpha\otimes\beta\otimes\gamma)]_k\\
&=(I)+(II).
\end{align*}

D'autre part, on a aussi
$$
{[\beta,\gamma]}=\hskip-0.5cm\sum_{\begin{smallmatrix}\sigma_1\in Bat(q,r)\\ p<\sigma_1^{-1}(k)\leq p+q<\sigma_1^{-1}(k+1)\end{smallmatrix}} \hskip-0.5cm[\sigma_1.(\beta\otimes\gamma)]_k.
$$
Donc
$$
\begin{aligned}
\alpha\otimes[\beta,\gamma]&=\hskip-0.5cm\sum_{\begin{smallmatrix}(id\otimes\sigma_1)\in Bat(p,q,r)\\ p<\sigma_1^{-1}(k)\leq p+q<\sigma_1^{-1}(k+1) \end{smallmatrix}}\hskip-0.5cm [(id\otimes\sigma_1).(\alpha\otimes\beta\otimes\gamma)]_k\\
&=\hskip-0.5cm\sum_{\begin{smallmatrix}(id\otimes\sigma_1)\in Bat(p,q,r)\\ p<\sigma_1^{-1}(k)\leq p+q<\sigma_1^{-1}(k+1) \end{smallmatrix}}\hskip-0.5cm\varepsilon\left(\begin{smallmatrix}a_1\hskip 0.5cm\dots\hskip0.5cm a_{p+q+r}\\ a_{(id\otimes\sigma_1)^{-1}(1)}\dots a_{(id\otimes\sigma_1)^{-1}(p+q+r)}\end{smallmatrix}\right)\beta_1^{\sigma_1}\otimes \dots\otimes \beta_{p+q+r-1}^{\sigma_1},
\end{aligned}
$$
o\`u on a pos\'e
$$
\beta_j^{\sigma_1}=\left\{
  \begin{array}{ll}
   \alpha_{(id\otimes\sigma_1)^{-1}(j)}, & \hbox{ si $1\leq j<k$,} \cr
   &\cr
   [\alpha_{(id\otimes\sigma_1)^{-1}(k)},\alpha_{(id\otimes\sigma_1)^{-1}(k+1)}], & \hbox{ si $j=k$,} \cr
   &\cr
   \alpha_{(id\otimes\sigma_1)^{-1}(j+1)}, & \hbox{si $k<j\leq p+q+r-1$.}
  \end{array}
\right.
$$
Par suite, on a
$$
\begin{aligned}
bat_{p,q+r-1}(\alpha,[\beta,\gamma])=&\hskip-0.3cm\sum_{\sigma_2\in Bat(p,q+r-1)}\hskip-0.5cm\sum_{\begin{smallmatrix}\sigma_1\in Ba(q,r)\\ p<\sigma_1^{-1}(k)\leq p+q< \sigma_1^{-1}(k+1) \end{smallmatrix}}\hskip-0.5cm\varepsilon\left(\begin{smallmatrix}\beta_1~~\dots~~\beta_{p+q+r-1}\\ \beta_{\sigma_2^{-1}(1)}\dots\beta_{\sigma_2^{-1}(p+q+r-1)}\end{smallmatrix}\right)\\&\varepsilon\left(\begin{smallmatrix}a_1\hskip 0.5cm\dots\hskip0.5cm a_{p+q+r}\\ a_{(id\otimes\sigma_1)^{-1}(1)}\dots a_{(id\otimes\sigma_1)^{-1}(p+q+r)}\end{smallmatrix}\right)\beta_{\sigma_2^{-1}(1)}^{\sigma_1}\otimes \dots\otimes \beta_{\sigma_2^{-1}(p+q+r-1)}^{\sigma_1}.
\end{aligned}
$$
En posant $k'=\sigma_2^{-1}(k)$, on voit que l'application $(\sigma_1,\sigma_2,k)\mapsto(\rho,k')$ est une bijection entre les ensembles
$$
\left\{(\sigma_1,\sigma_2,k)\in Bat(q,r)\times Bat(p,q+r-1)\times[1,p+q+r-1],~~p<\sigma_1^{-1}(k)\leq p+q< \sigma_1^{-1}(k+1)\right\}
$$
et
$$
\left\{(\rho,k')\in Bat(p,q,r)\times[1,p+q+r-1],~~p<\rho^{-1}(k')\leq p+q<\rho^{-1}(k'+1)\right\}.
$$
Donc
$$
\begin{aligned}
bat_{p,q+r-1}(\alpha,[\beta,\gamma])&=\hskip-0.5cm\sum_{\begin{smallmatrix}\rho\in Bat(p,q,r)\\ p<\rho^{-1}(k')\leq p+q< \rho^{-1}(k'+1)\end{smallmatrix}} \hskip-0.5cm [\varepsilon\left(\begin{smallmatrix}a_1~~\dots~~a_{p+q+r}\\ a_{\rho^{-1}(1)}\dots a_{\rho^{-1}(p+q+r)}\end{smallmatrix}\right)\alpha_{\rho^{-1}(1)}\otimes\dots\otimes\alpha_{\rho^{-1}(p+q+r)}]_{k'}\\
&=\hskip-0.5cm\sum_{\begin{smallmatrix}\rho\in Bat(p,q,r)\\ p<\rho^{-1}(k')\leq p+q< \rho^{-1}(k'+1)\end{smallmatrix}} \hskip-0.5cm [\rho.(\alpha\otimes\beta\otimes\gamma)]_{k'}\\
&=(I).
\end{aligned}
$$

On montre de m\^eme que
$$
(-1)^{ab}bat_{q,p+r-1}(\beta,[\alpha,\gamma])=\hskip-0.5cm\sum_{\begin{smallmatrix}\rho\in Bat(p,q,r)\\ \rho^{-1}(k')\leq p<p+q< \rho^{-1}(k'+1)\end{smallmatrix}} \hskip-0.5cm [\rho.(\alpha\otimes\beta\otimes\gamma)]_k=(II).
$$
D'o\`u le lemme. \hfill$\square$

\

Le lemme nous permet de d\'efinir $[~,~]$ par la m\^eme expression sur l'espace $\underline\bigotimes^p(\mathcal G[1])\wedge
\underline\bigotimes^q(\mathcal G[1])$. On obtient un crochet sur $\mathcal H$ qui v\'erifie les identit\'es de Jacobi et de Leibniz.\\

\begin{thm} {\rm ($\mathcal H$ est une alg\`ebre de Lie diff\'erentielle gradu\'ee)}

\

L'espace $\mathcal H$, muni du crochet $[~,~]$ et de l'op\'erateur $\mu$ est une alg\`ebre de Lie diff\'erentielle gradu\'ee: Pour tout $\alpha$, $\beta$ et $\gamma$ de $\mathcal H$, on a:

\begin{itemize}
\item[(i)] \quad$[\alpha,\beta]=-(-1)^{ab}[\beta,\alpha]$,

\item[(ii)] \quad$(-1)^{ac}\left[[\alpha,\beta],\gamma\right]+(-1)^{ba}\left[[\beta,\gamma],\alpha\right]+(-1)^{cb}\left[[\gamma,\alpha],\beta\right]=0$,

\item[(iii)] \quad$\mu\left([\alpha,\beta]\right)=\left[\mu(\alpha),\beta\right]+(-1)^a\left[\alpha,\mu(\beta)\right]$.
\end{itemize}
\end{thm}

\vskip 0.3cm
\noindent
{\bf Preuve}

\noindent
(i)  On sait que
$$
\begin{aligned}
{[\alpha_{[1,p]},\alpha_{[p+1,p+q]}]=}&\hskip-0.5cm\sum_{\begin{smallmatrix}\sigma\in Bat(p,q)\\ \sigma^{-1}(k)\leq p<\sigma^{-1}(k+1)\end{smallmatrix}}
\hskip-0.5cm \varepsilon\left(\begin{smallmatrix}a_1~\dots~a_{p+q}\\ a_{\sigma^{-1}(1)}\dots a_{\sigma^{-1}(p+q)}\end{smallmatrix}\right)\\
&\hskip 2cm\alpha_{\sigma^{-1}(1)}\underline\otimes\dots\underline\otimes[\alpha_{\sigma^{-1}(k)},\alpha_{\sigma^{-1}(k+1)}]\underline
\otimes\dots\underline\otimes\alpha_{\sigma^{-1}(p+q)}.
\end{aligned}
$$
Fixons un couple $(\sigma,k)$ dans cette somme (tel que $\sigma^{-1}(k)\leq p<\sigma^{-1}(k+1)$.) On d\'efinit trois permutations $\tau$, $\rho$ et $\nu$ de $S_{p+q}$ par:
$$
\tau(j)=\left\{
              \begin{array}{ll}
               j+p, & \hbox{si $1\leq j\leq q$} \\
                j-q, & \hbox{si $q< j\leq q+p$.}
              \end{array}
            \right.
,\qquad\rho=\sigma\circ\tau\quad\hbox{et}\quad\nu^{-1}(i)=\left\{
              \begin{array}{ll}
               \rho^{-1}(i), & \hbox{si }~i\notin\{k,k+1\} \\
               \rho^{-1}(k+1), & \hbox{si }~i=k,\\
               \rho^{-1}(k), & \hbox{si }~i=k+1.
              \end{array}
            \right.
$$

On v\'erifie imm\'ediatement que $\nu$ appartient \`a $Bat(q,p)$ et que $\nu^{-1}(k)\leq q<\nu^{-1}(k+1)$. De plus l'application $(\sigma,k)\mapsto(\nu,k)$ est une bijection sur les ensembles correspondants.

Posons maintenant $\beta_j=\alpha_{\tau(j)}$. On a:
$$
\varepsilon\left(\begin{smallmatrix}b_1~~\dots~~b_{p+q}\\ b_{\nu^{-1}(1)}\dots b_{\nu^{-1}(p+q)}\end{smallmatrix}\right)=
\varepsilon\left(\begin{smallmatrix}a_1~~\dots~~a_{p+q}\\ a_{\sigma^{-1}(1)}\dots a_{\sigma^{-1}(p+q)}\end{smallmatrix}\right)(-1)^{a_{\sigma^{-1}(k)}a_{\sigma^{-1}(k+1)}}(-1)^{a_{[1,p]}a_{[p+1,p+q]}}.
$$
Donc :
$$
\begin{aligned}
{[\alpha_{[1,p]},\alpha_{[p+1,p+q]}]=}&(-1)^{a_{[1,p]}a_{[p+1,p+q]}}\hskip-0.5cm\sum_{\begin{smallmatrix}\nu\in Bat(q,p)\\ \nu^{-1}(k)\leq q<\nu^{-1}(k+1)\end{smallmatrix}}
\hskip-0.5cm \varepsilon\left(\begin{smallmatrix}b_1~\dots~b_{p+q}\\ b_{\nu^{-1}(1)}\dots b_{\nu^{-1}(p+q)}\end{smallmatrix}\right)\\
&\hskip 0.5cm(-1)^{b_{\nu^{-1}(k)}b_{\nu^{-1}(k+1)}}\beta_{\nu^{-1}(1)}\underline\otimes\dots\underline\otimes[\beta_{\nu^{-1}(k+1)},\beta_{\nu^{-1}(k)}]\underline\otimes\dots\underline\otimes\beta_{\nu^{-1}(p+q)}\\
=&(-1)^{a_{[1,p]}a_{[p+1,p+q]}}[\alpha_{[p+1,p+q]},\alpha_{[1,p]}].
\end{aligned}
$$
D'o\`u le r\'esultat.

\

\noindent
(ii) Soient $\alpha=\alpha_1\underline\otimes\dots\underline\otimes\alpha_p$, $\beta=\beta_1\underline\otimes\dots\underline\otimes\beta_q$ et $\gamma=\gamma_1\underline\otimes\dots\underline\otimes\gamma_r$. Pour all\'eger les notations, pour toute permutation $\rho\in S_{p+q+r}$, notons $\varepsilon(\rho)$ le signe:
$$
\varepsilon(\rho)=\varepsilon\left(\begin{smallmatrix}x_1~~\dots~~x_{p+q+r}\\ x_{\rho^{-1}(1)}\dots x_{\rho^{-1}(p+q+r)}\end{smallmatrix}\right),
$$
si
$$
\xi_i=\left\{\begin{aligned}\alpha_i&\quad\text{si }~1\leq i\leq p\\
\beta_{i-p}&\quad\text{si }~p+1\leq i\leq p+q\\
\gamma_{i-p-q}&\quad\text{si }~p+q+1\leq i\leq p+q+r
\end{aligned}
\right.
$$

Gr\`ace \`a l'associativit\'e du produit battement, on a
$$
\begin{aligned}
(-1)^{ac}bat_{p,q,r}(\alpha\otimes\beta\otimes\gamma)&=(-1)^{ba}bat_{q,r,p}(\beta\otimes\gamma\otimes\alpha)\\
&=(-1)^{cb}bat_{r,p,q}(\gamma\otimes\alpha\otimes\beta).
\end{aligned}
$$
Plus pr\'ecis\'ement, il y a deux bijections canoniques:
$$
\begin{matrix}
Bat(p,q,r)&\longrightarrow&Bat(q,r,p)& \hskip0.5cm \hbox{et} \hskip0.5cm&Bat(p,q,r)&\longrightarrow&Bat(r,p,q)\\
\rho_1&\mapsto&\rho_2=\rho_1\circ\tau&&\rho_1&\longmapsto&\rho_3=\rho_1\circ\tau'.
\end{matrix}
$$
Avec
$$
\tau(i)=\left\{
              \begin{array}{ll}
               i+p, & \hbox{si $1\leq i\leq q+r$} \\
                i-(q+r), & \hbox{si $q+r< i\leq q+r+p$.}
              \end{array}
            \right.\quad\text{et}\quad\tau'(i)=\left\{
              \begin{array}{ll}
               i+p+q, & \hbox{si $1\leq i\leq r$} \\
                i-r, & \hbox{si $r< i\leq q+r+p$.}
              \end{array}
            \right.
$$
On a alors
$$
(-1)^{ac}\rho_1.(\alpha\otimes\beta\otimes\gamma)=(-1)^{ba}\rho_2.(\beta\otimes\gamma\otimes\alpha)=(-1)^{cb}\rho_3.(\gamma\otimes\alpha\otimes\beta).
$$
En \'ecrivant $(-1)^{ac}[[\alpha,\beta],\gamma]$, on trouve des termes de la forme:
$$
\begin{aligned}
&(1.1):\ (-1)^{ac}\varepsilon(\rho_{11})\xi_{i_1}\underline\otimes\dots\underline\otimes[\alpha_i,\beta_j] \underline\otimes\dots\underline\otimes [\beta_k,\gamma_l]\underline\otimes \dots\underline\otimes\xi_{i_{p+q+r}}\\
&(1.2):\ (-1)^{ac}\varepsilon(\rho_{12})\xi_{j_1}\underline\otimes\dots\underline\otimes[\beta_j,\gamma_l]\underline\otimes \dots\underline\otimes [\alpha_i,\beta_k]\underline\otimes \dots\underline\otimes \xi_{j_{p+q+r}}\\
&(1.3):\ (-1)^{ac}\varepsilon(\rho_{13})\xi_{k_1}\underline\otimes\dots\underline\otimes[\alpha_i,\beta_j]\underline\otimes \dots\underline\otimes[\alpha_k,\gamma_l]\underline\otimes \dots\underline\otimes \xi_{k_{p+q+r}}\cr
&(1.4):\ (-1)^{ac}\varepsilon(\rho_{14})\xi_{l_1}\underline\otimes\dots\underline\otimes[\alpha_i,\gamma_l]\underline\otimes \dots\underline\otimes[\alpha_k,\beta_j]\underline\otimes \dots\underline\otimes \xi_{l_{p+q+r}}\cr
&(1.5):\ (-1)^{ac}\varepsilon(\rho_{15})\xi_{s_1}\underline\otimes\dots\underline\otimes[[\alpha_i,\beta_j],\gamma_k]\underline\otimes \dots\underline\otimes \xi_{s_{p+q+r}}.
\end{aligned}
$$
En \'ecrivant $(-1)^{ba}[[\beta,\gamma],\alpha]$, on trouve des termes de la forme:
$$
\begin{aligned}
&(2.1):\ (-1)^{ba}\varepsilon(\rho_{21})\xi_{i_1}\underline\otimes\dots\underline\otimes[\beta_j,\alpha_i]\underline\otimes \dots\underline\otimes [\beta_k,\gamma_l]\underline\otimes \dots\underline\otimes \xi_{i_{p+q+r}}\\
&(2.2):\ (-1)^{ba}\varepsilon(\rho_{22})\xi_{j_1}\underline\otimes\dots\underline\otimes[\beta_j,\gamma_l]\underline\otimes \dots\underline\otimes[\beta_k,\alpha_i]\underline\otimes \dots\underline\otimes \xi_{j_{p+q+r}}\\
&(2.3):\ (-1)^{ba}\varepsilon(\rho_{23})\xi_{t_1}\underline\otimes\dots\underline\otimes[\beta_j,\gamma_k]\underline\otimes \dots\underline\otimes [\gamma_l,\alpha_i]\underline\otimes \dots\underline\otimes \xi_{t_{p+q+r}}\\
&(2.4):\ (-1)^{ba}\varepsilon(\rho_{24})\xi_{r_1}\underline\otimes\dots\underline\otimes[\gamma_k,\alpha_i]\underline\otimes \dots\underline\otimes[\beta_j,\gamma_l]\underline\otimes  \dots\underline\otimes \xi_{r_{p+q+r}}\\
&(2.5):\ (-1)^{ba}\varepsilon(\rho_{25})\xi_{s_1}\underline\otimes\dots \underline\otimes[[\beta_j,\gamma_k],\alpha_i]\underline\otimes \dots\underline\otimes \xi_{s_{p+q+r}}.
\end{aligned}
$$
En \'ecrivant $(-1)^{cb}[[\gamma,\alpha],\beta]$, on trouve des termes de la forme:
$$
\begin{aligned}
&(3.1):\ (-1)^{cb}\varepsilon(\rho_{31})\xi_{k_1}\underline\otimes\dots\underline\otimes[\alpha_i,\beta_j]\underline\otimes \dots\underline\otimes [\gamma_l,\alpha_k]\underline\otimes \dots\underline\otimes \xi_{k_{p+q+r}}\\
&(3.2):\ (-1)^{cb}\varepsilon(\rho_{32})\xi_{l_1}\underline\otimes\dots\underline\otimes[\gamma_l,\alpha_i]\underline\otimes \dots\underline\otimes[\alpha_k,\beta_j]\underline\otimes \dots\underline\otimes \xi_{l_{p+q+r}}\\
&(3.3):\ (-1)^{cb}\varepsilon(\rho_{33})\xi_{r_1}\underline\otimes\dots\underline\otimes[\gamma_k,\alpha_i]\underline\otimes \dots\underline\otimes [\gamma_l,\beta_j]\underline\otimes \dots\underline\otimes \xi_{r_{p+q+r}}\\
&(3.4):\ (-1)^{cb}\varepsilon(\rho_{34})\xi_{t_1}\underline\otimes\dots\underline\otimes[\gamma_k,\beta_j]\underline\otimes \dots\underline\otimes[\gamma_l,\alpha_i]\underline\otimes \dots\underline\otimes \xi_{t_{p+q+r}}\\
&(3.5):\ (-1)^{cb}\varepsilon(\rho_{35})\xi_{s_1}\underline\otimes\dots\underline\otimes[[\gamma_k,\alpha_i],\beta_j] \underline\otimes \dots\underline\otimes \xi_{s_{p+q+r}}.
\end{aligned}
$$

Consid\`erons les termes $(1.1)$ et $(2,1)$, on voit que $\rho_{21}$ est d\'efinie par $\rho_{21}^{-1}=(a_i,b_j)\circ(\rho_{11}\circ\tau)^{-1}$. On en d\'eduit que la somme de ces deux termes s'annule.

De m\^eme, on v\'erifie que: $(1.2)+(2.2)=0$, $(1.3)+(3.1)=0$, $(1.4)+(3.2)=0$, $(2.3)+(3.4)=0$ et $(2.4)+(3.3)=0$.

D'autre part, par construction
$$
\varepsilon(\rho_{25})=(-1)^{a_i(b_j+c_k)}(-1)^{a(b+c)}\varepsilon(\rho_{15})\quad\text{et }\quad \varepsilon(\rho_{35})=(-1)^{c_k(a_i+b_j)} (-1)^{c(a+b)}\varepsilon(\rho_{15}).
$$
Alors,
$$
\begin{aligned}
&(1.5)+(2.5)+(3.5)=\\
&=\varepsilon(\rho_{15})(-1)^{ac}\xi_{s_1}\underline\otimes\dots\underline\otimes\\
&\hskip 1cm\underline\otimes\Big([[\alpha_i,\beta_j],\gamma_k]+(-1)^{a_i(b_j+c_k)}[[\beta_j,\gamma_k],\alpha_i]
+(-1)^{c_k(a_i+b_j)}[[\gamma_k,\alpha_i],\beta_j]\Big)\underline\otimes \dots\underline\otimes \xi_{s_{p+q+r}}\cr
&=0.
\end{aligned}
$$

\

\noindent
(iii) On conserve la notation
$$
\varepsilon(\sigma)=\varepsilon\left(\begin{smallmatrix}a_1~~\dots~~a_{p+q}\\ a_{\sigma^{-1}(1)}\dots a_{\sigma^{-1}(p+q)}\end{smallmatrix}\right).
$$
Alors
$$
\begin{aligned}
\mu([\alpha_{[1,p]}&,\alpha_{[p+1,p+q]}])=\hskip-0.3cm\sum_{\begin{smallmatrix}\sigma\in Bat(p,q)\\ i<k-1 ; \ \sigma^{-1}(k)\leq p<\sigma^{-1}(k+1) \end{smallmatrix}}\hskip-0.3cm \varepsilon(\sigma)(-1)^{\sum_{r<i}a_{\sigma^{-1}(r)}}\\
&\alpha_{\sigma^{-1}(1)}\underline\otimes\dots\underline\otimes \mu\big(\alpha_{\sigma^{-1}(i)},\alpha_{\sigma^{-1}(i+1)}\big)\underline\otimes\dots\underline\otimes [\alpha_{\sigma^{-1}(k)},\alpha_{\sigma^{-1}(k+1)}]
\underline\otimes\dots\underline\otimes\alpha_{\sigma^{-1}(p+q)}\\
&+\hskip-0.3cm\sum_{\begin{smallmatrix}\sigma\in Bat(p,q)\\ i>k+1 ; \ \sigma^{-1}(k)\leq p<\sigma^{-1}(k+1)\end{smallmatrix}}\hskip-0.3cm
\varepsilon(\sigma)(-1)^{\sum_{r<i}a_{\sigma^{-1}(r)}}\\
&\alpha_{\sigma^{-1}(1)}\underline\otimes\dots\underline\otimes [\alpha_{\sigma^{-1}(k)},\alpha_{\sigma^{-1}(k+1)} \underline\otimes\dots \underline\otimes \mu\big(\alpha_{\sigma^{-1}(i)},\alpha_{\sigma^{-1}(i+1)}\big)\underline\otimes \dots\underline\otimes\alpha_{\sigma^{-1}(p+q)}\\
&+\hskip-0.3cm\sum_{\begin{smallmatrix}\sigma\in Bat(p,q)\\ k; \ \sigma^{-1}(k)\leq p<\sigma^{-1}(k+1)\end{smallmatrix}}\hskip-0.3cm
\varepsilon(\sigma)(-1)^{\sum_{r<k-1}a_{\sigma^{-1}(r)}}\\
&\alpha_{\sigma^{-1}(1)}\underline\otimes\dots\underline\otimes\mu\big(\alpha_{\sigma^{-1}(k-1)},[\alpha_{\sigma^{-1}(k)},\alpha_{\sigma^{-1}(k+1)}]\big)\underline\otimes\dots\underline\otimes \alpha_{\sigma^{-1}(p+q)}\\
&+\hskip-0.3cm\sum_{\begin{smallmatrix}\sigma\in Bat(p,q)\\ k; \ \sigma^{-1}(k)\leq p<\sigma^{-1}(k+1)\end{smallmatrix}}\hskip-0.3cm
\varepsilon(\sigma)(-1)^{\sum_{r<k}a_{\sigma^{-1}(r)}}\\
&\alpha_{\sigma^{-1}(1)}\underline\otimes\dots\underline\otimes\mu\big([\alpha_{\sigma^{-1}(k)},\alpha_{\sigma^{-1}(k+1)}],\alpha_{\sigma^{-1}(k+2)}\big)\underline\otimes\dots\underline\otimes\alpha_{\sigma^{-1}(p+q)}\\
&\hskip-0.5cm=(I)+(II)+(III)+(IV).
\end{aligned}
$$

\

Dans la somme $(I)$ (resp. $(II)$), on distingue quatre cas:
\begin{itemize}
\item[1)] si $\{\sigma^{-1}(i),\sigma^{-1}(i+1)\}\subset\{1,\dots,p\}$: on note $(I_1)$ (resp. $(II_1)$) la restriction de $(I)$ (resp. $(II)$)
\`a ce cas.
\item[2)] si $\{\sigma^{-1}(i),\sigma^{-1}(i+1)\}\subset\{p+1,\dots,p+q\}$: on note $(I_2)$ (resp. $(II_2)$) la restriction de $(I)$ (resp. $(II)$) \`a ce cas.
\item[3)] si $\sigma^{-1}(i)\leq p<\sigma^{-1}(i+1)$: on note $(I_3)$ (resp. $(II_3)$) la restriction de $(I)$ (resp. $(II)$) \`a ce cas.
\item[4)] si $\sigma^{-1}(i+1)\leq p<\sigma^{-1}(i)$: on note $(I_4)$ (resp. $(II_4)$) la restriction de $(I)$ (resp. $(II)$) \`a ce cas.
\end{itemize}

On v\'erifie directement que $(I_3)+(I_4)=0$ et $(II_3)+(II_4)=0$.

Dans la somme $(III)$, on distingue deux cas:
\begin{itemize}
\item[1)] si $\sigma^{-1}(k-1)\in\{1,\dots,p\}$: on note $(III_1)$ la restriction de $(III)$ \`a ce cas.
\item[2)] si $\sigma^{-1}(k-1)\in\{p+1,\dots,p+q\}$: on note $(III_2)$ la restriction de $(III)$ \`a ce cas.
\end{itemize}

Dans la somme $(IV)$, on distingue deux cas:
\begin{itemize}
\item[1)] si $\sigma^{-1}(k+2)\in\{1,\dots,p\}$: on note $(IV_1)$ la restriction de $(IV)$ \`a ce cas.
\item[2)] si $\sigma^{-1}(k+2)\in\{p+1,\dots,p+q\}$: on note $(IV_2)$ la restriction de $(IV)$ \`a ce cas.
\end{itemize}

Donc $\mu([X,Y])$ s'\'ecrit:
$$
\mu([\alpha_{[1,p]},\alpha_{[p+1,p+q]}])=(I_1)+(I_2)+(II_1)+(II_2)+(III_1)+(III_2)+(IV_1)+(IV_2).
$$

D'autre part, on sait que
$$
\mu(bat_{p,q}(\alpha_{[1,p]},\beta_{[p+1,p+q]}))=bat_{p-1,q}(\mu(\alpha_{[1,p]}),\alpha_{[p+1,p+q]})+(-1)^{a_{[1,p]}}bat_{p,q-1}(\alpha_{[1,p]},\mu(\alpha_{[p+1,p+q]})).
$$
Donc
$$
\begin{aligned}
&[\mu(\alpha_{[1,p]}),\alpha_{[p+1,p+q]}]=\hskip-0.3cm\sum_{\begin{smallmatrix}\sigma\in Bat(p,q)\\ i<k-1; \sigma^{-1}(k)\leq p<\sigma^{-1}(k+1)\cr\{\sigma^{-1}(i),\sigma^{-1}(i+1)\}\subset\{1,\dots,p\}\end{smallmatrix}}\hskip-0.3cm
\varepsilon(\sigma)(-1)^{\sum_{r<i}a_{\sigma^{-1}(r)}}\\
&\alpha_{\sigma^{-1}(1)}\underline\otimes\dots\underline\otimes \mu\big(\alpha_{\sigma^{-1}(i)},\alpha_{\sigma^{-1}(i+1)}\big)\underline\otimes\dots\underline\otimes [\alpha_{\sigma^{-1}(k)},
\alpha_{\sigma^{-1}(k+1)}]\underline\otimes\dots\underline\otimes\alpha_{\sigma^{-1}(p+q)}\\
&+\hskip-0.3cm\sum_{\begin{smallmatrix}\sigma\in Bat(p,q)\\ i>k+1;
\sigma^{-1}(k)\leq p<\sigma^{-1}(k+1)\\ \{\sigma^{-1}(i),\sigma^{-1}(i+1)\}\subset\{1,\dots,p\}\end{smallmatrix}}\hskip-0.3cm
\varepsilon(\sigma)(-1)^{\sum_{r<i}a_{\sigma^{-1}(r)}}\\
&\alpha_{\sigma^{-1}(1)}\underline\otimes\dots\otimes [\alpha_{\sigma^{-1}(k)},\alpha_{\sigma^{-1}(k+1)}\underline\otimes\dots\underline\otimes \mu\big(\alpha_{\sigma^{-1}(i)},\alpha_{\sigma^{-1}(i+1)}\big)\underline\otimes\dots\underline\otimes\alpha_{\sigma^{-1}(p+q)}\\
&+\hskip-0.3cm\sum_{\begin{smallmatrix}\sigma\in Bat(p,q)\\ k ; \sigma^{-1}(k)<\sigma^{-1}(k+1)\leq p<\sigma^{-1}(k+2)\end{smallmatrix}}\hskip-0.3cm
\varepsilon(\sigma)(-1)^{\sum_{r<k}a_{\sigma^{-1}(r)}}\\
&\alpha_{\sigma^{-1}(1)}\underline\otimes\dots\underline\otimes[\mu\big(\alpha_{\sigma^{-1}(k)},\alpha_{\sigma^{-1}(k+1)}\big),\alpha_{\sigma^{-1}(k+2)}]\underline\otimes\dots\underline\otimes\alpha_{\sigma^{-1}(p+q)}\\
&\hskip-0.5cm=(I')+(II')+(III').
\end{aligned}
$$
On v\'erifie que $(I')=(I_1)$, $(II')=(II_1)$. De plus, en appliquant l'identit\'e de Leibniz pour $[\mu\big(\alpha_{\sigma^{-1}(k)},\alpha_{\sigma^{-1}(k+1)}\big),\alpha_{\sigma^{-1}(k+2)}]$, on obtient,
$$
\begin{aligned}
(III')&=\hskip-0.3cm\sum_{\begin{smallmatrix}\sigma\in Bat(p,q)\\ k ;
\sigma^{-1}(k)<\sigma^{-1}(k+1)\leq p<\sigma^{-1}(k+2)\end{smallmatrix}}\hskip-0.3cm
\varepsilon(\sigma)(-1)^{\sum_{r<k}a_{\sigma^{-1}(r)}}\\
&\alpha_{\sigma^{-1}(1)}\underline\otimes\dots\underline\otimes \mu\big(\alpha_{\sigma^{-1}(k)},[\alpha_{\sigma^{-1}(k+1)}, \alpha_{\sigma^{-1}(k+2)}]\big)\underline\otimes\dots\underline\otimes\alpha_{\sigma^{-1}(p+q)}\\
&+\hskip-0.3cm\sum_{\begin{smallmatrix}\sigma\in Bat(p,q)\\ k ;\sigma^{-1}(k)<\sigma^{-1}(k+1)\leq p<\sigma^{-1}(k+2)\end{smallmatrix}}\hskip-0.3cm
\varepsilon(\sigma)(-1)^{\sum_{r<k}a_{\sigma^{-1}(r)}}(-1)^{a_{\sigma^{-1}(k+1)}a_{\sigma^{-1}(k+2)}}\\
&\alpha_{\sigma^{-1}(1)}\underline\otimes\dots\underline\otimes \mu\big([\alpha_{\sigma^{-1}(k)},\alpha_{\sigma^{-1}(k+2)}],\alpha_{\sigma^{-1}(k+1)}\big)
\underline\otimes\dots\underline\otimes\alpha_{\sigma^{-1}(p+q)}\\
&\hskip-0.5cm=(III_1)+(IV_1).
\end{aligned}
$$

De m\^eme, on v\'erifie que
$$
\begin{aligned}
(-1)^{a_{[1,p]}}[\alpha_{[1,p]}&,\mu(\alpha_{[p+1,p+q]})]=\hskip-0.3cm\sum_{\begin{smallmatrix}\sigma\in Bat(p,q)\\ i<k-1; \sigma^{-1}(k)\leq p<\sigma^{-1}(k+1)\\ \{\sigma^{-1}(i),\sigma^{-1}(i+1)\}\subset\{p+1,\dots,p+q\}\end{smallmatrix}}\hskip-0.3cm
\varepsilon(\sigma)(-1)^{\sum_{r<i}a_{\sigma^{-1}(r)}}\\
&\alpha_{\sigma^{-1}(1)}\underline\otimes\dots\underline\otimes \mu\big(\alpha_{\sigma^{-1}(i)},\alpha_{\sigma^{-1}(i+1)}\big)\underline\otimes\dots\underline\otimes [\alpha_{\sigma^{-1}(k)},
\alpha_{\sigma^{-1}(k+1)}]\underline\otimes\dots\underline\otimes\alpha_{\sigma^{-1}(p+q)}\\
&+\hskip-0.3cm\sum_{\begin{smallmatrix}\sigma\in Bat(p,q)\\ i>k+1; \sigma^{-1}(k)\leq p<\sigma^{-1}(k+1)\\ \{\sigma^{-1}(i),\sigma^{-1}(i+1)\}\subset\{p+1,\dots,p+q\}\end{smallmatrix}}\hskip-0.3cm \varepsilon(\sigma)(-1)^{\sum_{r<i}a_{\sigma^{-1}(r)}}\\
&\alpha_{\sigma^{-1}(1)}\underline\otimes\dots\underline\otimes [\alpha_{\sigma^{-1}(k)},\alpha_{\sigma^{-1}(k+1)}]\underline\otimes\dots\otimes \mu\big(\alpha_{\sigma^{-1}(i)},\alpha_{\sigma^{-1}(i+1)}\big)\underline\otimes\dots\underline\otimes\alpha_{\sigma^{-1}(p+q)}\\
&+\hskip-0.3cm\sum_{\begin{smallmatrix}\sigma\in Bat(p,q)\\ k ; \sigma^{-1}(k)\leq p<\sigma^{-1}(k+1)<\sigma^{-1}(k+2)\end{smallmatrix}}\hskip-0.3cm
\varepsilon(\sigma)(-1)^{\sum_{r<k+1}a_{\sigma^{-1}(r)}}\\
&\alpha_{\sigma^{-1}(1)}\underline\otimes\dots\underline\otimes [\alpha_{\sigma^{-1}(k)},\mu\big(\alpha_{\sigma^{-1}(k+1)},\alpha_{\sigma^{-1}(k+2)}\big)]
\underline\otimes\dots\underline\otimes\alpha_{\sigma^{-1}(p+q)}\\
&\hskip-0.5cm=(I'')+(II'')+(III'').
\end{aligned}
$$

On v\'erifie comme plus haut que $(I'')=(I_2)$, $(II'')=(II_2)$ et $(III'')=(III_2)+(IV_2)$, ce qui ach\`eve la preuve.
\hfill$\square$


\section{$G_\infty$ alg\`ebre}

\

\subsection{$L_\infty$ alg\`ebre associ\'ee \`a $\mathcal{H}$}

\

\

Par convention, dans toute la suite on notera le degr\'e  d'un \'el\'ement $\alpha_{[1,n]}=\alpha_1\underline{\otimes}\dots\underline{\otimes}\alpha_n\in\underline{\bigotimes }^n(\mathcal{G}[1])$ par $a_{[1,n]}=\sum_{i=1}^na_i$ et
on utilisera des lettres capitales pour les paquets, c'est \`a dire, pour un paquet $X=\alpha_1\underline{\otimes}\dots\underline{\otimes}\alpha_n\in\big(\underline{\bigotimes }^n(\mathcal{G}[1])\big)[1]\subset\mathcal{H}[1]$ le degr\'e sera not\'e $x=\sum_{i=1}^na_i-1=a_{[1,n]}-1$.\vskip0.2cm

Le degr\'e d'un \'el\'ement $X_1\pt\dots\pt X_n\in S^n\big(\mathcal{H}[1]\big) $ est alors $x_1+\dots+x_n$.
\vskip0.25cm

$(\mathcal{H}, \mu, [~,~])$ \'etant une alg\`ebre de Lie diff\'erentielle gradu\'ee. En suivant l'\'etude qu'on a fait dans la section $2$, on pourra construire une $L_\infty$ alg\`ebre associ\'ee \`a $\mathcal{H}$ not\'ee $\Big(S^+(\mathcal{H}[1]),\Delta, \ell+m\Big)$, avec $$(\ell+m)_2=\ell_2=[~,~]~~ \hbox{et} ~~(\ell+m)_1=m_1=\mu.$$

La comultiplication  $\Delta$ est d\'efinie sur $S^+(\mathcal{H}[1])$ par
$$
\Delta(X_1\pt\dots\pt X_n)=\sum_{\begin{smallmatrix}I\sqcup J=\{1,\dots n\}\\ \#I>0,\#J>0\end{smallmatrix}}\varepsilon\left(\begin{smallmatrix}x_{\{1,\dots,n\}}\\x_I,x_J\end{smallmatrix}\right)X_I\otimes X_J.
$$

Le crochet $[~,~]$ sur $\mathcal{H}$ \'etait antisy\'etrique de degr\'e $0$. Comme l'on veut une cod\'erivation de degr\'e $1$ pour $\Delta$, on pose $\ell_2(X,Y)=(-1)^{x}[X,Y]$. $\ell_2$ devient une application sym\'etrique sur $S^2(\mathcal{H}[1])$ de degr\'e $1$. Elle v\'erifie:
\begin{align*}&{\bf (i)} \
 \ell_2(X,Y)=(-1)^{xy}\ell_2(Y,X), \cr&{\bf (ii)} \ (-1)^{xz}\ell_2(\ell_2(X,Y),Z)+(-1)^{yx}\ell_2(\ell_2(Y,Z),X)+(-1)^{zy}\ell_2(\ell_2(Z,X),Y)=0,\cr&{\bf (iii)} \  m_1(\ell_2(X,Y))=-\ell_2(m_1(X),Y)
+(-1)^{1+x}\ell_2(X,m_1(Y)).\end{align*}

On prolonge le crochet $\ell_2$ \`a $S^+\big(\mathcal{H}[1]\big) $ de fa\c{c}on unique en une cod\'erivation de degr\'e $1$ de la cog\`ebre $(S^+(\mathcal{H}[1]),\Delta) $. Ce prolongement est donn\'e par:
$$
\ell(X_1\pt\dots\pt X_n)=\sum_{i<j}\varepsilon\left(\begin{smallmatrix}x_1\ \dots\ x_n\\ x_i x_j x_1\dots\hat{\imath}\hat{\jmath}\dots x_n\end{smallmatrix}\right)\ell_2(X_i,X_j)\pt X_1\pt\dots\hat\imath\dots \hat\jmath\dots\pt X_n.
$$
$\ell$ est de carr\'e nul $\ell\circ\ell=0$.

On prolonge, de m\^eme, $m_1$ \`a $S^n\big(\mathcal{H}[1]\big)$ en une d\'erivation $m$ de degr\'e $1$ par:
$$
m(X_1\pt\dots\pt X_n)=\sum_{j=1}^{n}(-1)^{\sum_{1\leq r<j}x_r}X_1\pt\dots\pt m_1(X_j)\pt\dots\pt X_n .
$$

$m$ est commutative et v\'erifie $m\circ m=0$.
\vskip0.3cm

De plus, gr\^ace \`a la propri\'et\'e (iii), on a $(\ell+m)\circ(\ell+m)=0$.

On peut voir alors $\Big(S^+\big(\mathcal{H}[1]\big),\Delta,\ell+m\Big) $ comme une $L_{\infty}$ alg\`ebre.

\


\subsection{Le cocrochet $\kappa$}

\

\

$(\mathcal{H},\delta , m_1)$ \'etant une $C_{\infty}$ alg\`ebre.
Le cocrochet $\delta$ sur $\underline{\bigotimes }^p(\mathcal{G}[1])$ devient un cocrochet $\kappa$ sur $\underline{\bigotimes }^p(\mathcal{G}[1])[1]$ d\'efini par:

Pour $X=\alpha_1\underline{\otimes}\dots\underline{\otimes}\alpha_p$,

 \begin{align*}\kappa(X)&=\sum_{j=1}^{p-1}
\Big((-1)^{a_{[1,j]}}\alpha_{[1,j]}\bigotimes \alpha_{[j+1,p]}-\varepsilon\left(\begin{smallmatrix}a_{[1,n]}\\ 
a_{[j+1,n]}~a_{[1,j]}\end{smallmatrix}\right)(-1)^{a_{[j+1,p]}}\alpha_{[j+1,p]}
\bigotimes\alpha_{\{1,\dots,j\}}\Big)\cr&=\sum_{j=1}^{p-1}(-1)^{u_j+1}\Big(
U_j\bigotimes V_j+
\varepsilon\left(\begin{smallmatrix}u_j~~v_j\\ v_j~ u_j\end{smallmatrix}\right)\alpha_{[j+1,p]}\bigotimes\alpha_{[1,j]}\Big).\end{align*}
o\`u $U_j=\alpha_{[1,j]}$, $V_j=\alpha_{[j+1,p]}$, $u_j=a_{[1,j]}-1$ et $v_j=a_{[j+1,p]-1}$.\vskip0.2cm

Autrement,
\begin{align*}\kappa(X)=\sum_{j=1}^{p-1}(-1)^{u_j+1}\Big(
U_j\bigotimes V_j+
\tau\big(U_j\bigotimes V_j\big)\Big).\end{align*}

Le cocrochet $\kappa$ sur $\underline{\bigotimes }^p(\mathcal{G}[1])[1]$ est alors cosym\'etrique ($\kappa=\tau\circ\kappa$) et de degr\'e $1$.

On prolonge $\kappa$ \`a $S^+\big(\mathcal{H}[1])$ par:
\begin{align*}
&\kappa(X_1\pt\dots\pt X_n)=\sum_{\begin{smallmatrix}1\leq s\leq n\\ I\cup J=\{1,\dots,n\}\setminus\{s\}\end{smallmatrix}} (-1)^{\sum_{i<s}x_i}\sum_{\begin{smallmatrix}U_s\otimes V_s=X_s\\ U_s,V_s\neq\emptyset\end{smallmatrix} } (-1)^{u_s+1}\times\cr
&\times\left(\varepsilon\left(\begin{smallmatrix} x_1\dots x_n\\ x_I~u_s~v_s~x_J\end{smallmatrix}\right)X_I\pt U_s\bigotimes V_s \pt X_J+\varepsilon\left(\begin{smallmatrix} x_1\dots x_n\\ x_I~v_s~u_s~x_J\end{smallmatrix}\right)X_I\pt V_s\bigotimes U_s \pt X_J\right),
\end{align*}
avec
\begin{align*}
\varepsilon\left(\begin{smallmatrix} x_1\dots x_n\\ x_I~u_s~v_s~x_J\end{smallmatrix}\right)=\varepsilon\left(\begin{smallmatrix} x_1\dots x_n\\ x_I~x_s~x_J\end{smallmatrix}\right)(-1)^{\sum_{\begin{smallmatrix}i<s\\ i\in J\end{smallmatrix}}x_i}(-1)^{\sum_{\begin{smallmatrix}i>s\\ i\in I\end{smallmatrix}}x_i}.
\end{align*}

En posant $\tau_{12}=\tau\otimes id$ et $\tau_{23}=id\otimes\tau$, $\kappa$ v\'erifie les identit\'es de coJacobi et de coLeibniz:

\begin{prop}

\

{\bf (i)} $\Big(id^{\otimes3}+\tau_{12}\circ\tau_{23}+\tau_{23}\circ\tau_{12}\Big)\circ(\kappa\otimes id)\circ\kappa=0$.
(identit\'e de coJacobi gradu\'ee)

 {\bf (ii)} $(id\otimes\Delta)\circ\kappa=(\kappa\otimes id)\circ\Delta+\tau_{12}\circ(id\otimes\kappa)\circ\Delta$. (identit\'e de coLeibniz gradu\'ee)

(Voir \cite{[BGHHW]})

\end{prop}
\vskip0.25cm

{\bf Preuve:}

\n (i) On calcule d'abord, $(\kappa\otimes id)\circ\kappa(X_1\pt\dots\pt X_n)$, on trouve pour $t\neq s$ des termes de la forme:

$$(1): \varepsilon_1. X_I\pt U_t\bigotimes V_t\pt X_J\pt U_s\bigotimes V_s\pt X_K$$
$$(2):\varepsilon_2.X_I\pt V_t\bigotimes U_t\pt X_J\pt U_s\bigotimes V_s\pt X_K$$
$$(3):\varepsilon_3 .X_I\pt U_t\bigotimes V_t\pt X_J\pt V_s\bigotimes U_s\pt X_K$$
$$(4):\varepsilon_4 .X_I\pt V_t\bigotimes U_t\pt X_J\pt V_s\bigotimes U_s\pt X_K$$
$$(5):\varepsilon_5 .X_I\pt U_s\pt U_t\bigotimes V_t\pt X_J\bigotimes V_s\pt X_K$$
$$(6):\varepsilon_6 .X_I\pt U_s\pt V_t\bigotimes U_t\pt X_J\bigotimes V_s\pt X_K$$
$$(7):\varepsilon_7 .X_I\pt V_s\pt U_t\bigotimes V_t\pt X_J\bigotimes U_s\pt X_K$$
$$(8):\varepsilon_8 .X_I\pt V_s\pt V_t\bigotimes U_t\pt X_J\bigotimes U_s\pt X_K$$
Et pour $t=s$, si $X_s=U_s\otimes V_s\otimes W_s$ on trouve des termes de la forme:
$$(9):\varepsilon_{9}.X_I\pt U_s\bigotimes V_s\pt X_J\bigotimes W_s\pt X_K$$
$$(10):\varepsilon_{10}.X_I\pt V_s\bigotimes U_s\pt X_J\bigotimes W_s\pt X_K$$
$$(11):\varepsilon_{11}.X_I\pt V_s\bigotimes W_s\pt X_J\bigotimes U_s\pt X_K$$
$$(12):\varepsilon_{12}.X_I\pt W_s\bigotimes V_s\pt X_J\bigotimes U_s\pt X_K$$
 On s'int\'eresse par exemple au terme $(1)=\varepsilon_1. X_I\pt U_t\bigotimes V_t\pt X_J\pt U_s\bigotimes V_s\pt X_K$ de $(\kappa\otimes id)\circ\kappa(X_1\pt\dots\pt X_n)$ et on cherche le terme correspondant dans $\Big(id^{\otimes3}+\tau_{12}\circ\tau_{23}+\tau_{23}\circ\tau_{12}\Big)\circ(\kappa\otimes id)\circ\kappa(X_1\pt\dots\pt X_n)$. Il apara\^it uniquement dans $\tau_{12}\circ\tau_{23}\circ(\kappa\otimes id)\circ\kappa(X_1\pt\dots\pt X_n)$ avec le signe $-\varepsilon_1$.

 En effet, le signe $\varepsilon_1$ dans $(1)$ est d\'etermin\'e par:

 \

- On part de $X_1,\dots,X_n$, on le ram\`ene en $X_I,X_t,X_J,X_s,X_K$ accompagn\'e du signe $\varepsilon\left(\begin{smallmatrix} x_1\dots x_n\\ x_I~x_t~x_J~x_s~x_K\end{smallmatrix}\right)$.

\

- On applique $\kappa$, le terme $X_I\pt X_t\pt X_J\pt U_s\bigotimes V_s\pt X_K$ appara\^it une seule fois avec le signe
$$(-1)^{x_I+x_t+x_J+u_s+1}\varepsilon\left(\begin{smallmatrix} x_1\dots x_n\\ x_I~x_t~x_J~x_s~x_K\end{smallmatrix}\right).$$

- Apr\`es en appliquant $(\kappa\otimes id)$, on obtient une seule fois le terme $X_I\pt U_t\bigotimes V_t\pt X_J\pt U_s\bigotimes V_s\pt X_K$ avec le signe
$$(-1)^{x_t+x_J+u_s+u_t}\varepsilon\left(\begin{smallmatrix} x_1\dots x_n\\ x_I~x_t~x_J~x_s~x_K\end{smallmatrix}\right)=(-1)^{v_t+u_s+x_J+1}\varepsilon\left(\begin{smallmatrix} x_1\dots x_n\\ x_I~x_t~x_J~x_s~x_K\end{smallmatrix}\right)=\varepsilon_1.$$

On cherche maintenant le signe du terme $(1)$ dans $\tau_{12}\circ\tau_{23}\circ(\kappa\otimes id)\circ\kappa(X_1\pt\dots\pt X_n)$.

\

- On part de $X_1,\dots,X_n$, on le ram\`ene en $X_J,X_s,X_K,X_t,X_I$ accompagn\'e du signe $\varepsilon\left(\begin{smallmatrix} x_1\dots x_n\\ x_J~x_s~x_K~x_t~x_I\end{smallmatrix}\right)$.

\

- On applique $\kappa$, le terme $X_J\pt X_s\pt X_K\pt V_t\bigotimes U_t\pt X_I$ appara\^it une seule fois avec le signe
$$(-1)^{x_J+x_s+x_K+u_t+1+u_tv_t}\varepsilon\left(\begin{smallmatrix} x_1\dots x_n\\ x_J~x_s~x_K~x_t~x_I\end{smallmatrix}\right)$$ qui s'\'ecrit encore $V_t\pt X_J\pt X_s\pt X_K\bigotimes X_I\pt U_t$ accompagn\'e du signe $$(-1)^{x_J+x_s+x_K+u_t+1+u_tv_t+u_tx_I+v_t(x_K+x_s+x_J)}\varepsilon\left(\begin{smallmatrix} x_1\dots x_n\\ x_J~x_s~x_K~x_t~x_I\end{smallmatrix}\right).$$

- Apr\`es en appliquant $(\kappa\otimes id)$, on obtient une seule fois le terme $V_t\pt U_s\bigotimes V_s\pt X_K\bigotimes X_I\pt U_t$ avec le signe
$$(-1)^{x_J+x_s+x_K+u_t+1+u_tv_t+u_tx_I+v_t(x_K+x_s+x_J)+x_J+v_t+u_s+1}\varepsilon\left(\begin{smallmatrix} x_1\dots x_n\\ x_J~x_s~x_K~x_t~x_I\end{smallmatrix}\right).$$

- On applique ensuite $\tau_{12}\circ\tau_{23}$, on obtient le terme $X_I\pt X_t\pt U_t\bigotimes V_t\pt X_J\pt U_s\bigotimes V_s\pt X_K$ avec le signe \begin{align*}&(-1)^{x_J+x_s+x_K+u_t+1+u_tv_t+u_tx_I+v_t(x_K+x_s+x_J)+x_J+v_t+u_s+1}\varepsilon\left(\begin{smallmatrix} x_1\dots x_n\\ x_J~x_s~x_K~x_t~x_I\end{smallmatrix}\right)\varepsilon\left(\begin{smallmatrix}v_t~x_J~u_s~v_s~x_K~x_I~u_t\\ x_I~u_t~v_t~x_J~u_s~v_s~x_K\end{smallmatrix}\right)\cr&
=(-1)^{x_J+x_s+x_K+u_t+1+u_tv_t+u_tx_I+v_t(x_K+x_s+x_J)+x_J+v_t+u_s+1}\times
\cr&\times\varepsilon\left(\begin{smallmatrix} x_1\dots x_n\\ x_J~x_s~x_K~x_t~x_I\end{smallmatrix}\right)\varepsilon\left(\begin{smallmatrix}x_J~x_s~x_K~x_t~x_I\\ x_I~x_t~x_J~x_s~x_K\end{smallmatrix}\right)(-1)^{x_I+(v_t+1)
(x_J+x_K+x_s)+u_t+v_t(x_I+u_t)+x_tx_I}\cr&=-\varepsilon_1.\end{align*}

Les autres termes se simplifient de fa\c{c}on pareille.

\

(ii) D'une part, on a

\begin{align*}&(id\otimes\Delta)\circ\kappa(X_1\pt\dots\pt X_n)=\sum_{\begin{smallmatrix}1\leq s\leq n\\ I\cup J\cup K=\{1,\dots,n\}\setminus\{s\}\end{smallmatrix}}(-1)^{\sum_{i<s}x_i}\sum_{\begin{smallmatrix}U_s\otimes V_s=X_s\\ U_s,V_s\neq\emptyset\end{smallmatrix} } (-1)^{u_s+1}\times\cr&\times\Big(\varepsilon\left(\begin{smallmatrix} x_1\dots x_n\\ x_J~u_s~x_I~v_s~x_K\end{smallmatrix}\right)X_J\pt U_s\bigotimes X_I\bigotimes V_s \pt X_K+\varepsilon\left(\begin{smallmatrix} x_1\dots x_n\\ x_J~u_s~v_s~x_K~x_I\end{smallmatrix}\right)X_J\pt U_s\bigotimes V_s\pt X_K\bigotimes X_I\cr&+\varepsilon\left(\begin{smallmatrix} x_1\dots x_n\\ x_J~v_s~x_I~u_s~x_K\end{smallmatrix}\right)X_J\pt V_s\bigotimes X_I\bigotimes U_s \pt X_K+\varepsilon\left(\begin{smallmatrix} x_1\dots x_n\\ x_J~v_s~u_s~x_K~x_I\end{smallmatrix}\right) X_J\pt V_s\bigotimes U_s\pt X_K\bigotimes X_I\Big)\cr&=(1)+(2)+(3)+(4).\end{align*}

D'autre part, on sait que $\Delta(X_1\pt\dots\pt X_n)=\displaystyle\sum_{\begin{smallmatrix}I\cup J=\{1,\dots,n\}\cr \#I,\# J>0\end{smallmatrix}}\varepsilon\left(\begin{smallmatrix} x_1\dots x_n\\ x_I~x_J\end{smallmatrix}\right)X_I\bigotimes X_J.$

Donc,
\begin{align*}&(id\otimes\kappa)\circ\Delta(X_1\pt\dots\pt X_n)=\sum_{\begin{smallmatrix}1\leq s\leq n\\ I\cup J\cup K=\{1,\dots,n\}\setminus\{s\}\end{smallmatrix}}\varepsilon\left(\begin{smallmatrix} x_1\dots x_n\\ x_I~x_J~x_s~x_K\end{smallmatrix}\right)(-1)^{x_I}(-1)^{\sum_{ i\in J}x_i}\sum_{\begin{smallmatrix}U_s\otimes V_s=X_s\\ U_s,V_s\neq\emptyset\end{smallmatrix} } (-1)^{u_s+1}\cr&\times\Big(\varepsilon\left(\begin{smallmatrix} x_I~x_J~ x_s~x_K\\ x_I~x_J~u_s~v_s~x_K\end{smallmatrix}\right)X_I\bigotimes X_J\pt U_s \bigotimes V_s \pt X_K +\varepsilon\left(\begin{smallmatrix} x_I~x_J~ x_s~x_K\\ x_I~x_J~v_s~u_s~x_K\end{smallmatrix}\right)X_I\bigotimes X_J\pt V_s\bigotimes U_s\pt X_K\Big).\end{align*}
Alors,
\begin{align*}&\tau_{12}\circ(id\otimes\kappa)\circ\Delta(X_1\pt\dots\pt X_n)=\sum_{\begin{smallmatrix}1\leq s\leq n\\ I\cup J\cup K=\{1,\dots,n\}\setminus\{s\}\end{smallmatrix}}\hskip-0.5cm\varepsilon\left(\begin{smallmatrix} x_1\dots x_n\\ x_I~x_J~x_s~x_K\end{smallmatrix}\right)(-1)^{x_I}(-1)^{\sum_{ i\in J}x_i}\sum_{\begin{smallmatrix}U_s\otimes V_s=X_s\\ U_s,V_s\neq\emptyset\end{smallmatrix} } (-1)^{u_s+1}\cr&\times\Big(\varepsilon\left(\begin{smallmatrix} x_I~x_J~ u_s~v_s~x_K\\ x_J~u_s~x_I~v_s~x_K\end{smallmatrix}\right)X_J\pt U_s\bigotimes X_I\bigotimes V_s\pt X_K +\varepsilon\left(\begin{smallmatrix} x_I~x_J~ u_s~v_s~x_K\\ x_J~v_s~x_I~u_s~x_K\end{smallmatrix}\right)X_J\pt V_s\bigotimes X_I\bigotimes U_s\pt X_K
\Big)\cr&=(5)+(6).\end{align*}

De plus, en \'ecrivant que $\Delta(X_1\pt\dots\pt X_n)=\displaystyle\sum_{\begin{smallmatrix}I\cup J=\{1,\dots,n\}\cr \#I,\# J>0\end{smallmatrix}}\varepsilon\left(\begin{smallmatrix} x_1\dots x_n\\ x_J~x_I\end{smallmatrix}\right)X_J\bigotimes X_I$, on a
\begin{align*}&(\kappa\otimes id)\circ\Delta(X_1\pt\dots\pt X_n)=\sum_{\begin{smallmatrix}1\leq s\leq n\\ I\cup J\cup K=\{1,\dots,n\}\setminus\{s\}\end{smallmatrix} }\varepsilon\left(\begin{smallmatrix} x_1\dots x_n\\ x_J~x_s~x_K~x_I\end{smallmatrix}\right)(-1)^{\sum_{ i\in J}x_i}\sum_{\begin{smallmatrix}U_s\otimes V_s=X_s\\ U_s,V_s\neq\emptyset\end{smallmatrix}} (-1)^{u_s+1}\cr&\times\Big(\varepsilon\left(\begin{smallmatrix}x_J~x_s~x_K~x_I\\ x_J~u_s~v_s~x_K~x_I\end{smallmatrix}\right)X_J\pt U_s\bigotimes V_s\pt X_K\bigotimes X_I +\varepsilon\left(\begin{smallmatrix}x_J~x_s~x_K~x_I\\ x_J~v_s~u_s~x_K~x_I\end{smallmatrix}\right)X_J\pt V_s\bigotimes U_s\pt X_K\bigotimes X_I\Big)\cr&=(7)+(8).\end{align*}

Un calcul nous montre que $(1)=(5)$, $(2)=(7)$, $(3)=(6)$ et $(4)=(8)$. Montrons par exemple que $(1)=(5)$.

En effet, dans $(1)$, le terme $X_J\pt U_s\bigotimes X_I\bigotimes V_s \pt X_K$ appara\^it avec le signe $$(-1)^{\sum_{i<s}x_i+u_s+1}\varepsilon\left(\begin{smallmatrix} x_1\dots x_n\\ x_J~u_s~x_I~v_s~x_K\end{smallmatrix}\right).$$
Dans $(5)$, ce terme appara\^it avec le signe
\begin{align*}
&\varepsilon\left(\begin{smallmatrix} x_1\dots x_n\\ x_I~x_J~x_s~x_K\end{smallmatrix}\right)(-1)^{\sum_{i\in J}x_i+x_I+u_s+1}\varepsilon\left(\begin{smallmatrix} x_I~x_J~ u_s~v_s~x_K\\ x_J~u_s~x_I~v_s~x_K\end{smallmatrix}\right)\cr
&=(-1)^{\sum_{\begin{smallmatrix} i>s\\ i\in I\cup J\end{smallmatrix}}x_i}(-1)^{\sum_{ \begin{smallmatrix}i<s\\ i\in K\end{smallmatrix}}x_i}(-1)^{\sum_{i\in J}x_i+x_I+u_s+1}\varepsilon\left(\begin{smallmatrix} x_1\dots x_n\\ x_I~x_J~u_s~v_s~x_K\end{smallmatrix}\right)\varepsilon\left(\begin{smallmatrix} x_I~x_J~ u_s~v_s~x_K\\ x_J~u_s~x_I~v_s~x_K\end{smallmatrix}\right)\cr
&=(-1)^{\sum_{i<s}x_i+u_s+1}\varepsilon\left(\begin{smallmatrix} x_1\dots x_n\\ x_J~u_s~x_I~v_s~x_K\end{smallmatrix}\right).\end{align*}
Ce qui montre que $(1)=(5)$.
\vskip0.3cm

Donc, on trouve que $(id\otimes\Delta)\circ\kappa=(\kappa\otimes id)\circ\Delta+\tau_{12}\circ(id\otimes\kappa)\circ\Delta$.\hfill$\square$

\

Finalement, l'espace $(S^+(\mathcal{H}[1]),\kappa)$ est une cog\`ebre de Lie. On v\'erifie que $m$ est une cod\'erivation de degr\'e $1$ de $S^+(\mathcal{H}[1])$ pour le cocrochet $\kappa$.
\begin{prop}
 $$\Big(id\otimes m+m\otimes id\Big)\circ\kappa=-\kappa\circ m.$$

\end{prop}

{\bf Preuve:}

On a \begin{align*}&\kappa\circ m(X_1\pt \dots\pt X_n)=\sum_{s=1}^n(-1)^{\sum_{i<s}x_i}\kappa(X_1\pt \dots\pt m(X_s)\pt\dots\pt X_n)\cr&=\sum_{t<s}(-1)^{\sum_{i<s}x_i+u_t+1+\sum_{i<t}x_i}\Big(\cr&\varepsilon\left(\begin{smallmatrix} x_1\dots u_tv_t\dots (x_s+1)\dots x_n\\ x_I~u_t~v_t~(x_s+1)~x_J\end{smallmatrix}\right)X_I\pt U_t\bigotimes V_t\pt m(X_s)\pt X_J+\varepsilon\left(\begin{smallmatrix} x_1\dots u_tv_t\dots (x_s+1)\dots x_n\\ x_I~(x_s+1)~u_t~v_t~x_J\end{smallmatrix}\right) X_I\pt m(X_s)\pt U_t\bigotimes V_t\pt  X_J\cr&+\varepsilon\left(\begin{smallmatrix} x_1\dots u_tv_t\dots (x_s+1)\dots x_n\\ x_I~v_t~u_t~(x_s+1)~x_J\end{smallmatrix}\right)X_I\pt V_t\bigotimes U_t\pt m(X_s)\pt X_J+\varepsilon\left(\begin{smallmatrix} x_1\dots u_tv_t\dots (x_s+1)\dots x_n\\ x_I~(x_s+1)~v_t~u_t~x_J\end{smallmatrix}\right)X_I\pt m(X_s)\pt V_t\bigotimes U_t\pt  X_J\Big)\cr&
+\sum_{t>s}(-1)^{\sum_{i<s}x_i+u_t+1+\sum_{i<t}x_i+1}\Big(\cr&\varepsilon\left(\begin{smallmatrix} x_1\dots u_tv_t\dots (x_s+1)\dots x_n\\ x_I~u_t~v_t~(x_s+1)~x_J\end{smallmatrix}\right)X_I\pt U_t\bigotimes V_t\pt m(X_s)\pt X_J+\varepsilon\left(\begin{smallmatrix} x_1\dots u_tv_t\dots (x_s+1)\dots x_n\\ x_I~(x_s+1)~u_t~v_t~x_J\end{smallmatrix}\right)X_I\pt m(X_s)\pt U_t\bigotimes V_t\pt  X_J\cr&+\varepsilon\left(\begin{smallmatrix} x_1\dots u_tv_t\dots (x_s+1)\dots x_n\\ x_I~v_t~u_t~(x_s+1)~x_J\end{smallmatrix}\right)X_I\pt V_t\bigotimes U_t\pt m(X_s)\pt X_J+\varepsilon\left(\begin{smallmatrix} x_1\dots u_tv_t\dots (x_s+1)\dots x_n\\ x_I~(x_s+1)~v_t~u_t~x_J\end{smallmatrix}\right)X_I\pt m(X_s)\pt V_t\bigotimes U_t\pt  X_J\Big)\cr&+\sum_{s=1}^n(-1)^{\sum_{i<s}x_i+(u_s+1)+1}\Big(\cr&\varepsilon\left(\begin{smallmatrix} x_1\dots (u_s+1)~v_s\dots x_n\\ x_I~(u_s+1)~v_s~x_J\end{smallmatrix}\right)X_I\pt m(U_s)\bigotimes V_s\pt X_J+ \varepsilon\left(\begin{smallmatrix} x_1\dots (u_s+1)~v_s\dots x_n\\ x_I~v_s~(u_s+1)~x_J\end{smallmatrix}\right)X_I\pt V_s\bigotimes m(U_s) \pt X_J\Big)\cr&+(-1)^{\sum_{i<s}x_i+
 u_s+1}\Big(\cr&\varepsilon\left(\begin{smallmatrix} x_1\dots u_s~ (v_s+1)\dots x_n\\ x_I~u_s ~(v_s+1)~x_J\end{smallmatrix}\right)X_I\pt U_s\bigotimes m(V_s)\pt X_J+ \varepsilon\left(\begin{smallmatrix} x_1\dots u_s~ (v_s+1)\dots x_n\\ x_I~(v_s+1)~u_s ~x_J\end{smallmatrix}\right)X_I\pt m(V_s)\bigotimes U_s \pt X_J\Big)\cr&=(1)+(2)+(3)+(4)+(5)+(6)+(7)+(8)+(9)+(10)+(11)+(12).\end{align*}

\

Apr\`es, en calculant $(id\otimes m)\circ\kappa(X_1\pt \dots\pt X_n)$, on obtient $-(1)-(3)-(5)-(7)-(10)-(11)$ et en calculant $(m\otimes id)\circ\kappa(X_1\pt \dots\pt X_n)$, on obtient $-(2)-(4)-(6)-(8)-(9)-(12)$. Par exemple, le terme $X_I\pt U_t\bigotimes V_t\pt m(X_s)\pt X_J$ de $(1)$ appara\^it dans $\kappa\circ m(X_1\pt \dots\pt X_n)$ avec le signe $$(-1)^{\sum_{i<s}x_i+u_t+1+\sum_{i<t}x_i}\varepsilon\left(\begin{smallmatrix} x_1\dots u_tv_t\dots (x_s+1)\dots x_n\\ x_I~u_t~v_t~(x_s+1)~x_J\end{smallmatrix}\right).$$
Ce m\^eme terme appa\^it dans $(id\otimes m)\circ\kappa(X_1\pt \dots\pt X_n)$ avec le signe \begin{align*}
&(-1)^{\sum_{i<t}x_i+u_t+1}(-1)^{x_I+u_t+v_t}\varepsilon\left(\begin{smallmatrix} x_1\dots u_tv_t\dots x_s\dots x_n\\ x_I~u_t~v_t~x_s~x_J\end{smallmatrix}\right)\cr
&=\varepsilon\left(\begin{smallmatrix} x_1\dots u_tv_t\dots (x_s+1)\dots x_n\\ x_I~u_t~v_t~(x_s+1)~x_J\end{smallmatrix}\right)(-1)^{\sum_{i<t}x_i+u_t+1}(-1)^{x_I+u_t+v_t}(-1)^{\sum_{\begin{smallmatrix}i<s\\ i\in J\end{smallmatrix}}x_i+\sum_{\begin{smallmatrix}i>s\\ i\in I\end{smallmatrix}}x_i}\cr&=(-1)^{v_t+\sum_{t< i<s}x_i+1}\varepsilon\left(\begin{smallmatrix} x_1\dots u_tv_t\dots (x_s+1)\dots x_n\\ x_I~u_t~v_t~(x_s+1)~x_J\end{smallmatrix}\right).\end{align*}
D'o\`u, le r\'esultat.\hfill$\square$
\vskip0.3cm

Aussi, on v\'erifie aussi que le crochet $\ell$ est une cod\'erivation de degr\'e $1$ de $S^+(\mathcal{H}[1])$ pour $\kappa$:
 \begin{prop}
 $$\Big(id\otimes \ell+\ell\otimes id\Big)\circ\kappa=-\kappa\circ \ell.$$\
\end{prop}

{\bf Preuve:}

On a d'une part,
\begin{align*}\kappa\circ\ell(X_1\pt\dots\pt X_n)&=\kappa\big(\sum_{i<j}\varepsilon\left(\begin{smallmatrix} x_1\dots x_n\\ x_i~x_j~x_1\dots\hat{i}\hat{j}\dots x_n\end{smallmatrix}\right)\ell_2(X_i,X_j)X_1\pt\dots\check{i}\dots\check{j}\dots\pt X_n\big)\cr&=\kappa\big(\sum_{\begin{smallmatrix}i<j\\ J=\{1,\dots,n\}\setminus\{i,j\}\end{smallmatrix}}\varepsilon\left(\begin{smallmatrix} x_1\dots x_n\\ x_i~x_jx_J\end{smallmatrix}\right)\ell_2(X_i,X_j)\pt X_J\big)
\end{align*}
Dans $\kappa\circ\ell(X_1\pt\dots\pt X_n)$, il appara\^it un terme $(I)$ de la forme $$\ell_2(X_i,X_j)\pt X_{J_1}\pt U_s\bigotimes V_s\pt X_{J_2}$$ avec le signe $$\varepsilon_1=(-1)^{x_i+x_j+x_{J_1}+u_s}\varepsilon\left(\begin{smallmatrix} x_1\dots x_n\\ x_i~x_jx_J\end{smallmatrix}\right)$$ o\`u on a pos\'e $J=J_1\cup\{s\}\cup J_2$.

D'autre part, cherchons le terme correspondant dans $(\ell\otimes id)\circ\kappa(X_1\pt \dots\pt X_n)$.

On a \begin{align*}&(\ell\otimes id)\circ\kappa(X_1\pt \dots\pt X_n)=(\ell\otimes id)\big(\sum_{i<j}\sum_{\begin{smallmatrix}s\\ J_1\cup J_2=\{1,\dots,n\}\setminus\{i,j,s\}\end{smallmatrix}}\sum_{\begin{smallmatrix}U_s\otimes V_s=X_s\\ U_s,V_s\neq\emptyset\end{smallmatrix}}(-1)^{\sum_{r<s}x_r+u_s+1}\times\cr&\times\varepsilon\left(\begin{smallmatrix} x_1\dots u_s~v_s\dots x_n\\ x_i~x_j~x_{J_1}~u_s~v_s~x_{J_2}\end{smallmatrix}\right) X_i\pt X_j\pt X_{J_1}\pt U_s\bigotimes V_s\pt X_{J_2}+\cdots\big)
\cr&=(-1)^{\sum_{r<s}x_r+u_s+1}\varepsilon\left(\begin{smallmatrix} x_1\dots u_s~v_s\dots x_n\\ x_i~x_j~x_{J_1}~u_s~v_s~x_{J_2}\end{smallmatrix}\right)\ell_2(X_i,X_j)\pt X_{J_1}\pt U_s\bigotimes V_s\pt X_{J_2}+\cdots.\end{align*}
Le premier terme $(1)=\ell_2(X_i,X_j)\pt X_{J_1}\pt U_s\bigotimes V_s\pt X_{J_2}$ appara\^it donc avec le signe
\begin{align*}&(-1)^{\sum_{r<s}x_r+u_s+1}(-1)^{\sum_{\begin{smallmatrix}r<s\\ r\in J_2\end{smallmatrix}}x_r}(-1)^{\sum_{\begin{smallmatrix}r>s\\ r\in J_1\cup\{i,j\}\end{smallmatrix}}x_r}\varepsilon\left(\begin{smallmatrix} x_1\dots x_s\dots x_n\\ x_i~x_j~x_{J_1}~x_s~x_{J_2}\end{smallmatrix}\right)
\cr&=(-1)^{x_i+x_j+x_{J_1}+u_s+1}\varepsilon\left(\begin{smallmatrix} x_1\dots x_n\\ x_i~x_j~x_{J}\end{smallmatrix}\right)=-\varepsilon_1.\end{align*}
Alors, $(I)=-(1)$. De m\^eme, les autres termes apparaissent dans le premier membre et le second membre \`a un signe $(-1)$ pr\`es.

\

Juste le cas o\`u on coupe $\ell_2(X_i,X_j)$ par $\kappa$, on va l'\'etudier comme le cas o\`u on a deux paquets.

En effet,
si $X=\alpha_1\underline{\otimes}\dots\underline{\otimes}\alpha_p$ et $Y=\alpha_{p+1}\underline{\otimes}\dots\underline{\otimes}\alpha_{p+q}$, on sait que
\begin{align*}&\ell_2(X,Y)=\cr&(-1)^x\sum_{\begin{smallmatrix}\sigma\in Bat(p,q)\\ \sigma^{-1}(k)\leq p<\sigma^{-1}(k+1)\end{smallmatrix}}\hskip-0.5cm
\varepsilon\left(\begin{smallmatrix}a_1~~\dots~~a_{p+q}\\ a_{\sigma^{-1}(1)}\dots a_{\sigma^{-1}(p+q)}\end{smallmatrix}\right) \alpha_{\sigma^{-1}(1)}\underline{\otimes}\dots\underline{\otimes}[\alpha_{\sigma^{-1}(k)},\alpha_{\sigma^{-1}(k+1)}]\underline{\otimes}
\dots\underline{\otimes}\alpha_{\sigma^{-1}(p+q)}.\end{align*}
Alors,
\begin{align*}&\kappa\circ\ell_2(X,Y)=(-1)^x\sum_{\begin{smallmatrix}\sigma\in Bat(p,q)\\ k ; \ \sigma^{-1}(k)\leq p<\sigma^{-1}(k+1)\end{smallmatrix}}
\varepsilon\left(\begin{smallmatrix}a_1~~\dots~~a_{p+q}\\ a_{\sigma^{-1}(1)}\dots a_{\sigma^{-1}(p+q)}\end{smallmatrix}\right) \sum_{t<k}\Big(\cr&(-1)^{\sum_{i\leq t}a_{\sigma^{-1}(i)}}\alpha_{\sigma^{-1}(1)}\underline{\otimes}\dots\underline{\otimes}\alpha_{\sigma^{-1}(t)} \bigotimes\alpha_{\sigma^{-1}(t+1)}\underline{\otimes}\dots\underline{\otimes} [\alpha_{\sigma^{-1}(k)},\alpha_{\sigma^{-1}(k+1)}]
\underline{\otimes}\dots\underline{\otimes}\alpha_{\sigma^{-1}(p+q)}\cr&-(-1)^{\sum_{i> t}a_{\sigma^{-1}(i)}}\varepsilon\left(\begin{smallmatrix}a_1~~\dots~~a_{p+q}\\ a_{[\sigma^{-1}(t+1), \sigma^{-1}(p+q)]},a_{[\sigma^{-1}(1), \sigma^{-1}(t)]}\end{smallmatrix}\right)\cr&\alpha_{\sigma^{-1}(t+1)}\underline{\otimes}\dots\underline{\otimes} [\alpha_{\sigma^{-1}(k)},\alpha_{\sigma^{-1}(k+1)}]
\underline{\otimes}\dots\underline{\otimes}\alpha_{\sigma^{-1}(p+q)}\bigotimes\alpha_{\sigma^{-1}(1)}\underline{\otimes}\dots\underline{\otimes}\alpha_{\sigma^{-1}(t)}\Big)\cr&+\sum_{t>k}\Big
(\cr&(-1)^{\sum_{i\leq t}a_{\sigma^{-1}(i)}}\alpha_{\sigma^{-1}(1)}\underline{\otimes}\dots\underline{\otimes} [\alpha_{\sigma^{-1}(k)},\alpha_{\sigma^{-1}(k+1)}]\underline{\otimes}\dots\underline{\otimes}\alpha_{\sigma^{-1}(t)} \bigotimes\alpha_{\sigma^{-1}(t+1)}
\underline{\otimes}\dots\underline{\otimes}\alpha_{\sigma^{-1}(p+q)}\cr&-(-1)^{\sum_{i> t}a_{\sigma^{-1}(i)}}\varepsilon\left(\begin{smallmatrix}a_1~~\dots~~a_{p+q}\\ a_{[\sigma^{-1}(t+1), \sigma^{-1}(p+q)]},a_{[\sigma^{-1}(1), \sigma^{-1}(t)]}\end{smallmatrix}\right)
\cr&\alpha_{\sigma^{-1}(t+1)}\underline{\otimes}\dots\underline{\otimes}\alpha_{\sigma^{-1}(p+q)}\bigotimes\alpha_{\sigma^{-1}(1)}\underline{\otimes}\dots\underline{\otimes} [\alpha_{\sigma^{-1}(k)},\alpha_{\sigma^{-1}(k+1)}]
\underline{\otimes}\dots\underline{\otimes}\alpha_{\sigma^{-1}(t)}\Big)\cr&=(1)+(2)+(3)+(4).\end{align*}
Sachant que $\kappa$ passe au quotient modulo les battements, alors, par exemple pour le terme $(1)=\alpha_{\sigma^{-1}(1)}\underline{\otimes}\dots\underline{\otimes}\alpha_{\sigma^{-1}(t)} \bigotimes\alpha_{\sigma^{-1}(t+1)}\underline{\otimes}\dots\underline{\otimes} [\alpha_{\sigma^{-1}(k)},\alpha_{\sigma^{-1}(k+1)}]
\underline{\otimes}\dots\underline{\otimes}\alpha_{\sigma^{-1}(p+q)}$ qui appara\^it dans $\kappa\circ\ell_2(X,Y)$, on a n\'ecessairement $\alpha_{\sigma^{-1}(1)},\dots,\alpha_{\sigma^{-1}(t)}$ appartiennent tous \`a $X$ ou \`a $Y$. Alors, on a  n\'ecessairement $\alpha_{\sigma^{-1}(1)}\underline{\otimes}\dots\underline{\otimes}\alpha_{\sigma^{-1}(t)}=\alpha_1\underline{\otimes}\dots\underline{\otimes}\alpha_t$ ou $\alpha_{\sigma^{-1}(1)}\underline{\otimes}\dots\underline{\otimes}\alpha_{\sigma^{-1}(t)}=\alpha_{p+1}\underline{\otimes}\dots\underline{\otimes}\alpha_{p+t}.$

Supposons, par exemple, que le terme $(1)$ s'\'ecrit: $$\alpha_{p+1}\underline{\otimes}\dots\underline{\otimes}\alpha_{p+t}\bigotimes\alpha_{\sigma^{-1}(t+1)}\underline{\otimes}\dots\underline{\otimes} [\alpha_{\sigma^{-1}(k)},\alpha_{\sigma^{-1}(k+1)}]
\underline{\otimes}\dots\underline{\otimes}\alpha_{\sigma^{-1}(p+q)}.$$
Il appara\^it dans $\kappa\circ\ell_2(X,Y)$ avec le signe $(-1)^x\varepsilon\left(\begin{smallmatrix}a_1~~\dots~~a_{p+q}\\ a_{\sigma^{-1}(1)}\dots a_{\sigma^{-1}(p+q)}\end{smallmatrix}\right)(-1)^{\sum_{p+1\leq i\leq p+t}a_{i}}$.

\

Cherchons le terme correspondant dans $(id\otimes\ell_2)\circ\kappa(X,Y)$.

\

Posons $U=\alpha_{p+1}\underline{\otimes}\dots\underline{\otimes}\alpha_{p+t}$ et $V=\alpha_{p+t+1}\underline{\otimes}\dots\underline{\otimes}\alpha_{p+q}$. On sait que \begin{align*}&\kappa(X,Y)=\sum_{t=1}^q(-1)^{x+u+1}\varepsilon\left(\begin{smallmatrix}x~u~v\\ u~v~x\end{smallmatrix}\right) U\bigotimes(\alpha_{p+t+1}
\underline{\otimes}\dots\underline{\otimes}\alpha_{p+q})\pt(\alpha_1\underline{\otimes}\dots\underline{\otimes}\alpha_{p})+\hbox{d'autres \ termes...}\end{align*}
En appliquant $(id\otimes \ell_2)$, on obtient:
 \begin{align*}&(id\otimes \ell_2)\circ\kappa(X,Y)=\sum_{t=1}^q(-1)^{u}(-1)^{x+u+1}\varepsilon\left(\begin{smallmatrix}x~u~v\\ u~v~x\end{smallmatrix}\right)U\bigotimes\ell_2(\alpha_{p+t+1}
\underline{\otimes}\dots\underline{\otimes}\alpha_{p+q}\pt\alpha_1\underline{\otimes}\dots\underline{\otimes}\alpha_{p})\cr&+\hbox{d'autres \ termes...}\end{align*}
On consid\`ere une permutation $\lambda $ de  $S_{p+q}$ d\'efinie par:
$$\lambda(i)=\left\{
\begin{array}{ll}
p+i, & \hbox{si $1\leq i\leq q$;} \\
i-q, & \hbox{si $q+1\leq i\leq p+q$.}
                     \end{array}
                          \right.$$
Posons $\beta_i=\alpha_{\lambda(i)}$. Alors,

\begin{align*}&(id\otimes \ell_2)\circ\kappa(X,Y)
=\sum_{\begin{smallmatrix} t<k\cr\eta\in Bat(q-t,p); \ \eta^{-1}(k)\leq q<\eta^{-1}(k+1)\end{smallmatrix}}(-1)^{x+v}\varepsilon\left(\begin{smallmatrix}x~u~v\\ u~v~x\end{smallmatrix}\right)\varepsilon\left(\begin{smallmatrix}b_{[t+1,p+q]}\\ b_{[\eta^{-1}(t+1),\eta^{-1}(p+q)]}\end{smallmatrix}\right)\cr&U\bigotimes\beta_{\eta^{-1}(t+1)}
\underline{\otimes}\dots\underline{\otimes}[\beta_{\eta^{-1}(k)},\beta_{\eta^{-1}(k+1)}]\underline{\otimes}\dots\underline{\otimes}\beta_{\eta^{-1}(p+q)}+\hbox{d'autres \ termes...}\end{align*}
o\`u $\eta$ est un battement de $Bat(q-t,p)$ d\'efinie sur $\{t+1,\dots,p+q\}$.

\

Dans la somme pr\'ec\'edente, fixant $\eta\in Bat(q-t,p)$ telle que $\eta(i)=\sigma\circ\lambda(i), \forall i\in\{t+1,\dots,p+q\}$. On construit, apr\`es, une permutation $\nu$ de $S_{p+q}$ d\'efinie par:
$$\nu^{-1}(i)=\left\{
\begin{array}{ll}
p+i, & \hbox{si $1\leq i\leq t$} \\
\eta^{-1}(i), & \hbox{si $t+1\leq i\leq p+q$.}
                     \end{array}
                          \right.
$$
On v\'erifie que $\nu$ appartient \`a $Bat(q,p)$, $\beta_{\nu^{-1}(i)}=\alpha_{\sigma^{-1}(i)}, \forall i\in\{1,\dots,p+q\}$ et que $$\varepsilon\left(\begin{smallmatrix}b_{[1,p+q]}\\ b_{[\nu^{-1}(1),\nu^{-1}(p+q)]}\end{smallmatrix}\right)=(-1)^{a_{[1,p]}a_{[p+1,p+q]}}\varepsilon\left(\begin{smallmatrix}a_{[1,p+q]}\\ a_{[\sigma^{-1}(1),\sigma^{-1}(p+q)]}\end{smallmatrix}\right).$$

Mais, par construction $\nu^{-1}(k)\leq q<\nu^{-1}(k+1)$. On construit, alors, une nouvelle permutation $\rho$ de $S_{p+q}$ d\'efinie par:

$\rho^{-1}(i)=\nu^{-1}(i), \ \forall i\notin\{k, k+1\}$, $\rho^{-1}(k)=\nu^{-1}(k+1)$ et $\rho^{-1}(k+1)=\nu^{-1}(k)$.

\vskip0.25cm

On v\'erifie que $\rho$ appartient \`a $Bat(q,p)$, $\beta_{\rho^{-1}(i)}=\alpha_{\sigma^{-1}(i)}, \forall i\notin\{k,k+1\}$, $\beta_{\rho^{-1}(k)}=\alpha_{\sigma^{-1}(k+1)}$, $\beta_{\rho^{-1}(k+1)}=\alpha_{\sigma^{-1}(k)}$, $\rho^{-1}(k)\leq p<\rho^{-1}(k+1)$ et que $$\varepsilon\left(\begin{smallmatrix}b_{[1,p+q]}\\ b_{[\rho^{-1}(1),\rho^{-1}(p+q)]}\end{smallmatrix}\right)=(-1)^{a_{[1,p]}a_{[p+1,p+q]}}\varepsilon\left(\begin{smallmatrix}a_{[1,p+q]}\\ a_{[\sigma^{-1}(1),\sigma^{-1}(p+q)]}\end{smallmatrix}\right)(-1)^{a_{\sigma^{-1}(k)}a_{\sigma^{-1}(k+1)}}.$$

Enfin, le terme \begin{align*}&U\bigotimes\beta_{\rho^{-1}(t+1)}\underline{\otimes}\dots\underline{\otimes} [\beta_{\rho^{-1}(k)},\beta_{\rho^{-1}(k+1)}]
\underline{\otimes}\dots\underline{\otimes}\beta_{\rho^{-1}(p+q)}\cr&=\alpha_{p+1}\underline{\otimes}\dots\underline{\otimes}\alpha_{p+t}\bigotimes\alpha_{\sigma^{-1}(t+1)}\underline{\otimes}\dots\underline{\otimes} [\alpha_{\sigma^{-1}(k+1)},\alpha_{\sigma^{-1}(k)}]
\underline{\otimes}\dots\underline{\otimes}\alpha_{\sigma^{-1}(p+q)}\end{align*}  appara\^it dans $(id\otimes \ell_2)\circ\kappa(X,Y)$ avec le signe
$$(-1)^{a_{[1,p]}+v+1}\varepsilon\left(\begin{smallmatrix}x u v\\ u v x \end{smallmatrix}\right)(-1)^{a_{[1,p]}a_{[p+1,p+q]}}\varepsilon\left(\begin{smallmatrix}a_{[1,p+q]}\\ a_{[\sigma^{-1}(1),\sigma^{-1}(p+q)]}\end{smallmatrix}\right)(-1)^{a_{\sigma^{-1}(k)}a_{\sigma^{-1}(k+1)}}.
$$
Donc, le terme $\alpha_{p+1}\underline{\otimes}\dots\underline{\otimes}\alpha_{p+t}\bigotimes\alpha_{\sigma^{-1}(t+1)}\underline{\otimes}\dots\underline{\otimes} [\alpha_{\sigma^{-1}(k)},\alpha_{\sigma^{-1}(k+1)}]
\underline{\otimes}\dots\underline{\otimes}\alpha_{\sigma^{-1}(p+q)}$ appara\^it dans $(id\otimes \ell_2)\circ\kappa(X,Y)$ avec le signe
\begin{align*}&(-1)^{x+v}\varepsilon\left(\begin{smallmatrix}x u v\\ u v x \end{smallmatrix}\right)(-1)^{a_{[1,p]}a_{[p+1,p+q]}}\varepsilon\left(\begin{smallmatrix}a_{[1,p+q]}\\ a_{[\sigma^{-1}(1),\sigma^{-1}(p+q)]}\end{smallmatrix}\right)(-1)\cr&=(-1)^{x+v+1}(-1)^{xa_{[p+1,p+q]}}(-1)^{a_{[1,p]}a_{[p+1,p+q]}}
\varepsilon\left(\begin{smallmatrix}a_{[1,p+q]}\\ a_{[\sigma^{-1}(1),\sigma^{-1}(p+q)]}\end{smallmatrix}\right)=(-1)^{x+u+1}
\varepsilon\left(\begin{smallmatrix}a_{[1,p+q]}\\ a_{[\sigma^{-1}(1),\sigma^{-1}(p+q)]}\end{smallmatrix}\right).
\end{align*}\hfill$\square$


\subsection{$G_\infty$ alg\`ebre }

\

\

On a construit deux cod\'erivations $m$ et $\ell$ de degr\'e $1$ pour la comultiplication $\Delta$ de la cog\`ebre cocommutative et coassociative $\big(S^+(\mathcal{H}[1]),\Delta\big)$ et pour le cocrochet $\kappa$ de la cog\`ebre de Lie $\big(S^+(\mathcal{H}[1]),\kappa\big)$ v\'erifiant $m\circ m=0$ et $\ell\circ\ell=0$. De plus, la comultiplication $\Delta$ et le cocrochet $\kappa$ v\'erifie l'identit\'e de coLeibniz $$(id\otimes\Delta)\circ\kappa=(\kappa\otimes id)\circ\Delta+\tau_{12}\circ(id\otimes\kappa)\circ\Delta.$$ On dit que $\big(S^+(\mathcal{H}[1]),\Delta,\kappa\big)$ est une cog\`ebre de Gerstenhaber gradu\'ee.

Comme $\ell+m$ est une cod\'erivation de degr\'e $1$ de $\big(S^+(\mathcal{H}[1]),\Delta,\kappa\big)$ v\'erifiant l'\'equation ma\^itresse $[\ell+m,\ell+m]=(\ell+m)^2=0$. Alors, $\big(S^+(\mathcal{H}[1]),\Delta,\kappa,\ell+m\big)$ est une cog\`ebre de Gerstenhaber diff\'erentielle gradu\'ee.

\begin{defn} {\rm ($G_\infty$ alg\`ebre)}

\

Une $G_\infty$ alg\`ebre est une bicog\`ebre de la forme $G(\mathcal G)=\left(S^+\left(\underline\bigotimes^+(\mathcal G[1])[1]\right),\Delta,\kappa\right)$ munie d'une cod\'erivation pour les deux structures $\Delta$ et $\kappa$ not\'ee $\ell+m$ et de carr\'e nul.

Soit $\mathcal G$ une alg\`ebre de Gerstenhaber, alors, $G(\mathcal G)=\left(S^+\left(\underline\bigotimes^+(\mathcal G[1])[1]\right),\Delta,\kappa,\ell+m\right)$ avec
$$
\ell_k=0~~\text{si }k\neq2, ~~~\ell_2(X,Y)=(-1)^x[X,Y],\quad m_k=0~~\text{si }k\neq1, ~~~ m_1(X)=\mu(X),
$$
$(X,Y\in\mathcal H=\underline\bigotimes^+(\mathcal G[1]))$ s'appelle la $G_\infty$ alg\`ebre enveloppante de l'alg\`ebre $\mathcal G$.
\end{defn}

\

\begin{rema}

\

Pour d\'efinir une cod\'erivation $Q$ de la bicog\`ebre $\left(S^+\left(\underline\bigotimes^+(\mathcal G[1])[1]\right),\Delta,\kappa\right)$, il suffit de se donner une suite d'applications
$$
Q^{(N)}_{p_1\dots p_n}:\underline\bigotimes^{p_1}(\mathcal G[1])[1]\pt\dots\pt\underline\bigotimes^{p_n}(\mathcal G[1])[1]\longrightarrow\mathcal G
\qquad(p_1+\dots+p_n=N).
$$
(cette construction est explicit\'ee dans la section suivante pour des morphismes de bicog\`ebres).

Une $G_\infty$ alg\`ebre est la $G_\infty$ alg\`ebre enveloppante d'une alg\`ebre de Gerstenhaber si et seulment si sa cod\'erivation $\ell+m$ v\'erifie $(\ell+m)^{(N)}_{p_1\dots p_n}=0$  sauf pour $(\ell+m)^{(2)}_{11}$ et $(\ell+m)^{(2)}_2$.
\end{rema}

\

\section{Cohomologie de Chevalley-Harrison des alg\`ebres de Gerstenhaber}\label{3}

\

\subsection{Morphismes de cog\`ebres entre $G_\infty$ alg\`ebres}

\

Rappelons que si $S^+(\mathfrak g[1])$ et $S^+(\mathfrak g'[1])$ sont deux $L_\infty$ alg\`ebres, resp.si $\underline\bigotimes^+(A[1])$ et $\underline\bigotimes(A'[1])$ sont deux $C_\infty$ alg\`ebres, un morphisme $F$ de $L_\infty$ alg\`ebre, resp de $C_\infty$ alg\`ebre entre ces deux alg\`ebres est un morphisme de cog\`ebre qui commute avec les codiff\'erentielles $\ell^{\mathfrak g}$ et $\ell^{\mathfrak g'}$, resp. les codiff\'erentielles $m^A$ et $m^{A'}$. De plus, les morphismes de cog\`ebres sont caract\'eris\'es par leurs projections $F_n$
$$
F_n:S^n({\mathfrak g}[1])\longrightarrow\mathfrak h\qquad\text{resp.}\quad F_n:\underline\bigotimes^n(A[1])\longrightarrow A'.
$$

Dans le cas de $G_\infty$ alg\`ebres, $S^+\left(\underline\bigotimes^+(\mathcal G[1])[1]\right)$ et  $S^+\left(\underline\bigotimes^+(\mathcal G'[1])[1]\right)$, on dispose de deux coproduits: $\Delta$ et $\kappa$.

\vskip 0.3cm
\begin{prop} {\rm (Les morphismes de cog\`ebres entre deux $G_\infty$-alg\`ebres)}

\

Un morphisme de cog\`ebre $F:S^+\left(\underline\bigotimes^+(\mathcal G[1])[1]\right)\longrightarrow S^+\left(\underline\bigotimes^+(\mathcal G'[1])[1]\right)$, c'est \`a dire une application lin\'eaire telle que
$$
(F\otimes F)\circ\Delta=\Delta\circ F\qquad\text{et}\quad(F\otimes F)\circ\kappa=\kappa\circ F,
$$
est uniquement caract\'eris\'e par ses projections dans $\mathcal G'$, que l'on note:
$$
f_{p_1\dots p_n}:\underline\bigotimes^{p_1}(\mathcal G[1])\pt\ldots\pt\underline\bigotimes^{p_n}(\mathcal G[1])\longrightarrow\mathcal G'.
$$
\end{prop}

\vskip 0.3cm

Nous ne d\'emontrerons pas cette proposition annonc\'ee dans \cite{[BGHHW]} et \cite{[GH]}. Nous allons seulement expliquer comment reconstruire $F$ \`a partir des $f_{p_1\dots p_n}$. On vient de voir qu'il suffit de retrouver les projections
$$
F_n:S^n\left(\underline\bigotimes^+(\mathcal G[1])[1]\right)\longrightarrow\underline\bigotimes^+(\mathcal G'[1])
$$
pour caract\'eriser $F$. On va donc se contenter de d\'ecrire la construction des applications $F_n$ \`a partir des applications $f_{p_1\dots p_n}$.

Soit donc $X_1\pt\dots\pt X_n$ un \'el\'ement de $S^n\left(\underline\bigotimes^+(\mathcal G[1])[1]\right)$, avec
$$
X_j=\alpha_1^j\underline\otimes\dots\underline\otimes\alpha^j_{p_j}.
$$
$F_n(X_1\pt\dots\pt X_n)$ est une somme de produits tensoriels modulo les battements de $f(Y_k)$ ($1\leq k\leq t$) o\`u les $Y_k$ sont des produit $\pt$ de parties des $X_j$ de la forme:
$$
U_i^j=\alpha^j_{r_i+1}\underline\otimes\alpha^j_{r_i+2}\underline\otimes\dots\underline\otimes\alpha^j_{r_{i+1}}\qquad(0\leq r_s\leq p_j).
$$
On envoie donc les produits $(U^1_1\underline\otimes\dots\underline\otimes U_{r_1}^1)\pt\dots\pt(U^n_1\underline\otimes\dots\underline\otimes U_{r_n}^n)$ sur des sommes de termes de la forme
$$
f(Y_1)\underline\otimes\dots\underline\otimes f(Y_t)=f(V_1^1\pt\dots\pt V^1_{s_1})\underline\otimes\dots\underline\otimes f(V^t_1\pt\dots\pt V_{s_t}^t).
$$
Les $V^k_\ell$ forment une permutation des $U_i^j$. Si un $X_j$ n'est pas d\'ecompos\'e ($r_j=1$), il ne peut appara\^\i tre qu'en facteur d'au moins une vraie partie $U^{j'}_i$ d'un autre $X$ ($r_{j'}>1$). Si un $X_j$ est d\'ecompos\'e ($r_j>1$) chacune de ses parties appara\^\i t dans un $Y$ diff\'erent. Enfin, il y a autant de $\pt$ et de $\underline\otimes$ dans l'expression de $X_1\pt\dots\pt X_n$ que dans celle des $f(Y_1)\underline\otimes\dots\underline\otimes f(Y_t)$.\\

\noindent
{\bf Etape 1:} Description des d\'ecoupages des $X_j$.\\

On utilise un tableau $T$ \`a $n$ lignes. On d\'ecoupe les $X_j$. Chaque ligne $j$ du tableau $T$ a $r_j$ cases. Dans chaque case, on place les portions de $X_j$ dans l'ordre du d\'ecoupage ainsi: si on \'ecrit $X_j=U_1^j\underline{\otimes}\dots\underline{\otimes}U_{r_j}^j$, on aura une ligne:

\begin{displaymath}
\begin{array}{c}
\begin{array}{|l|l|l|l|}
\hline
\makebox(15,15){$U_1^j$}&\makebox(15,15){$U_2^j$}&\makebox(15,15){$\dots$}&\makebox(15,15){$U_{r_j}^j$}\\
\hline
\end{array}\end{array}
\end{displaymath}
\vskip0.2cm

Notons $L=$ nombre de lignes, $C=$ nombre de colonnes, $N=$ nombre de cases de $T$.\vskip0.2cm

 Par exemple si on envisage de d\'ecomposer $X_1\pt X_2\pt X_3\pt X_4\pt X_5\pt X_6$ en
$$X_1\pt(U_2\underline\otimes V_2)\pt X_3\pt(U_4\underline\otimes V_4\underline\otimes W_4)\pt(U_5\underline\otimes V_5\underline\otimes W_5\underline\otimes R_5)\pt(U_6\underline\otimes V_6),
$$
on posera:
\begin{displaymath}
T~=~\begin{array}{c}
\begin{array}{|c|}
\hline
\makebox(15,15){$X_1$}\\
\end{array}\hfill\\
\begin{array}{|c|c|}
\hline
\makebox(15,15){$U_2$}&\makebox(15,15){$V_2$}\\
\hline
\end{array}\hfill\\
\begin{array}{|l|}
\makebox(15,15){$X_3$}\\
\end{array}\hfill\\
\begin{array}{|c|c|c|}
\hline
\makebox(15,15){$U_4$}&\makebox(15,15){$V_4$}&\makebox(15,15){$W_4$}\\
\end{array}\hfill\\
\begin{array}{|l|l|l|l|}
\hline
\makebox(15,15){$U_5$}&\makebox(15,15){$V_5$}&\makebox(15,15){$W_5$}&\makebox(15,15){$R_5$}\\
\hline
\end{array}\\
\begin{array}{|l|l|}
\makebox(15,15){$U_6$}&\makebox(15,15){$V_6$}\\
\hline
\end{array}\hfill\\
\end{array}
\end{displaymath}
Le tableau $T$ caract\'erise maintenant la d\'ecomposition de $X_1\pt\dots\pt X_n$. \\

\noindent{\bf Etape 2:} Suppression des lignes de longueur $1$\\

On appelle $T_1$ le sous-tableau obtenu en enlevant toutes les lignes de longueur 1 de $T$. Il est clair que l'on peut reconstruire $T$ \`a partir de $T_1$. Celui-ci caract\'erise donc aussi la d\'ecomposition de notre mon\^ome. Dans notre exemple:
\begin{displaymath}
T_1~=~\begin{array}{c}
\begin{array}{|c|c|}
\hline
\makebox(15,15){$U_2$}&\makebox(15,15){$V_2$}\\
\end{array}\hfill\\
\begin{array}{|c|c|c|}
\hline
\makebox(15,15){$U_4$}&\makebox(15,15){$V_4$}&\makebox(15,15){$W_4$}\\
\end{array}\hfill\\
\begin{array}{|l|l|l|l|}
\hline
\makebox(15,15){$U_5$}&\makebox(15,15){$V_5$}&\makebox(15,15){$W_5$}&\makebox(15,15){$R_5$}\\
\hline
\end{array}\\
\begin{array}{|l|l|}
\makebox(15,15){$U_6$}&\makebox(15,15){$V_6$}\\
\hline
\end{array}\hfill\\
\end{array}
\end{displaymath}

\

Notons $L_1=$ nombre de lignes, $C_1=$ nombre de colonnes, $N_1=$ nombre de cases de $T_1$ et $h=L_1-1$.\vskip0.2cm

\noindent
{\bf Etape 3:} Construction du coeur de $T'$.\\

On dessine tous les tableaux $T_2$ ayant $h$ cases vides tels que la longueur des lignes d\'ecro\^it de haut vers le bas.
Dans notre exemple, ce nombre est 3. Dans notre exemple, il y a quatre possibilit\'es:
\begin{displaymath}\begin{array}{ccccc}
T_2~=~\begin{array}{c}
\begin{array}{|c|c|c|}
\hline
\hbox{  }&\hbox{  }&\hbox{  }\\
\hline
\end{array}
\end{array}&\quad&
T_2~=~\begin{array}{c}
\begin{array}{|c|c|}
\hline
\hbox{  }&\hbox{  }\\
\hline
\end{array}\\
\begin{array}{|c|}
\hbox{  }\\
\hline
\end{array}\hfill
\end{array}&\quad&
T_2~=~\begin{array}{c}
\begin{array}{|c|c|}
\hline
\hbox{  }&\hbox{ }\\
\hline
\end{array}\\ \hskip0.47cm
\begin{array}{|c|}
\hbox{  }\\
\hline
\end{array}\hfill
\end{array}\quad\quad
T_2~=~\begin{array}{|c|}
\hline
\hbox{  }\\
\hline
\hbox{  }\\
\hline
\hbox{  }\\
\hline
\end{array}
\end{array}
\end{displaymath}
On ajoute ensuite  une case vide \`a $T_2$ en dessous de chacune des colonnes de $T_2$ (s'il n'y a aucune colonne, on ne fait rien).
 On obtient  un tableaux $T_3$. Par exemple si on retient le second $T_2$ ci dessus, on obtient,
 \begin{displaymath}
\begin{array}{ccc}
T_2~=~\begin{array}{c}
\begin{array}{|c|c|}
\hline
\hbox{  }&\hbox{  }\\
\hline
\end{array}\\
\begin{array}{|c|}
\hbox{  }\\
\hline
\end{array}\hfill
\end{array}&~\Longrightarrow~&
T_3~=~\begin{array}{c}
\begin{array}{|c|c|}
\hline
\hbox{  }&\hbox{  }\\
\hline
\end{array}\\
\begin{array}{|c|c|}
\hbox{  }&\hbox{  }\\
\hline
\end{array}\hfill\\
\begin{array}{|c|}
\hbox{  }\\
\hline
\end{array}\hfill\\
\end{array}
\end{array}
\end{displaymath}

\vskip0.2cm

\noindent
{\bf Etape 4:} Construction de $T_4$ vide\\

 On ajoute \`a la premi\`ere ligne du tableau ainsi obtenu autant de cases qu'il faut pour obtenir un tableau $T_4$ ayant le m\^eme nombre de cases que $T_1$, c'est \`a dire tel que $N_4=N_1$.

\begin{rema}

\

On a ajout\'e $N_1-N_3=N_1-(h+C_2)$ cases \`a la premi\`ere ligne de $T_3$. Alors, $T_4$ a donc $C_4$ colonnes o\`u
\begin{align*}C_4&=C_3+N_1-(h+C_2)\\&=N_1-h=N-h-\#\{X_j \ \hbox{non d\'ecompos\'es}\}\end{align*}
\end{rema}

 Par exemple si on retient le second $T_3$ ci dessus, on choisit, puisque $T_1$ poss\`ede 11 cases:
\begin{displaymath}
\begin{array}{ccc}
T_3~=~~\begin{array}{c}
\begin{array}{|c|c|}
\hline
\hbox{  }&\hbox{  }\\
\hline
\end{array}\\
\begin{array}{|c|c|}
\hbox{  }&\hbox{  }\\
\hline
\end{array}\hfill\\
\begin{array}{|c|}
\hbox{  }\\
\hline
\end{array}\hfill\\
\end{array}&~\Longrightarrow~&
T_4~=~\begin{array}{c}
\begin{array}{|c|c|c|c|c|c|c|c|}
\hline
\hbox{  }&\hbox{  }&\hbox{  }&\hbox{  }&\hbox{  }&\hbox{  }&\hbox{  }&\hbox{  }\\
\hline
\end{array}\\
\hskip0.47cm\begin{array}{|c|c|}
\hbox{  }&\hbox{  }\\
\hline
\end{array}\hfill\\
\hskip0.47cm\begin{array}{|c|}
\hbox{  }\\
\hline
\end{array}\hfill\\
\end{array}
\end{array}
\end{displaymath}

\noindent
{\bf Etape 5:} Remplissage de $T_4$\\

On remplit ensuite $T_4$ en mettant dans chaque case une des lettres de $T_1$ ainsi:

On lit $T_1$, ligne par ligne, de la gauche vers la droite. Pour la ligne $j$ de $T_1$, de longueur $r_j$, on choisit $r_j$ colonnes de $T_4$ telles qu'il y ait dans chaque colonne au moins une case vide. Disons que ces colonnes sont $c_{i_1},\dots,c_{i_{r_j}}$ avec $i_1<\dots<i_{r_j}$. On place la $k^{i\grave eme}$ entr\'ee $U_k^j$ de la ligne $j$ de $T_1$ dans la premi\`ere case libre de la colonne n$^\circ$ $i_k$ lorsqu'on parcourt la colonne de haut en bas. On obtient ainsi un tableau rempli $T_4$. Par exemple:
\begin{displaymath}
T_4~=~\begin{array}{c}
\begin{array}{|l|l|l|l|l|l|l|l|}
\hline
\makebox(15,15){$U_2$}&\makebox(15,15){$U_4$}&\makebox(15,15){$V_2$}&\makebox(15,15){$W_4$}&\makebox(15,15){$V_5$}&\makebox(15,15){$W_5$}&\makebox(15,15){$R_5$}&\makebox(15,15){$V_6$}\\
\hline
\end{array}\\
\hskip0.87cm\begin{array}{|l|l|}
\makebox(15,15){$U_5$}&\makebox(15,15){$V_4$}\\
\hline
\end{array}\hfill\\
\hskip0.87cm\begin{array}{|l|}
\makebox(15,15){$U_6$}\\
\hline
\end{array}\hfill\\
\end{array}
\end{displaymath}

\noindent
{\bf Etape 6:} Ajout des $X_j$ ind\'ecompos\'es

On ajoute des cases au tableau $T_4$ en dessous des colonnes existante de $T_4$ (exactement $N-N_4$ cases). On remplit ces cases en mettant les $X_j=U_1^j$ ind\'ecompos\'es de $T$. On obtient un tableau $T_5$. Par exemple:
\begin{displaymath}
T_4 ~\Longrightarrow~ T_5~=~\begin{array}{c}
\begin{array}{|l|l|l|l|l|l|l|l|}
\hline
\makebox(15,15){$U_2$}&\makebox(15,15){$U_4$}&\makebox(15,15){$V_2$}&\makebox(15,15){$W_4$}&\makebox(15,15){$V_5$}&\makebox(15,15){$W_5$}&\makebox(15,15){$R_5$}&\makebox(15,15){$V_6$}\\
\hline
\end{array}\\
\hskip0.87cm\begin{array}{|l|l|l|}
\makebox(15,15){$U_5$}&\makebox(15,15){$V_4$}&\makebox(15,15){$X_1$}\\
\hline
\end{array}\hfill\\
\hskip0.87cm\begin{array}{|l|l|}
\makebox(15,15){$U_6$}&\makebox(15,15){$X_3$}\\
\hline
\end{array}\hfill\\
\end{array}
\end{displaymath}\vskip0.25cm

\noindent
{\bf Etape 7:} D\'efinition de $T'$

 On range chaque colonne du tableau $T_5$ ainsi obtenu dans l'ordre croissant du haut vers le bas. C'est \`a dire, on construit des colonnes $c'_i$ de la forme:

\begin{displaymath}\begin{array}{ccccc}
 \begin{array}{|c|}
\hline
\makebox(15,13){ $U_{a_1}^{j_1}$ }\\
\hline
\makebox(15,13){ $U_{a_2}^{j_2}$ }\\
\hline
\makebox(15,13){ $\vdots$ }\\
\hline
\makebox(15,13){ $U_{a_s}^{j_t}$ }\\
\hline
\end{array}
\end{array} \hskip1cm \hbox{avec $j_1<\dots<j_t.$}
\end{displaymath}

 On obtient ainsi un tableau not\'e $T'$.

Dans notre exemple:
\begin{displaymath}
T'~=~\begin{array}{c}
\begin{array}{|c|c|c|c|c|c|c|c|}
\hline
\makebox(15,15){$U_2$}&\makebox(15,15){$U_4$}&\makebox(15,15){$V_2$}&\makebox(15,15){$X_1$}&\makebox(15,15){$V_5$}&\makebox(15,15){$W_5$}
&\makebox(15,15){$R_5$}&\makebox(15,15){$V_6$}\\
\hline
\end{array}\\
\hskip0.88cm\begin{array}{|c|c|c|}
\makebox(15,15){$U_5$}&\makebox(15,15){$X_3$}&\makebox(15,15){$W_4$}\\
\hline
\end{array}\hfill\\
\hskip0.88cm\begin{array}{|c|c|}
\makebox(15,15){$U_6$}&\makebox(15,15){$V_4$}\\
\hline
\end{array}\hfill\\
\end{array}
\end{displaymath}

\

\noindent
{\bf Etape 8:} Calcul des $F(T')$.\\

A chaque colonne $c_i'$ de $T'$ ($1\leq i\leq s$), on associe un \'el\'ement $f(c_i')$ de $\mathcal G'$ qui est simplement l'image par $f$ du produit $\pt$ des termes de la colonne. Finalement on d\'efinit $F(T')$ comme le produit $\varepsilon(T,T')\underline\otimes_{i=1}^sf(c_i')$ modulo les battements avec le signe obtenu \`a partir de la permutation passant de $T$ \`a $T'$.

Dans notre exemple, on pose donc
$$
\aligned
f(Y_1)&=f(U_2)\\
f(Y_2)&=f(U_4\pt U_5\pt U_6)\\
f(Y_3)&=f(V_2\pt X_3\pt V_4)\\
f(Y_4)&=f(X_1\pt W_4)\\
f(Y_5)&=f(V_5)\\
f(Y_6)&=f(W_5)\\
f(Y_7)&=f(R_5)\\
f(Y_8)&=f(V_6).
\endaligned
$$
Donc $$
F(T')=\varepsilon(T,T')f(Y_{1})\underline\otimes f(Y_{2})\underline\otimes\dots\underline\otimes f(Y_{8}).
$$

\

\noindent
{\bf Etape 9:} Calcul de $F(X_1\pt\dots\pt X_n)$.\\

La quantit\'e $F(X_1\pt\dots\pt X_n)$ s'obtient en faisant la somme pour des $F(T')$ toutes les d\'ecompositions $T$ de $X_1\pt\dots\pt X_n$ et pour chaque $T$ pour tous les tableaux $T'$ qu'on peut construire \`a partir de $T$ :
$$
F(X_1\pt\dots\pt X_n)=\sum_{T,T'}\varepsilon( T,T')f(Y_{1})\underline\otimes f(Y_{2})\underline\otimes\dots\underline\otimes f(Y_{C_{T'}}).
$$

\

\subsection{Cohomologie de Chevalley-Harrison}

\

\

On a la d\'efinition naturelle:

\vskip 0.3cm
\begin{defn}{\rm (Morphismes de $G_\infty$ alg\`ebres)}

\

Soit $\left(S^+\left(\underline\bigotimes^+(\mathcal G[1])[1]\right),\Delta,\kappa,\ell+m\right)$ et $\left(S^+\left(\underline\bigotimes^+(\mathcal G'[1])[1]\right),\Delta',\kappa',\ell'+m'\right)$ deux $G_\infty$ alg\`ebres. Une application $F: S^+\left(\underline\bigotimes^+(\mathcal G[1])[1]\right)\longrightarrow S^+\left(\underline\bigotimes^+(\mathcal G'[1])[1]\right)$ est un morphisme de $G_\infty$ alg\`ebres si $F$ est un morphisme de cog\`ebres:
$$
(F\otimes F)\circ\Delta=\Delta'\circ F,\qquad (F\otimes F)\circ\kappa=\kappa'\circ F,
$$
de degr\'e $0$ qui pr\'eserve les codiff\'erentielles:
$$
F\circ(\ell+m)=(\ell'+m')\circ F.
$$
\end{defn}

\

Soient $(\mathcal G,\wedge,[~,~])$ et $(\mathcal G',\wedge',[~,~]')$ deux alg\`ebres de Gerstenhaber. Comme pr\'ec\'e\-dem\-ment, cherchons \`a construire un morphisme de $G_\infty$ alg\`ebres $F$. On a vu qu'en tant que morphisme de cog\`ebres, $F$ est caract\'eris\'e par la suite des applications $(f_{p_1\dots p_r})$. Donnons-nous un morphisme d'alg\`ebres de Gerstenhaber $f_1:\mathcal G\longrightarrow\mathcal G'$ et cherchons \`a construire des applications $f_{p_1\dots p_n}$ suivantes. Dans cette partie, nous noterons ces applications $f^{(N)}_{p_1\dots p_n}$ si $\sum p_j=N$.

Supposons connues tous les $f^{(k)}_{p_1\dots p_r}$ avec $k=p_1+\dots+p_r<N$, on cherche les $f^{(N)}_{p'_1\dots p'_n}$.\\

Si $X_j$ est un \'el\'ement de $\underline\bigotimes^{p_j}(\mathcal G[1])[1]$, on notera aussi $p(X_j)$ le nombre $p_j$.\\

Si on applique $(\ell'+m')\circ F-F\circ(\ell+m)$ \`a $X_1\pt\dots\pt X_n$ avec $X_j\in\underline\bigotimes^{p(X_j)}(\mathcal G[1])[1]$ et $\sum_j p(X_j)\leq N$, aucun terme en $f^{(N)}_{p'_1\dots p'_r}$ n'appara\^\i t.\\

Ces termes apparaissent lorsqu'on applique $(\ell'+m')\circ F- F\circ(\ell+m)$ \`a $X_1 \pt\dots \pt X_n$ avec $\sum_j p(X_j)=N+1$. On trouve les termes suivants:

\

Dans $(F\circ\ell)(X_1\pt\dots\pt X_n)$, on a:
$$
\begin{aligned}
&(F\circ\ell)(X_1 \pt\dots \pt X_n)\\
&=\sum_{j<k}\varepsilon\left(\begin{smallmatrix}x_1~~\dots~~x_n\\ x_jx_kx_1\dots\hat{\jmath}\dots\hat{k}\dots x_n\end{smallmatrix}\right)
F\big(\ell_2(X_j,X_k)\pt X_1\pt\dots\hat{\jmath}\dots\hat{k}\dots\pt X_n\big)\\
&=\sum_{j<k}\varepsilon\left(\begin{smallmatrix}x_1~~\dots~~x_n\\ x_jx_kx_1\dots\hat{\jmath}\dots\hat{k}\dots x_n\end{smallmatrix}\right)
f^{(N)}_{p(X_j)+p(X_k)-1,p(X_1),\dots\hat{\jmath}\dots\hat{k}\dots,p(X_n)}\big(\ell_2(X_j,X_k)\pt X_1\pt\dots\hat{\jmath}\dots\hat{k}\dots\pt X_n\big)
\end{aligned}
$$

Dans $( F\circ m)(X_1\pt\dots\pt X_n)$, on a:
$$
\begin{aligned}
(F\circ m)&(X_1\pt\dots\pt X_n)=\sum_{j/p(X_j)>1}(-1)^{\sum_{i<j}x_i} F(X_1\pt\dots\pt m_1(X_j)\pt\dots\pt X_n)\\
&=\sum_{j/p(X_j)>1}(-1)^{\sum_{i<j}x_i}f_{p(X_1),\dots,p(X_j)-1,\dots,p(X_n)}^{(N)}(X_1\pt\dots\pt m_1(X_j)\pt\dots\pt X_n)
\end{aligned}
$$

Par ailleurs, dans $(\ell'\circ F)(X_1\pt\dots\pt X_n)$, les termes en $f^{(N)}$ n'apparaissent que dans les termes d'ordre 2 du d\'eveloppement de $F$. Avec les notations ci-dessus, il faut, en effet, au moins un produit de deux paquets et que l'un contienne $N$ vecteurs de $\mathcal G'$. Il reste seulement:
$$
\begin{aligned}
(\ell'\circ F)(X_1\pt\dots\pt X_n)&=\sum_{j/p(X_j)=1}\ell'\Big(\varepsilon\left(\begin{smallmatrix}x_1~~\dots~~x_n\\ x_jx_1\dots \hat{\jmath}\dots x_n\end{smallmatrix}\right) f_1(X_j) \pt f_{p(X_1),\dots\hat\jmath\dots,p(X_n)}^{(N)}(X_1\pt\dots\hat\jmath\dots\pt X_n)\Big)\\
&+\sum_{r<N}~\text{termes en }f^{(r)}.
\end{aligned}
$$

De m\^eme pour $(m'\circ F)(X_1 \pt\dots \pt X_n)$, les seuls termes en $f^{(N)}$ apparaissant sont:
$$
\begin{aligned}
&(m'\circ F)(X_1\pt\dots\pt X_n)\\
&=\sum_{j/p(X_j)>1}\varepsilon\left(\begin{smallmatrix}x_1~~\dots~~x_n\\ x_jx_1\dots \hat{\jmath}\dots x_n\end{smallmatrix}\right) m'\Big(f^{(1)}_1(\alpha^j_1)\underline\otimes f_{p_j-1,p_1,\dots\hat\jmath\dots,p_r}^{(N)} ((\alpha^j_2\underline\otimes \dots\underline\otimes\alpha_{p_j}^j)\pt X_1\pt\dots\hat\jmath\dots\pt X_n)\Big)\\
&+\varepsilon\left(\begin{smallmatrix}x_1~~\dots~~x_n\\ x_1\dots \hat{\jmath}\dots x_nx_j\end{smallmatrix}\right) m'\Big(f_{p_1,\dots\hat\jmath\dots,p_r,p_j-1}^{(N)}(X_1\pt\dots\hat\jmath\dots\pt X_n \pt(\alpha^j_1\underline\otimes\dots\underline\otimes\alpha_{p_j-1}^j))\underline\otimes f^{(1)}_1(\alpha^j_{p_j})\Big)\\
&+\sum_{r<N}~\text{termes en }f^{(r)}.
\end{aligned}
$$

L'op\'erateur de cobord $d_{CH}$ dit de Chevalley-Harrison associ\'e \`a cette construction de $ F$ et \`a $f_1^{(1)}$ est donc le suivant.

L'espace des cocha\^ines est l'espace
$$
C^N=\sum_{p_1+\dots+p_n=N}Hom\left(\underline\bigotimes^{p_1}(\mathcal G[1])[1]\pt\dots\pt\underline\bigotimes^{p_n}(\mathcal G[1])[1],\mathcal G'[1-N-n]\right).
$$

Le cobord $d_{CH}:C^N\longrightarrow C^{N+1}$ est de la forme $d_{CH}=d_m+d_\ell$ avec
$$
\begin{aligned}
&(d_mf^{(N)}_{p_1\dots p_n})(X_1\pt\dots\pt X_n)=(d_mf)^{(N+1)}_{p_1\dots(p_j+1)\dots p_n}(X_1\pt\dots\pt X_n)\\
&=(-1)^{\alpha^j_1\sum_{i<j}x_i} m'\left(f^{(1)}_1(\alpha^j_1)\underline\otimes f_{p_1,\dots,p_j-1,\dots,p_r}^{(N)}(X_1 \pt\dots\pt(\alpha^j_2\underline\otimes\dots\underline\otimes\alpha_{p_j}^j)\pt\dots\pt X_n)\right)\\
&+(-1)^{\alpha^j_{p_j}.\sum_{i>j}x_i} m'\left(f_{p_1,\dots,p_j-1,\dots,p_r}^{(N)}(X_1\pt\dots\pt(\alpha^j_1\underline\otimes\dots\underline\otimes\alpha_{p_j-1}^j)\dots\pt X_n)\underline\otimes f^{(1)}_1(\alpha^j_{p_j})\right)\\
&-(-1)^{\sum_{i<j}x_i}f^{(N)}_{p_1\dots p_n}(X_1\pt\dots\pt m(X_j)\pt\dots\pt X_n).
\end{aligned}
$$

De m\^eme, $d_\ell:C^N_{p_1\dots p_n}\longrightarrow\displaystyle\sum_{j,~q_1+q_2=p_j+1}C^{N+1}_{q_1,q_2,p_1\dots\hat\jmath\dots p_n}$ s'\'ecrit
$$
(d_\ell f^N_{p_1\dots p_n})=\sum_{\begin{smallmatrix}j\\ q_1+q_2=p_j+1\end{smallmatrix}}(d_\ell f)^N_{q_1,q_2,p_1\dots \hat\jmath\dots p_n}.
$$
Avec
\begin{itemize}
\item[1.] Si $q_1>1$ et $q_2>1$, alors
$$
\begin{aligned}
&(d_\ell f)^{(N+1)}_{q_1,q_2,p_1\dots\hat\jmath\dots p_n}(Y_1\pt Y_2\pt X_1\pt\dots\hat\jmath\dots\pt X_n)=\\
&-\varepsilon\left(\begin{smallmatrix}y_1y_2x_1\dots\hat\jmath\dots x_n\\ x_1\dots y_1y_2\dots x_n\end{smallmatrix}\right)(-1)^{\sum_{i<j}x_i}f^{(N)}_{p_1\dots p_n}(X_1\pt\dots\pt\ell(Y_1,Y_2)\pt\dots\pt X_n).
\end{aligned}
$$
\item[2.] Si $q_1=1$ et $q_2=p_j>1$, alors
$$
\begin{aligned}
&(d_\ell f)^{(N+1)}_{q_1,q_2,p_1\dots\hat\jmath\dots p_n}(Y_1\pt Y_2\pt X_1\pt\dots\hat\jmath\dots\pt X_n)=\\
&-\varepsilon\left(\begin{smallmatrix}y_1y_2x_1\dots\hat\jmath\dots x_n\\ x_1\dots y_1y_2\dots x_n\end{smallmatrix}\right)(-1)^{\sum_{i<j}x_i} f^{(N)}_{p_1\dots p_n}(X_1\pt\dots\pt\ell(Y_1,Y_2)\pt\dots\pt X_n)\\
&+\varepsilon\left(\begin{smallmatrix}y_1y_2x_1\dots\hat\jmath\dots x_n\\ y_1x_1\dots y_2\dots x_n\end{smallmatrix}\right) \ell'\Big(f_1^{(1)}(Y_1),f^{(N)}_{p_1\dots p_N}(X_1\pt\dots\pt Y_2\pt\dots\pt X_n)\Big).
\end{aligned}
$$
\item[3.] On a la m\^eme formule par sym\'etrie si $q_2=1$ et $q_1=p_j>1$.
\item[4.] Enfin, si $q_1=q_2=1$, alors
$$
\begin{aligned}
&(d_\ell f)^{(N+1)}_{1,1,p_1\dots\check{j}\dots p_n}(Y_1\pt Y_2\pt X_1\pt\dots\hat\jmath\dots\pt X_n)=\\
&\varepsilon\left(\begin{smallmatrix}y_1y_2x_1\dots\hat\jmath\dots x_n\\ y_1x_1\dots y_2\dots x_n\end{smallmatrix}\right)
\ell'\left(f_1^{(1)}(Y_1)\pt f^{(N)}_{p_1\dots p_N}(X_1\pt\dots\pt Y_2\pt\dots\pt X_n)\right)\\
&+\varepsilon\left(\begin{smallmatrix}y_1y_2x_1\dots\hat\jmath\dots x_n\\ x_1\dots y_1\dots x_ny_2\end{smallmatrix}\right)
\ell'\left(f^{(N)}_{p_1\dots p_N}(X_1\pt\dots\pt Y_1\pt\dots\pt X_n)\pt f_1^{(1)}(Y_2)\right)\\
&-(-1)^{\sum_{i<j}x_i}\varepsilon\left(\begin{smallmatrix}y_1y_2x_1\dots\hat\jmath\dots x_n\\ x_1\dots y_1y_2\dots x_n\end{smallmatrix}\right) f^{(N)}_{p_1\dots p_n}(X_1\pt\dots\pt\ell(Y_1\pt Y_2)\pt\dots\pt X_n).
\end{aligned}
$$
\end{itemize}

\

\begin{prop} {\rm ($d_{CH}$ est un op\'rateur de cobord)}

\

Soit $f^{(1)}_1:\mathcal G\longrightarrow\mathcal G'$ un morphisme d'alg\`ebres de Gerstenhaber. Alors
\begin{itemize}
\item[(i)] Pour tout $N$, $f^{(1)}_1+\sum_{p_1+\dots+p_n=N}f^{(N)}_{p_1\dots p_n}$ d\'efinit un morphisme de $G_\infty$ alg\`ebres \`a l'ordre $N+1$ de $G(\mathcal G)$ dans $G(\mathcal G')$ si et seulement si:
$$
d_{CH}\left(\sum_{p_1+\dots+p_n=N}f^{(N)}_{p_1\dots p_n}\right)=0
$$
\item[(ii)] Pour tout $g=\sum_{p_1+\dots+p_n=N-1}g^{(N-1)}_{p_1\dots p_n}$, $f+d_{CH}g$ est un morphisme \`a l'ordre $N+1$.
\item[(iii)] On a donc
$$
d_{CH}\circ d_{CH}=0.
$$
\end{itemize}
\end{prop}


\

\section{Un exemple}

\

\

Dans cette section, nous allons montrer comment le premier cocycle fondamental de la cohomologie de Chevalley des champs de vecteurs \`a valeurs dans les fonctions survit dans la cohomologie de Chevalley-Harrison d'une sous alg\`ebre de Gerstenhaber naturelle de $T_{poly}(\mathbb{R}^d)$ \`a valeurs dans le corps de base.

On note $\mathcal G=T_{poly}^{hom}(\mathbb{R}^d)$ l'espace des tenseurs totalement antisym\'etriques
$$
\alpha=\sum_{i_1<\dots<i_k}\alpha^{i_1\dots i_k}\partial_{i_1}\wedge\dots\wedge\partial_{i_k}
$$
tels que chaque coefficient $\alpha^{i_1\dots i_k}$ est un polyn\^ome homog\`ene de degr\'e $k$.

$T_{poly}^{hom}(\mathbb{R}^d)$ est une sous alg\`ebre de Gerstenhaber de $T_{poly}(\mathbb{R}^d)$. En effet, si $\alpha, \beta\in T_{poly}^{hom} (\mathbb{R}^d)$, alors, $\alpha\wedge\beta$ et $[\alpha,\beta]_S$ sont encore des tenseurs homog\`enes et ils appartiennent \`a $T_{poly}^{hom} (\mathbb{R}^d)$.

D'autre part, $\mathcal G'=\mathbb{R}$ muni de la multiplication usuelle $\alpha\wedge\beta=\alpha\beta$ et du crochet nul $[\alpha,\beta]=0$ est une alg\`ebre de Gerstenhaber pour la graduation $degr\acute{e}(\alpha)=0$, quel que soit $\alpha\in\mathbb{R}$. L'application
$$
f_1^1:T_{poly}^{hom}(\mathbb{R}^d)\longrightarrow\mathbb{R}
\quad\alpha\longmapsto\left\{
   \begin{array}{ll}
   \alpha, & \hbox{si $\alpha\in \big(T_{poly}^{hom}(\mathbb{R}^d)\big)^{0}$ ;} \\
   0, & \hbox{sinon.}
   \end{array}\right.
$$
est un morphisme d'alg\`ebres de Gerstenhaber.

Ceci d\'efinit donc une cohomologie de Chevalley-Harrison sur les espaces
$$
C^N_{p_1\dots p_n}=Hom\left(\underline\bigotimes^{p_1}(\mathcal G[1])[1]\pt\dots\pt\underline\bigotimes^{p_n}(\mathcal G[1])[1],\mathcal G'[1-N-n]\right).
$$

On sait dans \cite{[AAC1]} ou \cite{[AAC2]}, que le premier cocycle fondamental de Chevalley pour les champs de vecteurs ou les tenseurs lin\'eaires est un $3$-cocycle qui s'\'ecrit:
$$
f(\alpha_1,\alpha_2,\alpha_3)=\partial_{i_3}\alpha_1^{i_1}\partial_{i_1}\alpha_2^{i_2}\partial_{i_2}\alpha_3^{i_3}
-\partial_{i_2}\alpha_1^{i_1}\partial_{i_1}\alpha_3^{i_3}\partial_{i_3}\alpha_2^{i_2}.
$$

Consid\`erons donc l'application $f^3_{111}\in C^3_{111}$ d\'efinie ainsi:
$$
f^3_{111}\big((\alpha_1)\pt(\alpha_2)\pt(\alpha_3)\big)=\partial_{i_3}\alpha_1^{i_1}\partial_{i_1}\alpha_2^{i_2}\partial_{i_2}\alpha_3^{i_3}
-\partial_{i_2}\alpha_1^{i_1}\partial_{i_1}\alpha_3^{i_3}\partial_{i_3}\alpha_2^{i_2},
$$
si tous les $\alpha_j$ sont des champs de vecteurs, 0 sinon. (On a utilis\'e la notation $(\alpha)$ pour repr\'esenter un paquet contenant le seul tenseur $\alpha$.)

L'application $f=f^3_{111}$ n'est pas nulle car $f^3_{111}\big((x_1\partial_{2})\pt(x_2\partial_{3})\pt(x_3\partial_{1})\big)=1$. Elle est bien d\'efinie sur $G(\mathcal G)$. De plus, elle est un cocycle car $d_m(f)\in C^4_{211}$ mais $f(m(\alpha\underline\otimes\beta),\alpha_2,\alpha_3)$ est non nul seulement si $\alpha\wedge\beta$ est un champ de vecteurs, c'est \`a dire si $\alpha$ est une constante et $\beta$ est un champ de vecteurs ou $\beta$ est une constante et $\alpha$ est un champ de vecteurs.

$m'\big(f_1^1(\alpha),f(\beta,\alpha_2,\alpha_3)\big)$ est non nul seulement si $\alpha$ est une constante et $\beta$ est un champ de vecteurs.

De m\^eme, $m'\big(f(\alpha,\alpha_2,\alpha_3),f_1^1(\beta)\big)$ est non nul seulement si $\beta$ est une constante et $\alpha$ est un champ de vecteurs.

Il nous reste $$(d_mf)\big((\alpha),(\beta),(\alpha_2),(\alpha_3)\big)=-f^3_{111}(\alpha\beta,\alpha_2,\alpha_3)+\alpha f^3_{111}(\beta,\alpha_2,\alpha_3)=0$$
ou
$$(d_mf)\big((\alpha),(\beta),(\alpha_2),(\alpha_3)\big)=\beta f^3_{111}(\alpha,\alpha_2,\alpha_3)-f^3_{111}(\alpha\beta,\alpha_2,\alpha_3)=0.$$

D'autre part, $d_\ell(f)=0$. En effet, on a $\ell'=0$ et n\'ecessairement $d_\ell(f)\in C^4_{1111}$.

Alors, les seuls termes restant sont de la forme $f^3_{111}\big(\ell((\alpha)\pt(\beta))\pt(\gamma)\pt(\delta)\big)$. Ces termes sont diff\'erents de $0$ seulement si $\alpha$, $\beta$, $\gamma$ et $\delta$ sont des champs de vecteurs. Il ne reste que:
$$
\begin{aligned}&(d_\ell f)^4_{1111}\big((\alpha_1)\pt(\alpha_2)\pt(\alpha_3)\pt(\alpha_4)\big)
=f([\alpha_1,\alpha_2]_S,\alpha_3,\alpha_4)-f([\alpha_1,\alpha_3]_S,\alpha_2,\alpha_4)\cr
&+f([\alpha_1,\alpha_4]_S,\alpha_2,\alpha_3)+f([\alpha_2,\alpha_3]_S,\alpha_1,\alpha_4)
-f([\alpha_2,\alpha_4]_S,\alpha_1,\alpha_3)+f([\alpha_3,\alpha_4]_S,\alpha_1,\alpha_2)\cr
&=(d_Cf)(\alpha_1,\alpha_2,\alpha_3,\alpha_4)=0.
\end{aligned}
$$

De plus, $f$ ne peut pas \^etre un cobord car la seule possibilit\'e serait $f^3_{111}=d_\ell(g^2_{11})$ avec $d_m(g^2_{11})=0$. Mais, puisque $f^3_{111}$ s'annule sur les tenseurs qui ne sont pas des champs de vecteurs et puisque $\ell((\alpha)\pt(\beta))=-[\alpha,\beta]_S$ n'est un champ de vecteurs que si $\alpha$ et $\beta$ sont des champs de vecteurs, on peut supposer que $g^2_{11}$ s'annule sur tous les paquets $(\alpha), (\beta)$ sauf si $\alpha$ et $\beta$ sont des champs de vecteurs.

Alors, l'\'equation $f^3_{111}=d_\ell(g^2_{11})$ s'\'ecrit $f=d_C(g^2_{11})$ et on sait que $f$ n'est pas un cobord pour la cohomologie de Chevalley.

\begin{thm}
\

La cohomologie de Chevalley-Harrison de $T_{poly}^{hom}(\mathbb{R}^d)$ \`a valeurs dans $\mathbb{R}$ n'est pas triviale.
\end{thm}

\begin{rema}
\

Ce cocycle est du type des op\'erateurs d\'efinis par des ``graphes avec paquets" introduits par Gammella et Halbout (voir \cite{[GaHa]}). Il est associ\'e au graphe:

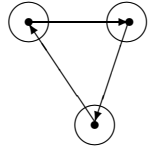
\begin{center}
\begin{picture}(40,60)(10,10)
\put(43,59){\vector(-1,-3){13}}\put(43,58){\circle*{3}}\put(42,58){\circle{15}}
\put(30,20){\vector(-2,3){25}}\put(30,19){\circle*{3}}\put(30,19){\circle{15}}
\put(5,58){\vector(1,0){38}}\put(5,58){\circle*{3}}\put(5,58){\circle{15}}
\end{picture}
\end{center}

Le fait de se restreindre \`a $T_{poly}^{hom}(\mathbb{R}^d)$ nous a permis d'\'eliminer dans le calcul de $d_m(f^3_{111})$ le graphe suivant:

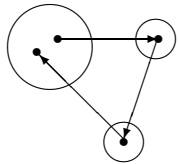
\begin{center}
\begin{picture}(40,60)(10,10)
\put(43,59){\vector(-1,-3){13}}\put(43,58){\circle*{3}}\put(42,58){\circle{15}}
\put(30,20){\vector(-1,1){32}}\put(30,19){\circle*{3}}\put(30,19){\circle{15}}
\put(5,58){\vector(1,0){38}}\put(-3,53){\circle*{3}}\put(5,58){\circle*{3}}\put(2,55){\circle{30}}
\end{picture}
\end{center}

\end{rema}

\



\begin{thebibliography}{12}



\vskip 0.8cm
\bibitem[AAC1] {[AAC1]} W. Aloulou, D. Arnal, R. Chatbouri, {\sl Cohomologie de Chevalley des graphes vectoriels}, Pacific J of Math, vol 229, no
2, (2007) 257-292.\\

\bibitem[AAC2] {[AAC2]} W. Aloulou, D. Arnal, R. Chatbouri, {\sl Chevalley cohomlogy for linear graphs}, Lett. Math. Phys. vol 80 (2007) 139-154.\\

\bibitem[AMM] {[AMM]} D. Arnal, D. Manchon, M. Masmoudi, {\sl Choix des signes pour la formalit\'e de M. Kontsevich}, Pacific J of Math, vol 203, no 1 (2002),
23-66.\\

\bibitem[BGHHW] {[BGHHW]} M. Bordemann, G. Ginot, G. Halbout, H.C. Herbig, S. Waldmann, {\sl Formalit\'e $G_{\infty}$ adapt\'ee et star-repr\'esentations sur des sous vari\'et\'es co\"{\i}sotropes}, math.QA/0504276 v 1 13 Apr 2005.\\

\bibitem[G] {[G]} G. Ginot, {\sl Homologie et mod\`ele minimal des alg\`ebres de Gerstenhaber}, Ann. Math. Blaise Pascal, vol 11, no 1 (2004),
95-126.\\

\bibitem[GaHa] {[GaHa]} A. Gamella, G. Halbout, {\sl $G_\infty$-Formality Theorem in Terms of Graphs and Associated Chevalley-Eilenberg-Harrison Cohomology }, Comm. in algebra Vol 33. N$^\circ$ 10, sep 2005, 3515-3528.\\

\bibitem[GH] {[GH]} G. Ginot, G. Halbout, {\sl A formality theorem for Poisson manifolds}, Lett. Math. Phys. Vol 66 (2003), 37-64.\\

\bibitem[K] {[K]} M. Kontsevich, {\sl Deformation quantization of Poisson manifolds}, Lett. Math. Phys. Vol 66 (2003), no. 3, 157-216.\\

\bibitem[L] {[L]} J.L. Loday, {\sl Cyclic Homology}, Second Edition, Volume 301,
Grundlerhren der Mathematischern Wissenschaften, A series of compehensive studies in mathematics, Springer.\\

\bibitem[T] {[T]} D. Tamarkin, {\sl Another proof of M. Kontsevich's formality theorem}, math.QA/9803025.\\

\end{thebibliography}
\end{document}